# Stochastic Hartree NLS in 3d coming from a Many-Body Quantum System with White Noise Potential

by Francesco Carlo De Vecchi[a], Xiaohao Ji[b], Immanuel Zachhuber[b]

*a*. Department of Mathematics,
University of Pavia, Italy

*b*. Department of Mathematics,
Freie Universität Berlin, Germany

## Abstract

In this paper, we consider the defocusing Hartree NLS with white noise external potential on $\mathbb{T}^3$ i.e. the Hartree NLS whose linear part is given by the Anderson Hamiltonian. A Strichartz-type estimate is established for the Anderson Hamiltonian using perturbative arguments and the local and global well-posedness of the NLS is considered with initial data in the domain and form-domain of the Anderson Hamiltonian under different regularity assumptions on the Hartree interaction.

Furthermore, we establish the Anderson Hartree NLS as an effective equation describing many-body Bosonic systems and, in particular, we prove the convergence of the linear Schrödinger equation for the many body system to the BBGKY hierarchy for the Coulomb interaction.

**A.M.S. subject classification:** 60H17, 81V70, 35Q55

## Table of contents

# 1 Introduction

In this paper we study a *non-linear Schrödinger equation* where the linear part is given by the *Anderson Hamiltonian* operator and the non-linearity is of *Hartree-Fock* type. More precisely, we consider the following formal (stochastic) PDE

$$i\partial_t u(t,x) = -\frac{1}{2}(\Delta + \xi(x) - \infty)u(t,x) + u(t,x)V * |u(t,\cdot)|^2(x), \quad (t,x) \in \mathbb{R}_+ \times \mathbb{T}^3 \quad (1.1)$$

where $\xi$ is *white noise* on $\mathbb{T}^3$, the symbol $-\infty$ indicates a suitable *renormalization* procedure, and $V \geqslant 0$ is an interaction potential on $\mathbb{T}^3$, whose regularity and integrability will be discussed below.

The *continuum Anderson Hamiltonian* can be thought of describing the motion of particles in a random environment. The analysis of such random operators is also motivated by the phenomenon called *Anderson localization* after the seminal work [And58], which describes the absence of diffusion of waves in a disordered medium. Anderson showed that in 3D the bottom of the spectrum of the discrete Anderson Hamiltonian consists of localized eigenfunctions. Mathematically, the (de)localization of the Anderson Hamiltonian is more well-understood in the discrete[1.1] setting, c.f. [GMP77], [KS80], [FS83], [AM93], [BK01], [BK05], [GKM07] and [CGK10], and the nice books and lecture notes [Kir07], [Sto11], [CL12] and [Kön16] and references therein for a more complete overview of the vast field of research. Meanwhile, the rigorous justification of the Anderson (de)localization with the continuum operator has been well-studied only in 1D with the recent breakthroughs [DL20] and [DL24].

Using the theory of Dirichlet forms, Fukushima and Nakao [FN77] constructed the continuum Anderson Hamiltonian on $[-L,L]$ with Dirichlet boundary condition for $0 < L < \infty$, while the construction is more challenging for dimensions two and three. This was only resolved after the introduction of the powerful tools of singular SPDEs, like *regularity structures* [Hai14] and *paracontrolled calculus* [GIP15]. The (pathwise) solution theory of singular SPDEs has witnessed tremendous advances in recently years, and we refer the readers to [Kup16], [CHS18], [Duc21], [LOTT24], [BC24], [CCHS24] and references therein for discussions on various equations and different perspectives. Roughly speaking, those methods allow to solve via fixed-point *scaling subcritical* singular SPDEs in a suitably constructed space by treating the nonlinearity of the equation perturbatively in some sense.

---

[1.1]. which we mean to be either the discretization on lattices $\mathbb{Z}^d$ or the Schrödinger operators on the continuum spaces with discretized random potentials.

While *parabolic* singular SPDEs are quite well understood, less is known outside that regime. There have been some interesting developments on singular dispersive SPDEs, an interesting line of inquiry that largely started with [GKO23] and [GKOT22], where nonlinear wave equations with additive space-time white noise are considered. See also [OOR20] and [ORSW21] amongst others in a similar vein.

The first construction of the (singular) continuum Anderson Hamiltonian was in [AC15] on $\mathbb{T}^2$ using paracontrolled calculus. This construction was extended and reformulated in [GUZ20] also to $\mathbb{T}^3$ using paracontrolled calculus combined with an exponential ansatz inspired by [HL15] and [DW18]. Labbé [Lab19] used regularity structures to construct the operator in dimensions up to three with periodic and Dirichlet boundary conditions. Let us also mention the following related constructions and results on bounded domains/compact manifolds [CvZ21][MvZ22], [BDM22], [Mou22],[MO25] and on the full space [Uek25],[HL24].

In a related but separate development, dispersive equations with random data have been studied widely, mainly as part of studying invariant (Gibbs) measures following the pioneering works [Fri85], [Zhi91],[LRS88] and [Bou94] in dimension one and [Bou96] in dimension two. More recently, in [DNY24], [DNY22] and [BDNY24] the new theory of *random tensors* was introduced and implemented to vastly extend Bourgain's results on 2d NLS with general power nonlinearity and to 3d cubic nonlinear wave.

Schrödinger equations of the form (1.1) with a Hartree-Fock nonlinearity (without the noise term $\xi$) physically describe the so-called *Bose-Einstein* condensation in either the mean-field limit, dilute gas or the Gross-Pitaevskii limits. The Schrödinger equation combining the Anderson Hamiltonian and a Hartree-Fock non-linearity can then be used to model the behaviour of a Bose-Einstein condensate in a random environment. Such a model has been tested in some experiments having the main aim of showing and studying the Anderson localization in a Bose-Einstein condensate [RDF+08]. In this paper we rigorously study such a limit jointly with a limit in the noise regularization as in order to justify the study of the Anderson Hartree NLS (1.1).

The convergence of linear many-body Schrödinger equation to the Hartree-Fock nonlinear Schrödinger equation *without the presence of a random external potential* is investigated extensively both in the physics and mathematics liteature and we do not claim to give an exhaustive overview here, see for example the nice survey [Rou21]. The first proof of qualitative convergence of many-body Bosonic system was obtained by Spohn [Spo80] under the assumption that the potential is bounded under the framework of the *BBGKY (Bogoliubov–Born–Green–Kirkwood–Yvon) hierarchy*, which was then further extended to Coulomb interactions in the seminal work of Erdős and Yau [EY01], c.f. also [BGM00]. The Gross-Pitaevskii limit was extensively studied in a series of important works of Elgart-Erdős-Schlein-Yau [EESY06], [ESY06], [ESY07], [ESY09], [ESY10], and the uniqueness of the BBGKY hierarchy was also considered [KM08], [CP14], [CH16] and many others. While those works above exploit the BBGKY hierarchy, a quantitative convergence can be derived by introducing the Bosonic Fock space following, see [Hep74]. Various methods have beed developed ever since, including the technique of coherent states [GV79a], [GV79b]; Pickl's method of excited particles [KP10], [Pic11], [Pic15]; Egorov's method [FGS07], [FKP07], see also [ADU17, ADRU20, ADRU22] where a probabilistic approach using Nelson's stochastic mechanics and sto-chastic optimal control was used to prove the convergence of the $N$-particle system to the Bose-Einstein condensate.

The study of well-posedness results of nonlinear Schrödinger equations with multiplicative spatial white noise, which is studied in the current paper, was initiated by [DW18] who studied the cubic equation on $\mathbb{T}^2$ and was extended to the sub-cubic case on $\mathbb{R}^2$ in [DM19]. In a similar vein, i.e. using the same kind of exponential transform (inspired by [HL15]), [TV23b] and [TV23a] extended the results of [DW18] on $\mathbb{T}^2$ to arbitrary power nonlinearities using Strichartz estimates and modified energies introduced in [PTV17] and in [DLTV24] the analogous result was proved on $\mathbb{R}^2$. The same kind of equation was studied in [GUZ20] in slightly different spaces using a paracontrolled description of the Anderson Hamiltonian on $\mathbb{T}^2$ and $\mathbb{T}^3$ obtaining global well-posedness of the cubic equation on $\mathbb{T}^2$ and local well-posedness on $\mathbb{T}^3$– the latter result we recall in Proposition 4.1– as well as global existence of energy solutions on $\mathbb{T}^2$. These results were extended to compact surfaces[Mou22],[MZ24] where in the latter also Strichartz estimates were proved and local well-posedness was proved in low-regularity spaces independent of the noise. Also using Strichartz estimates and modified energies, global well-posedness in the energy space and in the domain of the Anderson Hamiltonian was shown in [Zac19] for general nonlinearities. In a very similar way, the cubic nonlinear wave equation with multiplicative white noise was shown to be globally well-posed in [GUZ20] on $\mathbb{T}^2$ and $\mathbb{T}^3$ in both the strong and energy settings. Using finite speed of propagation, this equation was studied on the full space in [Zac21] and its invariant Gibbs measure on $\mathbb{T}^2$ was constructed in [BDZ23] combining the singular SPDE ideas of [GUZ20] with the methods of [OT20] to show invariance of Gibbs measures for wave equations. We also refer to the papers [CG15], [XYW24] and the references therein for similar results for (Hartree) Schrödinger equations with multiplicative noise which is *white in time* but more regular in space.

The multiplicative stochastic nonlinear Schrödinger equations in the (physical) dimension three has not been widely studied in the literature, where to the best of the authors' knowledge, the local well-posedness result in [GUZ20] is the only result in this direction. Broadly speaking, there are several technical difficulties compared to the more well-studied two dimensional setting.

- The noise is $\frac{1}{2}$ less regular, meaning that one needs a more involved ansatz and one has more irregular error terms. In particular, if one uses an exponential transform as in [DW18] one has formally $(-\Delta + \xi)e^X : H_{\mathbb{T}^2}^2 \to H_{\mathbb{T}^2}^{-\varepsilon}$ meaning one is $\varepsilon > 0$ from being in $L_{\mathbb{T}^2}^2$ but in three dimensions one has no better than (formally) $(-\Delta + \xi)e^W : H_{\mathbb{T}^3}^2 \to H_{\mathbb{T}^3}^{-\frac{1}{2}-\varepsilon}$, see Section 2 for the correct choice of $W$, so if one naively inserts such an ansatz into the equation, one gets that $\partial_t u \in H_{\mathbb{T}^3}^{-\frac{1}{2}-\varepsilon}$ while $e^{-X} u$ needs to be $H_{\mathbb{T}^3}^2$ so it is rather difficult to bridge this gap.

- One has worse Sobolev/Besov embeddings in 3d, for example one has $H_{\mathbb{T}^2}^1 \hookrightarrow L_{\mathbb{T}^2}^p$ for all $1 \leq p < \infty$ but only $H_{\mathbb{T}^3}^1 \hookrightarrow L_{\mathbb{T}^3}^p$ for $1 \leq p \leq 6$, so treating higher powers in the nonlinearity becomes more difficult. Relatedly, in two dimensions one "almost" has that the $H_{\mathbb{T}^2}^1$ norm controls the $L_{\mathbb{T}^2}^\infty$ norm (this is true up to a logarithmic term, see [DW18] and [GUZ20] where this Brezis-Gallouet inequality was essential), but in three dimensions, one has no better than $\|u\|_{L^\infty_{\mathbb{T}^3}} \lesssim \|u\|_{H^1_{\mathbb{T}^3}}^{\frac{1}{2}}\|u\|_{H^2_{\mathbb{T}^3}}^{\frac{1}{2}}$ so one again "loses" 1/2 derivative compared to two dimensions.

- The Strichartz estimates on $\mathbb{T}^3$ give less integrability than the corresponding ones on $\mathbb{T}^2$, namely $\|e^{-it\Delta}u\|_{L^4_{t;[0,1]}L^4_{\mathbb{T}^2}} \lesssim \|u\|_{H^\varepsilon_{\mathbb{T}^2}}$ and $\|e^{-it\Delta}u\|_{L^{\frac{10}{3}}_{t;[0,1]}L^{\frac{10}{3}}_{\mathbb{T}^3}} \lesssim \|u\|_{H^\varepsilon_{\mathbb{T}^3}}$ for $\varepsilon > 0$. Moreover, since the proof of Strichartz estimates for the operator $\mathcal{H}\text{``=''}-\Delta + \xi$ from [MZ24],[Zac19] is perturbative and one loses more derivatives the more irregular the noise.

The chief accomplishment of the main part of the paper is the first global well-posedness of NLS with Anderson Hamiltonian on $\mathbb{T}^3$. Because of the above difficulties, we restrict ourselves to Hartree nonlinearities with a range of parameters and separately consider well-posedness separately in the natural energy and strong settings. Finally, in Section 6 we give a justification for the study of the type of equation, namely that it arises as the mean field limit of $N$–body Quantum Dynamics with smooth potentials converging to white noise. To the best of our knowledge, this is the first result in this direction, however let us mention the recent preprint [NZZ25] which deals with a different problem but also combines singular SPDEs and many body Quantum Dynamics and [SS15] which studies the Gross-Pitaevskii hierarchy with random initial data.

One of the main differences of the setting considered here and the ones mentioned above is that the Anderson Hamiltonian operator has a peculiar domain (see Section 2.1) due to the roughness of the random potential (which is only in the Hölder-Besov space $C^{-3/2-}_{\mathbb{T}^3}$ a.s.) which intersects with the space of smooth functions only in zero and thus there is no control on the commutator of the operator with the multiplication by the interaction potential $V$. Moreoever, the resolvent of the Anderson Hamiltonian operator does not regularise as that of Laplacian does. We summarize the two main aims of this paper as follows

I. to prove the global well-posedness of equation (1.1) in three spatial dimension for a certain family of interaction potentials $V$;

II. to prove the convergence of a suitable N-particle bosonic system to infinitely many independent bosonic particles described by a wave function solution to equation (1.1) which in turn gives a physical motivation to study this kind of equation.

## 1.1 Description of the main results

The main results of this paper are briefly summarized in this section and we refer the readers to the fully detailed and rigorous definitions and statements in the corresponding later sections. In Corollary 5.1, Theorem 5.12, Theorem 5.14 and Proposition 5.15, we establish the well-posedness results of the Anderson NLS (1.1) under different conditions. The rigorous definitions of the domain and form domain of the operator are given in Section 2. Another thing to remark is that since we follow the general methodology of [GUZ20], all of the solutions we obtain are also limits of solutions of regularised equations as was proved rigorously therein. This kind of continuity in the noise is only implicit apart from in the final section where we consider the joint limit of regularization parameter and the number of particles.

**Theorem. (Global Well-posedness)**

a) **Strong GWP:** Let $T > 0$ and $u_0$ be in the domain $\mathcal{D}(\mathcal{H})$ of the Anderson Hamiltonian, and suppose that either $0 \leq V \in L^\infty_{\mathbb{T}^3}$ or
$$0 \leq V \in W^{\beta,1}_{\mathbb{T}^3}$$
with $0.733 < \beta$. Then there is a unique (mild) global in time solution to equation (1.1) with initial condition $u_0$ in $C([0,T]; \mathcal{D}(\mathcal{H}))$ depending continuously on the data and whose norm grows algebraically for large $T > 0$.

b) **Energy GWP:** Take again $T > 0$ and assume $u_0$ belongs to the form-domain $\mathcal{D}(\sqrt{-\mathcal{H}})$ and
$$0 \leq V \in W^{\beta,1}_{\mathbb{T}^3} + L^\infty_{\mathbb{T}^3}$$
with $\beta > \frac{17}{20}$ then there exists a unique (mild) solution to equation (1.1) in $C([0,T]; \mathcal{D}(\sqrt{-\mathcal{H}}))$ which depends continuously on the data.

**Remark 1.1.** A particular type of potential of interest for the well-posedness results above is of the form
$$V(x) = \frac{1}{|x|^{3-\beta}},$$
around the origin such that the above theorem can be summarized in the following table, writing LWP/GWP to mean local/global well-posedness:[1.2]

|  | LWP | GWP |
| --- | --- | --- |
| $u_0 \in \mathcal{D}(\mathcal{H})$ | $\beta > 0$ | $\beta > 0.733$ |
| $u_0 \in \mathcal{D}(\sqrt{-\mathcal{H}})$ | $\beta > 0.85$ | $\beta > 0.85$ |

**Figure 1.1.** Well-posedness for (6.10)

In the second part of the paper, we focus on the convergence of the solution to $N$-particle Schrödinger equations to the tensorizations of solution to the NLS with Hartree nonlinearity.

We follow the strategy of the BBGKY hierarchy, that is we write down the Heisenberg picture of the $N$-body non-relativistic spinless bosonic system

$$i\partial_t \rho_{N,\delta}(t) = [H_{N,\delta}, \rho_{N,\delta}(t)], \quad \text{where } H_{N,\delta} := -\frac{1}{2}\sum_{j=1}^N \mathcal{H}_{\delta,x_j} + \frac{1}{N}\sum_{1 \leq j < k \leq N} V_{j,k}, \tag{1.2}$$

where $(\rho^n_{N,\delta})_{n\in\mathbb{N}}$ are the marginal distributions, $\mathcal{H}_{\delta,x_l}$ is the $\delta$-regularized Anderson Hamiltonian falling on the particle $l$ and $V_{j,k}$ is the multiplication operator of the interaction between the particles $j$ and $k$. In the case where $\rho_{N,\delta}(t)$ is a pure state, namely there is $\Psi_{N,\delta} \in L^\infty(\mathbb{R}_+; L^2_s((\mathbb{T}^3)^N))$ s.t. $\rho_{N,\delta}(t) = \Pi_{\Psi_{N,\delta}(t,\cdot)}$, the wave function $\Psi_{N,\delta}$ solves the standard Schrödinger equation associated with the Hamiltonian $H_{N,\delta}$ i.e.

$$i\partial_t \Psi_{N,\delta}(t,x_{1:N}) = -\frac{1}{2}\sum_{j=1}^N \mathcal{H}_{\delta,x_j}\Psi(t,x_{1:N}) + \frac{1}{N}\sum_{1 \leq j < k \leq N} V(x_j - x_k)\Psi_{N,\delta}(x_{1:N}),$$

---

[1.2]. We even have local well-posedness for initial data in $\mathcal{D}(\mathcal{H})$ for the Anderson NLS with power nonlinearity.

where $x_{1:N} = (x_1, \ldots, x_N) \in (\mathbb{T}^3)^N$.

Clearly $H_{N,\delta}$ is invariant under any permutation of particles. We consider the $n$-particle marginal distribution of the particle system. By assuming indistinguishability of the particles, fixing the $n$ parameter and taking $N \to \infty$, we derive formally a hierachy equation

$$i\partial_t \rho_\infty^n = -\frac{1}{2}\sum_{j=1}^n [\mathcal{H}_{x_j}, \rho_\infty^n] + \sum_{j=1}^n \text{Tr}_{n+1}([V_{j,n+1}, \rho_\infty^{n+1}]), \tag{1.3}$$

since we assume mean-field interaction in (1.2), where $(\rho_\infty^n)_{n\in\mathbb{N}}$ are the limits of the marginal distributions, $\mathcal{H}_j$ is the Anderson Hamiltonian falling on the particle $j$, $\text{Tr}_{n+1}$ is the partial trace with the particle labelled $n+1$ and the square bracket is the commutator.

As in the classical mean-field theory of BBGKY hierarchies, we can establish (1.1) as the effective equation of the Bosonic many-body system if we can prove that the formal limit we mentioned above is rigorous and (1.3) has a unique solution. The next result of this paper is thus *the rigorous justification of the limit to (1.3), when we assume Coulomb interaction* in Theorem 6.3:

**Theorem. (Convergence towards BBGKY)** *Suppose that the interaction potential V is of the form*

$$0 \leq V(x) = \frac{V_1(x)}{|x|^{3-\beta}} + V_2(x), \quad x \in \mathbb{T}^3,$$

*for $V_1 \in L^\infty_{\mathbb{T}^3}$ and $V_2 \in L^{3+}_{\mathbb{T}^3}$, where $\beta \geq 2$ and*

$$\sup_{N\in\mathbb{N}, \delta>0} \left\| \sqrt{-\mathcal{H}_{\delta,x_1}} \Psi_{N,\delta}(0) \right\|_{L^2((\mathbb{T}^3)^N)} < +\infty,$$

*then up to a subsequence, we have for any $n \in \mathbb{N}$ that*

$$\lim_{N\to\infty, \delta\to 0^+} \rho_{N,\delta}^n = \rho_\infty^n, \quad s.t. \ \{\rho_\infty^n\}_{n\in\mathbb{N}} \ solves \ (1.3).$$

Lastly, we show in Theorem 6.4 the uniqueness of solution to equation (1.3) when *V is bounded*, as a direct consequence of which, the convergence of the Bosonic system to the solution of the equation (1.1) is validated. This implies that (1.1) indeed describes effectively the many-body Bosonic system in random environment under consideration, when the number of particles is sufficiently large, which is the content of Corollary 6.5.

**Theorem. (Anderson NLS as Effective Equation)** *If we suppose that $V \in L^\infty_{\mathbb{T}^3}$, then equation (1.3) admits a unique solution. If we further assume that*

$$\sup_{N\in\mathbb{N}, \delta>0} \left\| \sqrt{-\mathcal{H}_{\delta,x_1}} \Psi_{N,\delta}(0) \right\|_{L^2_{\mathbb{T}^3}} < +\infty \quad and \quad \lim_{N\to+\infty, \delta\to 0} \text{Tr}_{n+1:N}[\Pi_{\Psi_{N,\delta}(0)}] = \Pi_{u_0^{\otimes n}}, \quad \forall n \in \mathbb{N},$$

*for some $u_0 \in \mathcal{D}(\sqrt{-\mathcal{H}})$. Then we have*

$$\lim_{N\to+\infty, \delta\to 0} \sup_{t\in\mathbb{R}_+} \|\rho_{N,\delta}^n(t,\cdot) - \Pi_{u(t,\cdot)^{\otimes n}}\|_{L^2((\mathbb{T}^3)^n)} = 0, \quad \forall n \in \mathbb{N},$$

*where $u$ is the unique solution to equation (1.1) with initial condition $u_0$.*

## 1.2 Notations

### 1.2.1 Function Spaces

We denote by $\mathbb{T}^d = (\mathbb{R}/\mathbb{Z})^d$ the $d$-dimensional torus and write $L^p_{\mathbb{T}^d} = L^p(\mathbb{T}^d)$, $p \in [1,\infty]$ for the usual Lebesgue spaces. By $H^\alpha_{\mathbb{T}^d}$ we mean the Sobolev space with index $\alpha \in \mathbb{R}$ defined as $H^\alpha_{\mathbb{T}^d} := \{u \in \mathcal{S}'_{\mathbb{T}^d} : \|(1-\Delta)^{\alpha/2}u\|_{L^2_{\mathbb{T}^3}} < \infty\}$.

Moreover, we write $W^{s,p}_{\mathbb{T}^d} := B^s_{p,p}(\mathbb{T}^d)$ which we call *fractional Sobolev spaces* and $C^s_{\mathbb{T}^d} := B^s_{\infty,\infty}(\mathbb{T}^d)$ called *Hölder-Besov spaces* for $s \in \mathbb{R}, p \in [1,\infty]$ which are special cases of Besov spaces whose definition is recalled in Appendix A. Occasionally we will write, by a slight abuse, $W^{0,p}_{\mathbb{T}^d} = L^p_{\mathbb{T}^d}$ and $W^{s,\infty}_{\mathbb{T}^d} = C^s_{\mathbb{T}^d}$.

We also sometimes use the following notation $\|f(t)\|_{L^p_{t;[a,b]}} := \|f\|_{L^p_{[a,b]}} := \left(\int_a^b |f(t)|^p dt\right)^{\frac{1}{p}}$ for an interval $-\infty \leq a < b \leq \infty$ and the obvious modification for $p = \infty$.

### 1.2.2 Density Matrices

Let $\mathcal{L}^1_n = \mathcal{L}^1(L^2_s((\mathbb{T}^3)^n))$ be the set of trace-class operators on the symmetric space $L^2_s((\mathbb{T}^3)^n)$, the set of symmetric $L^2$ functions on $n$-copies of the three dimensional torus $\mathbb{T}^3$ for $n \in \mathbb{N}$. The class of *density matrices* is then the subset of $\rho \in \mathcal{L}^1_n$ such that $\rho$ is self-adjoint, non-negative and normalized by $\text{Tr}(\rho) = 1$. If $S$ is a linear operator on $L^2(\mathbb{T}^3)$, we define $S_{x_k} = S_k$ to be the linear operator on the $k$-th component of $L^2_s((\mathbb{T}^3)^n)$, where the tensor product is symmetrized. We also denote the family of compact operators on $L^2_s((\mathbb{T}^3)^n)$ as $\mathcal{K}_n = \mathcal{K}(L^2_s((\mathbb{T}^3)^n))$.

More generally, we denote by $\mathcal{L}^p(L^2_s((\mathbb{T}^3)^n))$ the space of symmetric compact operators with finite *$p$-Schatten norm*, namely

$$A \in \mathcal{L}^p(L^2_s((\mathbb{T}^3)^n)) \iff \|A\|^p_{\mathcal{L}^p} = \text{Tr}(|A|^p) < +\infty.$$

If $A \in \mathcal{L}^p(L^2_s((\mathbb{T}^3)^n))$ for some $p \leq 2$, $A$ is in particular a Hilbert-Schmidt operator and thus, we can –and will– identify $A$ with its kernel in $L^2((\mathbb{T}^3)^{2n})$. We further define the partial trace $\text{Tr}_{k+1:n}$ as a bounded operator from $\mathcal{L}^1(L^2_s((\mathbb{T}^3)^n))$ into $\mathcal{L}^1(L^2_s((\mathbb{T}^3))^k)$ for $k \leq n$ to be such that

$$\text{Tr}(\text{Tr}_{k+1:n} A K) = \text{Tr}(AK \otimes \text{Id}_{k+1:n}), \tag{1.4}$$

for any $K \in \mathcal{K}(L^2_s((\mathbb{T}^3)^{n-k}))$. Equivalently, the partial trace is defined in terms of kernels by

$$\text{Tr}_{k+1:n}(A)(x_{1:k}, x'_{1:k}) = \int_{\mathbb{T}^{3(n-k)}} A((x_{1:k}, z_{k+1:n}), (x'_{1:k}, z_{k+1:n})) dz_{k+1:n}$$

where, hereafter, $w_{h:\ell}$ denotes the vector $w_{h:\ell} \in (\mathbb{T}^3)^{(\ell-h+1)}$ formed by the components from $h$-th to $\ell$-th of an arbitrary vector $w_{1:n} \in (\mathbb{T}^3)^n$. We call a family of trace-class operators $(\rho^n)_{n \in \mathbb{N}}$ with $\rho^n \in \mathcal{L}^1(L^2_s(\mathbb{T}^{3n}))$ *compatible*, if for any $n \in \mathbb{N}$ and $k \leq n$ we have

$$\rho^k = \text{Tr}_{k+1:n}(\rho^n).$$

We use the same name for a finite set $(\rho_n)_{n \leq N}$ satisfying the equation as above. We further define the Sobolev spaces of operators with respect to the self-adjoint, negative and densely defined $\mathcal{H}$ and $\mathcal{H}_\delta$ to be

$$\mathcal{W}_N^{\alpha,\beta} = \{\rho \in \mathcal{L}^\beta((\mathbb{T}^3)^N) : \|\rho\|_{\mathcal{W}_N^{\alpha,\beta}} < +\infty\}, \quad \|\rho\|_{\mathcal{W}_N^{\alpha,\beta}} = \left(\sum_{k=1}^N \mathrm{tr}(|(-\mathcal{H}_{x_k})^{\alpha/2}\rho(-\mathcal{H}_{x_k})^{\alpha/2}|^\beta)\right)^{\frac{1}{\beta}} \quad (1.5)$$

and respectively

$$\mathcal{W}_{N,\delta}^{\alpha,\beta} = \{\rho \in \mathcal{L}^\beta((\mathbb{T}^3)^N) : \|\rho\|_{\mathcal{W}_{N,\delta}^{\alpha,\beta}} < +\infty\}, \quad \|\rho\|_{\mathcal{W}_{N,\delta}^{\alpha,\beta}} = \left(\sum_{k=1}^N \mathrm{tr}(|(-\mathcal{H}_{\delta,x_k})^{\alpha/2}\rho(-\mathcal{H}_{\delta,x_k})^{\alpha/2}|^\beta)\right)^{\frac{1}{\beta}}, \quad (1.6)$$

where the Anderson Hamiltonian $\mathcal{H}$ is given in Proposition 2.3, and its mollified version $\mathcal{H}_\delta$ is defined in (6.1).

## 1.3 Other notational conventions

We use some usual conventions, such as writing

$$a \lesssim b$$

if there exists a constant $C > 0$ s.t. $a \leq Cb$ which does not depend on relevant parameters. We also write

$$a \asymp b \text{ if } a \lesssim b \text{ and } b \lesssim a$$

and if we want to keep track that the implicit constant depends on a parameter $\gamma$ we write $a \lesssim_\gamma b$ etc.

Since the problem we are considering is *path-wise*, meaning that we fix a realization $\omega \in \Omega$ of the noise which has an enhanced noise $\Xi(\omega)$ with $\|\Xi\|_\mathfrak{X} < \infty$, see Lemma 2.1. Since all our constructions depend on this realisation so we allow *every* implicit constant to depend on this norm, schematically written as

$$a \lesssim b \Longleftrightarrow a \lesssim_\Xi b.$$

On a related note, we often write things like $f \in C_{\mathbb{T}^d}^{\alpha-}$ to mean $f \in C_{\mathbb{T}^d}^\beta$ for all $\beta < \alpha$ and in some cases we take it to mean $f \in C_{\mathbb{T}^d}^{\alpha-\kappa}$ for an arbitrarily small $\kappa > 0$ which can be chosen independently of other parameters. Similarly, we write something of the form

$$\|Af\|_X \lesssim \|f\|_{H_{\mathbb{T}^d}^{\alpha+}} \Longleftrightarrow \|Af\|_X \lesssim \|f\|_{H_{\mathbb{T}^d}^{\alpha+\varepsilon}} \text{ for all } \varepsilon > 0$$

to mean that the bound is true for every $\varepsilon > 0$ and the implicit constant blows up as $\varepsilon \downarrow 0$.

## 1.4 Flattened, sharpened and natural variables

We try to stay consistent in labelling the variables by the "flat", "sharp", and "natural" superscript to denote the transformation one is considering which are succesively more refined and complicated. For the convenience of the reader, we briefly summarize their definitions and where they are used.

- **Flat transform**: multiplication by an exponential of a noise term, see (2.9), which characterizes the form domain of $\mathcal{H}$ i.e. $u^\flat = e^{-W} u$ and $\|u\|_{\mathcal{D}(\sqrt{-\mathcal{H}})} \asymp \|u^\flat\|_{H_{\mathbb{T}^3}^1}$.

- **Sharp transform**: paracontrolled ansatz (inverse of $\Gamma$, see Lemma 2.5) performed after the exponential ansatz i.e. $u^\sharp = \Gamma^{-1} e^{-W} u$ which characterizes the domain of $\mathcal{H}$ which we also transform as $\mathcal{H}^\sharp := \Gamma^{-1} e^{-W} \mathcal{H} e^W \Gamma$ via

$$\|u\|_{\mathcal{D}(\mathcal{H})} \asymp \|u^\sharp\|_{H^2_{\mathbb{T}^3}} \asymp \|\mathcal{H}^\sharp u^\sharp\|_{L^2_{\mathbb{T}^3}}.$$

We use this transformation in Section 4.2 and to obtain the embeddings and norm equivalences in Lemma 2.6.

- **Natural transform:** The $\Theta$ transform, see Lemma 2.8, applied after the sharp transform i.e. $u^\natural = \Theta u^\sharp = \Theta \Gamma^{-1} e^{-W} u$ which is the inverse of $\Lambda$ from Lemma 2.9. We set also $\mathcal{H}^\natural := \Lambda^{-1} \mathcal{H} \Lambda$ and this transformation has the important properties

$$\mathcal{H}^\natural - \Delta: H^{\frac{3}{2}+}_{\mathbb{T}^3} \to H^{0+}_{\mathbb{T}^3} \quad \text{and} \quad \Lambda - \mathrm{id}: H^{0-}_{\mathbb{T}^3} \to H^{\frac{1}{2}-}_{\mathbb{T}^3},$$

which are properties that are not true without this additional refinement. We heavily use this transform in Section 5.2 and this transform is one of the novelties of the current paper, whereas the other two already appeared in [GUZ20].

## 2 Anderson Hamiltonian on the 3D Torus

### 2.1 Construction

For the reader's convenience, we review briefly the construction and basic properties of the Anderson Hamiltonian on the 3D torus. We start, using largely the same notation as in [GUZ20], with a *spatial white noise* on $\mathbb{T}^3$, that is $\xi \in \mathcal{C}^{-3/2-\varepsilon}(\mathbb{T}^3)$, a.s. for $\varepsilon > 0$ arbitrarily small. We set $\xi_\delta := a_\delta * \xi$ with mollifiers $(a_\delta)_{\delta > 0}$ and the *mollified enhancement* as

$$X_\delta := (1-\Delta)^{-1} \xi_\delta, \quad X_\delta^{\vee} := (1-\Delta)^{-1}(|\nabla X_\delta|^2 - c_\delta^1),$$
$$X_\delta^{\vee\vee} := 2(1-\Delta)^{-1}(\nabla X_\delta \cdot \nabla X_\delta^{\vee}), \quad W_\delta := X_\delta + X_\delta^{\vee} + X_\delta^{\vee\vee}, \tag{2.1}$$

where $(c_\delta^1)_{\delta > 0} \sim 1/\delta$ is family of *renormalization* constants which tends to $+\infty$ as $\delta \to 0^+$ and which should counteract the divergence of the stochastic objects in the limit. We make the following exponential ansatz inspired by [HL15]:

$$u_\delta = \exp(W_\delta^M) u_\delta^\flat, \tag{2.2}$$

where $W_\delta^M := P_{>M} W_\delta$ is the high-frequency cut-off of $W_\delta$ defined in (2.1) by Fourier multiplier $P_{>M} = \mathbf{1}_{|\cdot|>2^M}(D)$. It follows then by construction that

$$\mathcal{H}_\delta u_\delta = (\Delta + \xi_\delta - c_\delta) u_\delta = \exp(W_\delta^M)(\Delta u_\delta^\flat + 2\nabla W_\delta^M \cdot \nabla u_\delta^\flat + [(1-\Delta) Z_\delta] u_\delta^\flat), \tag{2.3}$$

where $c_\delta := c_\delta^1 + c_\delta^2$ for another family of *renormalization* constant $(c_\delta^2)_{\delta > 0} \sim \ln(1/\delta)$ that also tends to $+\infty$ as $\delta \to 0^+$, and $Z$ depends on $M$ as well for now.

We emphasize that the $M$ parameter is *not* present in [GUZ20] and it is needed in order to refine the mixed exponential-paracontrolled ansatz. This is the only place where we differ from the [GUZ20] construction of the Anderson Hamiltonian. We specify the following random objects as

$$\begin{aligned}(1-\Delta)Z_\delta^M =\ & |\nabla X_\delta^{\vee}|^2 - c_\delta^2 + |\nabla X_\delta^{\vee\vee}|^2 + 2\nabla X_\delta \cdot \nabla X_\delta^{\vee\vee} + 2\nabla X_\delta^{\vee} \cdot \nabla X_\delta^{\vee\vee} + W_\delta \\ & -2\nabla(P_{\leq M} W_\delta) \cdot \nabla W_\delta + \Delta[\exp(-P_{\leq M} W_\delta)] \exp(P_{\leq M} W_\delta),\end{aligned}$$

and
$$R_\delta^1 = \nabla W_\delta \odot (1-\Delta)^{-1} \nabla^2 W_\delta, \quad R_\delta^2 = \nabla Z_\delta \odot \nabla W_\delta. \tag{2.4}$$

We will also call the family of stochastic objects
$$\Xi_\delta(\omega) := \left(X_\delta(\omega), X_\delta^{\vee}(\omega), X_\delta^{\vee}(\omega), Z_\delta(\omega), R_\delta^1(\omega), R_\delta^2(\omega)\right) \tag{2.5}$$

the mollified enhancement and refer the readers to [GUZ20] for the detailed construction of Anderson Hamiltonian on $\mathbb{T}^3$ and [CC18] for the construction of the enhancement and the proof of the convergence of the enhancement established in the space
$$\mathfrak{X}_\varepsilon := C_{\mathbb{T}^3}^{\frac{1}{2}-\varepsilon} \times C_{\mathbb{T}^3}^{1-\varepsilon} \times C_{\mathbb{T}^3}^{\frac{3}{2}-\varepsilon} \times C_{\mathbb{T}^3}^{\frac{3}{2}-\varepsilon} \times C_{\mathbb{T}^3}^{-\varepsilon}(\mathbb{R}^3) \times C_{\mathbb{T}^3}^{-\varepsilon}, \tag{2.6}$$

equipped with the product topology. Note that the formulation is slightly different from the definition of the enhanced noise in [GUZ20] but is essentially equivalent, since we have only subtracted the smooth objects $P_{\leq M} W_\delta$.

**Lemma 2.1.** *There exists a choice of $(c_\delta^1)_{\delta>0}$ and $(c_\delta^2)_{\delta>0}$ such that for any $\varepsilon > 0$, $\omega$-a.s. it holds that up to a subsequence the mollified enhancement in (2.5) converge in $\mathfrak{X}_\varepsilon$ as $\delta \to 0^+$, whose limit we denote as*
$$\Xi_\omega := (X(\omega), X^{\vee}(\omega), X^{\vee}(\omega), Z(\omega), R^1(\omega), R^2(\omega)). \tag{2.7}$$

**Notation 2.2.** *Since the value of $\varepsilon > 0$ in (2.6) is irrelevant, we will usually drop it from our notation. Furthermore, we drop the fixed realization $\omega$ so we will write for brevity*
$$\|\Xi\|_{\mathfrak{X}} := \|\Xi(\omega)\|_{\mathfrak{X}_\varepsilon}$$
*and so on.*

With all the probabilistic results at hand, it is crucially observed in [GIP15] that, by imposing the paracontrolled Ansatz[2.1], i.e.
$$u^\flat = P_{>N}[u^\flat \oslash (1-\Delta)^{-1} Z + 2\nabla u^\flat \oslash (1-\Delta)^{-1} \nabla W^M + B_\Xi(u^\flat)] + u^\sharp, \quad u^\sharp \in H_{\mathbb{T}^3}^2, \tag{2.8}$$

for some $N \in \mathbb{N}$ large enough and some continuous and linear random function $B_\Xi : H_{\mathbb{T}^3}^{3/2-} \to H_{\mathbb{T}^3}^2$, whose explicit form we postpone to (B.1), and
$$u = \exp(W^M) u^\flat, \tag{2.9}$$

one can exploit the smoothing effect of certain well-chosen commutators and then rigorously define the limit of $(\Delta + \xi_\delta - c_\delta)_{\delta>0}$ in a suitable sense.

**Proposition 2.3.** ([GUZ20]) *Given a realization $\Xi$ of the enhancement as in Lemma 2.1, we define the domain of the Anderson Hamiltonian $\mathcal{D}(\mathcal{H})$ to be*
$$\mathcal{D}(\mathcal{H}) = \exp(W^M) \mathcal{D}^\flat, \quad \mathcal{D}^\flat = \{u^\flat \in L_{\mathbb{T}^3}^2 : u^\flat \text{ satisfies (2.8) with some } u^\sharp \in H_{\mathbb{T}^3}^2\}, \tag{2.10}$$

*and the renormalized Anderson Hamiltonian is defined on $\mathcal{D}(\mathcal{H}) \ni u$ by*
$$\mathcal{H} u = \exp(W^M)(\Delta u^\sharp + [(1-\Delta)Z] \odot u^\flat + 2\nabla W^M \odot \nabla u^\sharp + G(u^\flat)) - K_\Xi u, \tag{2.11}$$

---

2.1. See Appendix A for the definition of Bony's paraproduct $\oslash$ and resonant product $\odot$.

where the linear functions $B_\Xi(\cdot)$ and $G(\cdot)$ are specified in Appendix B and $K_\Xi < \infty$ is defined to be large enough only depending on $\|\Xi\|_\mathfrak{x}$, such that

$$(u, -\mathcal{H}u) \geq \|u\|_{L^2_{\mathbb{T}^3}}^2. \tag{2.12}$$

It holds that $(\mathcal{H}, \mathcal{D}(\mathcal{H}))$ is a densely defined, closed and self-adjoint operator on $L^2_{\mathbb{T}^3}$ with discrete spectrum.

**Remark 2.4.** One can show– using paraproduct estimates–that any $u^\flat$ satisfying (2.8) is automatically in $H^{3/2-}_{\mathbb{T}^3}$. Moreover, it is straightforward that the operator and its domain do not depend on the parameters $M$ and $N$, but they only appear in the parametrization of the domain. The argument can be found in Remark 2.44 in [GUZ20].

## 2.2 Perturbative Characterizations

We now address the Fourier cutoffs $P_{>N}$ in (2.8) giving rise to the so-called $\Gamma$-*map* introduced in [GUZ20], which serves as the homeomorphism between a para-controlled variable and its smooth remainder in suitable spaces, when $N \in \mathbb{N}$ is chosen to be large enough only depending on the realization of the enhanced noise $\Xi_\omega$.

**Lemma 2.5. ($\Gamma$-map, [GUZ20])** *For a fixed realisation $\Xi$ as in Lemma 2.1, there exists an $N = N(\|\Xi\|_\mathfrak{x})$ s.t. there is a linear map $\Gamma$ given implicitly by*

$$\Gamma u^\sharp = P_{>N}[\Gamma u^\sharp \otimes (1-\Delta)^{-1}Z + 2\nabla \Gamma u^\sharp \otimes (1-\Delta)^{-1}\nabla W^M + B_\Xi(\Gamma u^\sharp)] + u^\sharp$$

*which satisfies for any $\gamma \in \left(-\frac{1}{2}, 1\right]$, $\beta \in [0, 3/2)$ and $p \in [1, \infty]$ that*

- *$\Gamma$ is a homeomorphism from $H^2_{\mathbb{T}^3}$ to $e^{-W}\mathcal{D}(\mathcal{H})$[2.2],*
- *$\Gamma$ is a homeomorphism on $W^{\beta,p}_{\mathbb{T}^3}$,*
- *$\Gamma - \mathrm{id}$ is regularizing at low regularities, more precisely*

$$\Gamma - \mathrm{id}: W^{\gamma,p}_{\mathbb{T}^3} \to W^{\gamma+\frac{1}{2}-,p}_{\mathbb{T}^3} \text{ is bounded.} \tag{2.13}$$

The proof of Lemma 2.5 (as a slight extension of Proposition 2.46 in [GUZ20]) is included in Appendix B.1 for completeness and we need to choose $N$ to be large enough due to the implicit definition of $\Gamma$ via fixed point argument.

We now collect some facts about the approximation of the Anderson Hamiltonian $\mathcal{H}$ by $\mathcal{H}_\delta$ defined in (2.3) and refer the readers to [GUZ20] for further details.

**Remark. (Mollification)** We recall from (2.2) the ansatz for mollified enhancement (2.5) reads $u = \exp(W^M_\delta)u^\flat$, and we define the $\Gamma_\delta$-map implicitly by

$$u^\flat = \Gamma_\delta u^\sharp = P_{>N}[\Gamma_\delta u^\sharp \otimes (1-\Delta)^{-1}Z_\delta + 2\nabla \Gamma_\delta u^\sharp \otimes (1-\Delta)^{-1}\nabla W^M_\delta + B_{\Xi_\delta}(\Gamma_\delta u^\sharp)] + u^\sharp, \tag{2.14}$$

---

[2.2]. where $e^{-W}\mathcal{D}(\mathcal{H}) = \mathcal{D}^\flat$ is equipped with the inner product $\|u^\flat\|_{\mathcal{D}^\flat}^2 = \|u^\flat\|_{H^{3/2}_{\mathbb{T}^3}}^2 + \|u^\sharp\|_{H^2_{\mathbb{T}^3}}^2$, which is Definition 2.42 in [GUZ20].

or equivalently

$$\Gamma_\delta^{-1} u^\flat = u^\flat - P_{>N}[u^\flat \otimes (1-\Delta)^{-1} Z_\delta + 2\nabla u^\flat \otimes (1-\Delta)^{-1}\nabla W_\delta^M + B_{\Xi_\delta}(u^\flat)].$$

It follows by linearity that

$$(\Gamma_\delta^{-1} - \Gamma^{-1})u^\flat = P_{>N}[u^\flat \otimes (1-\Delta)^{-1}(Z - Z_\delta) + 2\nabla u^\flat \otimes (1-\Delta)^{-1}\nabla(W^M - W_\delta^M) + B_{\Xi - \Xi_\delta}(u^\flat)],$$

which implies by the regularity of (2.5) and (B.3) that for any small $\varepsilon > 0$,

$$\|\Gamma_\delta^{-1} - \Gamma^{-1}\|_{\mathcal{L}(H_{\mathbb{T}^3}^{-1/2+\varepsilon}, L_{\mathbb{T}^3}^2)} \to 0, \quad \exp(W_\delta^M) \to \exp(W^M) \text{ in } C_{\mathbb{T}^3}^{\frac{1}{2}-\varepsilon} \quad \delta \to 0^+. \tag{2.15}$$

On the other hand, we can write

$$\mathcal{H}_\delta^\sharp u^\sharp = \Gamma_\delta^{-1}(\Delta u^\sharp + [(1-\Delta)Z_\delta] \odot u^\sharp + 2\nabla W_\delta^M \odot \nabla u^\sharp + G_\delta(u^\flat)) - K_\Xi u^\sharp. \tag{2.16}$$

so that

$$\begin{aligned}(\mathcal{H}_\delta^\sharp - \mathcal{H}^\sharp) u^\sharp &= (\Gamma_\delta^{-1} - \Gamma^{-1})(\Delta u^\sharp + [(1-\Delta) Z_\delta] \odot u^\sharp + 2\nabla W_\delta^M \odot \nabla u^\sharp + G_\delta(u^\flat)) \\ &\quad + \Gamma^{-1}([(1-\Delta)(Z_\delta - Z)] \odot u^\sharp + 2\nabla(W_\delta^M - W^M) \odot \nabla u^\sharp + (G_\delta - G)(u^\flat))\end{aligned}$$

It follows from (2.15), the regularity of (2.5) and (B.4) that for any small $\varepsilon > 0$

$$\|\mathcal{H}_\delta^\sharp - \mathcal{H}^\sharp\|_{\mathcal{L}(H_{\mathbb{T}^3}^{3/2+\varepsilon}, L_{\mathbb{T}^3}^2)} \to 0, \quad \delta \to 0^+. \tag{2.17}$$

A very useful consequence of Lemma 2.5 are the following bounds which we will use liberally. We denote the *form/energy domain* $\mathcal{D}(\sqrt{-\mathcal{H}})$ of the Anderson Hamiltonian to be the closure of the domain under the norm defined by $\|u\|^2_{\mathcal{D}(\sqrt{-\mathcal{H}})} = (u, -\mathcal{H} u)$.

**Lemma 2.6. (Embeddings and Equivalence of Norms)** *Recall that $(\cdot, \cdot)$ denotes the inner product in $L^2_{\mathbb{T}^3}$. The following norm equivalences are true:*

$$(-u, \mathcal{H} u)^{\frac{1}{2}} =: \|u\|_{\mathcal{D}(\sqrt{-\mathcal{H}})} \asymp \|e^{-W} u\|_{H^1_{\mathbb{T}^3}} \asymp \|e^{-W^M} u\|_{H^1_{\mathbb{T}^3}}, \tag{2.18}$$

$$\|\mathcal{H} u\|_{L^2} =: \|u\|_{\mathcal{D}(\mathcal{H})} \asymp \|(e^W \Gamma)^{-1} u\|_{H^2_{\mathbb{T}^3}} \asymp \|(e^{W^M} \Gamma)^{-1} u\|_{H^2_{\mathbb{T}^3}}, \tag{2.19}$$

*Moreover, one has the embeddings*

$$e^{-W^M} H^1_{\mathbb{T}^3} = \mathcal{D}(\sqrt{-\mathcal{H}}) \hookrightarrow W^{s,p}_{\mathbb{T}^3}, \text{ for } 0 \leq s < \frac{1}{2}, \quad p = \frac{6}{1+2s}, \tag{2.20}$$

and

$$e^{-W^M} \mathcal{D}(\mathcal{H}) \hookrightarrow W^{s,p}_{\mathbb{T}^3}, \text{ for } \frac{1}{2} < s < \frac{3}{2}, \quad p = \frac{6}{2s-1}, \tag{2.21}$$

$$H^2_{\mathbb{T}^3} = \Gamma^{-1} e^{-W^M} \mathcal{D}(\mathcal{H}) \hookrightarrow W^{s,p}_{\mathbb{T}^3}, \text{ for } \frac{1}{2} < s < 2, \quad p = \frac{6}{2s-1}, \tag{2.22}$$

and lastly,

$$\Gamma^{-1} e^{-W^M} \mathcal{D}(\mathcal{H}) \hookrightarrow C_{\mathbb{T}^3}^{\frac{1}{2}-}, e^{-W^M} \mathcal{D}(\mathcal{H}) \hookrightarrow C_{\mathbb{T}^3}^{\frac{1}{2}-}, \text{ and } \mathcal{D}(\mathcal{H}) \hookrightarrow C_{\mathbb{T}^3}^{\frac{1}{2}-}. \tag{2.23}$$

**Proof.** The first two points were already proven in [GUZ20], Proposition 2.57, and the embeddings follow from the Besov embeddings (c.f. e.g. [vZ22], Chapter 16), and the boundedness of multiplication of $\exp(\pm W^M)$ in $W^{s,p}_{\mathbb{T}^3}$ for $s < 1/2$ and the boundedness of $\Gamma$ on $W^{s,p}$ for $s < 3/2$ from Lemma 2.5. □

**Corollary 2.7.** *We have the following uniform in $\delta > 0$ bounds for the regularized objects*

$$\|f\|_{W^{s,p}_{\mathbb{T}^3}} \lesssim \|f\|_{\mathcal{D}\left((-\mathcal{H}_\delta)^{\frac{1}{2}}\right)} \asymp \|e^{-W_\delta}f\|_{H^1_{\mathbb{T}^3}}$$

$$\text{for } 0 \leq s < \frac{1}{2}, \quad p = \frac{6}{1+2s}$$

$$\|f\|_{C^{\frac{2}{3}-}_{\mathbb{T}^3}} + \|f\|_{W^{s,p}_{\mathbb{T}^3}} \lesssim \|f\|_{\mathcal{D}(\mathcal{H}_\delta)} \asymp \|\Gamma_\delta^{-1}e^{-W_\delta}f\|_{H^2_{\mathbb{T}^3}}$$

$$\text{for } \frac{1}{2} < s \leq 2, \quad p = \frac{6}{2s-1}.$$

So far the construction is essentially identical to that in [GUZ20], since the exponential Ansatz (2.2) differs from that in [GUZ20] only by a smooth function $\exp(P_{\leq M}W)$, while the merit of considering a high-frequency cutoff is shown in the following lemmata. We will use essentially such a transformation in the control of the second-order energy in Section 5.2, the use of which will be clarified in Remark 5.7.

**Lemma 2.8.** **($\Theta$-map)** *Let $\Theta: \mathcal{S}'(\mathbb{T}^3) \to \mathcal{S}'(\mathbb{T}^3)$ be defined by*

$$\Theta(u^\sharp) := -(1-\Delta)^{-1}([\exp(W^M) - 1] \oslash \Delta u^\sharp) + u^\sharp. \tag{2.24}$$

*We have for any $s \in \mathbb{R}$ that $\Theta: W^{s,p}_{\mathbb{T}^3} \to W^{s,p}_{\mathbb{T}^3}$ is continuous for $p \in [1,\infty)$, and further there exists an $M = M(\|\Xi\|_{\mathfrak{X}})$ large enough, such that $\Theta$ is invertible and $\Theta^{-1}: W^{s,p}_{\mathbb{T}^3} \to W^{s,p}_{\mathbb{T}^3}$ is continuous.*

**Proof.** It suffices to notice that by the paraproduct estimates from Lemma A.1 that

$$\|(1-\Delta)^{-1}([\exp(W^M)-1]\oslash \Delta u^\sharp)\|_{W^{s,p}_{\mathbb{T}^3}} \lesssim \|\exp(W^M)-1\|_{L^\infty_{\mathbb{T}^3}} \|u^\sharp\|_{W^{s,p}_{\mathbb{T}^3}}$$

$$\lesssim \|W^M\|_{L^\infty_{\mathbb{T}^3}} \|u^\sharp\|_{W^{s,p}_{\mathbb{T}^3}},$$

where we can further control $\|W^M\|_{L^\infty_{\mathbb{T}^3}} = \|P_{>M}W\|_{L^\infty_{\mathbb{T}^3}} \lesssim 2^{-M/2+}\|W\|_{C^{1/2-}_{\mathbb{T}^3}}$ by Bernstein's inequality and derive

$$\|(1-\Delta)^{-1}([\exp(W^M)-1]\oslash \Delta u^\sharp)\|_{W^{s,p}_{\mathbb{T}^3}} \leq \frac{1}{2}\|u^\sharp\|_{W^{s,p}_{\mathbb{T}^3}}, \tag{2.25}$$

*by choosing $M$ large enough. This shows that $\Theta$ is in fact invertible.* □

We now state two perturbative results which are true if one includes the $\Theta$ but which are not clear if one simply uses the mixed exponential-paracontrolled ansatz from [GUZ20]. In the two-dimensional case, one often uses a similar result e.g. in [Zac19] and [MZ24]. The transformation that accomplishes both (2.26) and (2.27) is the one of the technical novelties of this paper.

**Lemma 2.9.** *Let $p \in [2,\infty]$, given a realization $\Xi$ of the enhancement (2.7) of the spatial white noise $\xi$, we have that the map $\Lambda := \exp(W^M)\Gamma\Theta^{-1}$ has the following perturbative properties:*

$$\Lambda - \mathrm{id}: \quad W^{\beta,p}_{\mathbb{T}^3} \to W^{\frac{1}{2}+\beta-,p}_{\mathbb{T}^3}, \tag{2.26}$$

$$\Lambda^{-1}\mathcal{H}\Lambda - \Delta: \quad W^{\frac{3}{2}+\beta,p}_{\mathbb{T}^3} \to W^{\beta-,p}_{\mathbb{T}^3} \tag{2.27}$$

are bounded for $\beta \in \left(-\frac{1}{2}, \frac{1}{2}\right)$. The map $\Lambda$ equivalently parametrises the domain of $\mathcal{H}$, i.e. $\mathcal{D}(\mathcal{H}) = e^W \Gamma(H^2_{\mathbb{T}^3}) = \Lambda(H^2_{\mathbb{T}^3})$.

**Proof.** For the first bound, we use Lemmas 2.5 and 2.8 on $(1-\Delta)[\exp(W^M)\Gamma u^\sharp - \Theta u^\sharp]$ to derive[2.3]

$$
\begin{aligned}
(1-\Delta)[\exp(W^M)\Gamma u^\sharp - \Theta u^\sharp] \stackrel{(2.24)}{=}{} & (1-\Delta)([\exp(W^M)-1]\Gamma u^\sharp) + (1-\Delta)\Gamma u^\sharp \\
& - (1-\Delta)u^\sharp + [\exp(W^M)-1]\otimes\Delta u^\sharp \\
={} & (1-\Delta)([\exp(W^M)-1] \succcurlyeq \Gamma u^\sharp) + (1-\Delta)(\Gamma-\mathrm{id})u^\sharp \\
& + [\exp(W^M)-1]\otimes\Gamma u^\sharp - \Delta[\exp(W^M)]\otimes\Gamma u^\sharp \\
& - 2\nabla[\exp(W^M)]\otimes\nabla\Gamma u^\sharp \\
& + [\exp(W^M)-1]\otimes\Delta[(\Gamma-\mathrm{id})u^\sharp],
\end{aligned}
$$

together with (2.13) which implies that the second and the last term on the right gains $1/2-$ regularity as desired. Now notice that the left hand side above is $(1-\Delta)$ applied to

$$\exp(W^M)\Gamma u^\sharp - \Theta u^\sharp = [\exp(W^M)\Gamma\Theta^{-1} - \mathrm{id}]\Theta u^\sharp = (\Lambda - \mathrm{id})\Theta u^\sharp,$$

and (2.26) follows using also the boundedness of $\Theta^{\pm 1}$ from Lemma 2.8.

For the second bound, we have

$$
\begin{aligned}
[\mathcal{H}\exp(W^M)\Gamma - \Delta\Theta]u^\sharp \stackrel{(2.24)}{=}{} & \exp(W^M)(\Delta u^\sharp + [(1-\Delta)Z]\odot u^\sharp + 2\nabla W^M \odot \nabla u^\sharp + G(u^\flat)) \\
& + (1-\Delta)u^\sharp - [\exp(W^M)-1]\otimes\Delta u^\sharp - \Theta(u^\sharp) \\
={} & [\exp(W^M)-1] \succcurlyeq \Delta u^\sharp - (1-\Delta)^{-1}([\exp(W^M)-1]\otimes\Delta u^\sharp) \\
& + \exp(W^M)([(1-\Delta)Z]\odot u^\sharp + 2\nabla W^M \odot \nabla u^\sharp + G(u^\flat)).
\end{aligned}
$$

We have used crucially the transformation (2.24) such that the worst terms $(1-\Delta)([\exp(W^M)-1]\otimes\Gamma u^\sharp)$ and $[\exp(W^M)-1]\otimes\Delta u^\sharp$ have a cancellation, since

$$(1-\Delta)\Theta(u^\sharp) = (1-\Delta)u^\sharp - [\exp(W^M)-1]\otimes\Delta u^\sharp,$$

and

$$u^\sharp - \Theta(u^\sharp) = (1-\Delta)^{-1}([\exp(W^M)-1]\otimes\Delta u^\sharp). \qquad \square$$

**Notation 2.10.** *We see in Lemma 2.8 that we only need to fix the value of M globally depending on the norm of the noise enhancement and similarly we will fix N and we will largely drop these harmless parameters in order to not clutter the notation. Moreover, we usually do not refer to the fixed realisation $\omega$ of the noise. In particular, M and N can both be chosen to hold uniformly in $\delta$ if one has a smooth approximation $\Xi_\delta$ as in (2.7).*

To conclude the discussion on the transformations of Anderson Hamiltonian as a perturbation of the Laplacian, we prove that $(-\mathcal{H}_\delta)^{-1/2} \to (-\mathcal{H})^{-1/2}$ as $\delta \to 0^+$ w.r.t. the operator norm, which plays a role in the proof of the convergence to the BBGKY hierarchy in Section 6.2. For that purpose, we need the short-time heat kernel estimate of the $\mathcal{H}^\sharp$ in Lemma 2.12 and we refer the readers to [BDM22] for the detailed analysis performed with the Anderson Hamiltonian on surfaces, without seeking a complete description of the heat flow generated by $\mathcal{H}$ on $\mathbb{T}^3$.

---

2.3. We recall the abbreviation $f \succcurlyeq g := f \otimes g + f \odot g$ for distributions $f$ and $g$.

**Notation 2.11.** *Thanks to Lemma 2.3, we can define a bounded Functional Calculus $F(-\mathcal{H}) \in \mathcal{L}(L^2_{\mathbb{T}^3})$ with $-\mathcal{H} \geq c_\Sigma$ and a bounded continuous function $F$ on $[c_\Sigma, \infty)$ for some $c_\Sigma > 0$. For any such function, we will also write*

$$F(-\mathcal{H}^\natural) := \Lambda^{-1} F(-\mathcal{H}) \Lambda \in \mathcal{L}(L^2_{\mathbb{T}^3}). \tag{2.28}$$

*Furthermore, we have*

$$F(-\mathcal{H}^\natural) = (\mathcal{H}\Lambda)^{-1} F(-\mathcal{H}) \mathcal{H}\Lambda \in \mathcal{L}(H^2_{\mathbb{T}^3}), \tag{2.29}$$

*as $\mathcal{H}\Lambda : H^2_{\mathbb{T}^3} \to L^2_{\mathbb{T}^3}$ and $(\mathcal{H}\Lambda)^{-1} : L^2_{\mathbb{T}^3} \to H^2_{\mathbb{T}^3}$ are continuous. Therefore, $F(-\mathcal{H}^\natural)$ can be defined as bounded linear operators on $H^s_{\mathbb{T}^3}$ for $s \in [0, 2]$ via interpolation.*

**Lemma 2.12.** *The following heat-kernel estimate of $(e^{t\mathcal{H}^\natural})_{t \geq 0}$ holds for $t > 0$:*

$$\|e^{t\mathcal{H}^\natural} u\|_{H^\gamma_{\mathbb{T}^3}} \lesssim t^{-\frac{\gamma}{2}} \|u\|_{L^2_{\mathbb{T}^3}} \quad 0 \leq \gamma \leq 2. \tag{2.30}$$

*We have also the additional regularization of the flow*

$$\|e^{t\mathcal{H}^\natural} u\|_{W^{2,\infty}_{\mathbb{T}^3}} \lesssim \lceil t \rceil^{1+} t^{-\frac{7}{4}-} \|u\|_{L^2_{\mathbb{T}^3}}, \tag{2.31}$$

*for any $t > 0$ and the following boundedness property*

$$\|e^{t\mathcal{H}^\natural} u\|_{W^{s,p}_{\mathbb{T}^3}} \lesssim \|u\|_{W^{s,p}_{\mathbb{T}^3}} \tag{2.32}$$

*for $0 \leq t \leq 1$ and $\frac{3}{2} < s \leq 2$, $p \in [2, \infty]$.*

The proof is postponed to Appendix B.1. The lemma on the convergence of square-root resolvent of the Anderson Hamiltonian now follows.

**Lemma 2.13.** *For any given realization of the enhancement $\Xi$ we have the operator-norm convergence*

$$\|(-\mathcal{H}_\delta)^{-1/2} - (-\mathcal{H})^{-1/2}\|_{\mathcal{L}(L^2_{\mathbb{T}^3})} \to 0, \tag{2.33}$$

*as $\delta \to 0^+$.*

**Proof.** This follows by the functional calculus for self-adjoint operators, see e.g. Theorem VIII.20 in [RS80]. □

## 3 Strichartz Estimates

In this section, we prove the Strichartz-type estimates for the Anderson Hamiltonian on the three-dimensional torus following closely [Zac19] and the Strichartz estimates established in [DJLM23], which in turn uses the celebrated decoupling result of Bourgain-Demeter [BD15].

**Theorem 3.1.** ([DJLM23]) *For any $\varepsilon > 0$ and Strichartz pair $(q, r)$ such that*

$$\frac{2}{q} + \frac{3}{r} \geq \frac{3}{2}, \quad q \in \left[\frac{10}{3}, \infty\right], \quad r \in \left[2, \frac{10}{3}\right], \tag{3.1}$$

we have the homogeneous Strichartz estimate for any $s \in \mathbb{R}$ that

$$\|e^{-it\Delta} u\|_{L^q_{t;[0,1]} W^{s,r}_{\mathbb{T}^3}} \lesssim_{q,r,\varepsilon} \|u\|_{H^{s+\varepsilon}_{\mathbb{T}^3}}. \tag{3.2}$$

It follows immediately from (3.2) and Minkowski inequality that for any $0 \leq t_0 < t_1 \leq 1$

$$\left\| \int_{t_0}^{t} e^{-i(t-s)\Delta} v(s) \mathrm{d}s \right\|_{L^q_{t;[t_0,t_1]} W^{s,r}_{\mathbb{T}^3}} \leq \int_{t_0}^{t_1} \|e^{-i(t-s)\Delta} v(s)\|_{L^q_{t;[s,t_1]} W^{s,r}_{\mathbb{T}^3}} \mathrm{d}s \overset{(3.2)}{\lesssim} \|v\|_{L^1_{[t_0,t_1]} H^{s+\varepsilon}_{\mathbb{T}^3}}. \tag{3.3}$$

**Remark 3.2.** Since the argument is perturbative, it is natural to consider the linear flow of the transformed equations $(e^{-it\mathcal{H}^\#})_{t \in \mathbb{R}}$ and respectively $(e^{-it\mathcal{H}^\natural})_{t \in \mathbb{R}}$, which are defined to be the solution maps of

$$i\partial_t u^\# = \mathcal{H}^\# u^\# = (e^W \Gamma)^{-1} \mathcal{H}(e^W \Gamma) u^\#, \text{ so } e^{-it\mathcal{H}^\#} := (e^W \Gamma)^{-1} e^{-it\mathcal{H}} \Gamma e^W \tag{3.4}$$

and respectively

$$i\partial_t u^\natural = \mathcal{H}^\natural u^\natural := \Lambda^{-1} \mathcal{H} \Lambda u^\natural, \quad \text{so } e^{-it\mathcal{H}^\natural} := \Lambda^{-1} e^{-it\mathcal{H}} \Lambda. \tag{3.5}$$

One derives immediately the equivalence of equations

$$\begin{cases} i\partial_t u = \mathcal{H} u \\ u(0) = u_0 \end{cases} \iff \begin{cases} i\partial_t u^\# = \mathcal{H}^\# u^\# \\ u^\#(0) = (e^W \Gamma)^{-1} u_0 \end{cases} \iff \begin{cases} i\partial_t u^\natural = \mathcal{H}^\natural u^\natural \\ u^\natural(0) = \Lambda^{-1} u_0 \end{cases}$$

since the transforms are all invertible in the spaces under consideration.

The essential observation for the perturbative argument is the following bound for any $\varepsilon > 0$, $\beta \in [\varepsilon, 1/2]$ and $p \in [1, \infty]$:

$$\|(\mathcal{H}^\natural - \Delta) u^\natural\|_{W^{\beta-\varepsilon,p}_{\mathbb{T}^3}} \lesssim \|u^\natural\|_{W^{\beta+3/2,p}_{\mathbb{T}^3}}, \tag{3.6}$$

see Lemma 2.9.

We first need to show the continuity of the flows $(e^{-it\mathcal{H}^\#})_{t \in \mathbb{R}}$ and $(e^{-it\mathcal{H}^\natural})_{t \in \mathbb{R}}$:

**Lemma 3.3.** *The linear flows $(e^{-it\mathcal{H}^\#})_{t \in \mathbb{R}}$ and $(e^{-it\mathcal{H}^\natural})_{t \in \mathbb{R}}$ are strongly continuous groups on $H^s_{\mathbb{T}^3}$ for $s \in [0, 2]$.*

**Proof.** By the construction of the unitary group $(e^{-it\mathcal{H}})_{t \in \mathbb{R}}$ in [GUZ20], $t \mapsto e^{-it\mathcal{H}}$ is continuous both on the domain $\mathcal{D}(\mathcal{H})$ and $L^2_{\mathbb{T}^3}$, which implies that $(e^{-it\mathcal{H}^\#})_{t \in \mathbb{R}}$ and $(e^{-it\mathcal{H}^\natural})_{t \in \mathbb{R}}$ are strongly continuous on $H^2_{\mathbb{T}^3}$ due to (2.19) and continuous on $L^2_{\mathbb{T}^3}$, since $e^W, \Gamma, \Theta$ and their inverses are continuous on $L^2_{\mathbb{T}^3}$. The claim then follows by interpolation as in [Zac19] and [MZ24]. □

As alluded to above, the first Strichartz-type estimate we have for the sharpened Anderson Hamiltonian $\mathcal{H}^\#$ follows from Theorem 3.1 and a perturbative argument inspired by [BGT04], which however only works in the *low-regularity* regime.

**Lemma 3.4.** *For any $\varepsilon > 0$, $s \in \left[ \frac{3}{2q} + \varepsilon, \frac{3}{2q} + \frac{1}{2} \right]$ and Strichartz pair $(q, r)$ as in (3.1), we have the homogeneous Strichartz-type estimate of $\mathcal{H}^\#$, that is*

$$\|e^{-it\mathcal{H}^\natural} u^\#\|_{L^q_{[0,1]} W^{s-\frac{3}{2q}-\varepsilon,r}_{\mathbb{T}^3}} \lesssim \|u^\natural\|_{H^s_{\mathbb{T}^3}}. \tag{3.7}$$

**Proof.** We have due to

$$[e^{-i(t-t_0)\mathcal{H}^\natural} - e^{-i(t-t_0)\Delta}]u^\natural = -i\int_{t_0}^t e^{-i(t-s)\Delta}(\mathcal{H}^\natural - \Delta)e^{-i(s-t_0)\mathcal{H}^\natural}u^\sharp\,\mathrm{d}s$$

that for each $\varepsilon > 0$, any $j^{\text{th}}$ Littlewood-Paley block $\Delta_j u$ of $u$ for $j \geq -1$ and arbitrary partion $\{I\}$ of $[0,1]$,

$$\begin{aligned}
\|e^{-it\mathcal{H}^\natural}\Delta_j u^\natural\|_{L^q_{[0,1]}W^{s,r}_{\mathbb{T}^3}}^q &= \sum_{I=[t_0,t_1]} \int_I \|e^{-i(t-t_0)\mathcal{H}^\natural}e^{-it_0\mathcal{H}^\natural}\Delta_j u^\natural\|_{W^{s,r}_{\mathbb{T}^3}}^q \mathrm{d}t \\
&\lesssim \sum_{I=[t_0,t_1]} \int_I \|e^{-i(t-t_0)\Delta}e^{-it_0\mathcal{H}^\natural}\Delta_j u^\natural\|_{W^{s,r}_{\mathbb{T}^3}}^q \mathrm{d}t \\
&\quad + \sum_{I=[t_0,t_1]} \int_I \left\|\int_{t_0}^t e^{-i(t-s)\Delta}(\mathcal{H}^\natural - \Delta)e^{-is\mathcal{H}^\natural}\Delta_j u^\natural\,\mathrm{d}s\right\|_{W^{s,r}_{\mathbb{T}^3}}^q \mathrm{d}t \\
&\overset{(3.2)}{\lesssim} \sum_{I=[t_0,t_1]} \|e^{-it_0\mathcal{H}^\natural}\Delta_j u^\natural\|_{H^{s+\varepsilon}_{\mathbb{T}^3}}^q + \sum_{I=[t_0,t_1]} \|e^{-is\mathcal{H}^\natural}\Delta_j u^\natural\|^q_{L^1_s;H^{s+\frac{3}{2}+2\varepsilon}_{\mathbb{T}^3}} \\
&\lesssim \#(I)\left[\|\Delta_j u^\natural\|_{H^{s+\varepsilon}_{\mathbb{T}^3}}^q + \sum_{I=[t_0,t_1]} 2^{\frac{3}{2}jq}|I|^q\|\Delta_j u^\natural\|_{H^{s+2\varepsilon}_{\mathbb{T}^3}}^q\right],
\end{aligned}$$

for any $s \in [0, 1/2 - 2\varepsilon]$, where we have used Lemma 3.3 for the last inequality. If we further choose for each frequency $j \geq -1$ the partition of the time interval $|I_j| \sim 2^{-3j/2}$, we have

$$\|e^{-it\mathcal{H}^\natural}u^\natural\|_{L^q_{[0,1]}W^{s,r}_{\mathbb{T}^3}} \lesssim \sum_{j\geq -1}\|e^{-it\mathcal{H}^\natural}\Delta_j u^\natural\|_{L^q_{[0,1]}W^{s,r}_{\mathbb{T}^3}} \lesssim \sum_{j\geq -1} 2^{\frac{3}{2q}j}\|\Delta_j u^\natural\|_{H^{s+2\varepsilon}} \lesssim \|u^\natural\|_{H^{\frac{3}{2q}+s+3\varepsilon}}, \qquad (3.8)$$

where we have applied Bernstein inequality. This then gives (3.12) for $s \in [3/(2q) + \varepsilon, 3/(2q) + 1/2]$. □

Clearly the same statement is true for $e^{-it\mathcal{H}^\sharp}$ by the boundedness of $\Theta$. In order to extend Lemma 3.4 up to $s = 2$, we need the following characterization of Sobolev norms with the transformed operator $\mathcal{H}^\natural$.

**Lemma 3.5.** *For any $\beta$ such that $3/2 < \beta \leq 2$ and $p \in [2,\infty)$, the following bounds for the operator*

$$(-\mathcal{H}^\natural)^{\beta/2} := \Lambda^{-1}(-\mathcal{H})^{\beta/2}\Lambda$$

*hold:*

$$\|u^\natural\|_{W^{\beta-,p}_{\mathbb{T}^3}} \lesssim \left\|(1-\Delta)^{\frac{\beta}{2}}u^\natural\right\|_{L^p_{\mathbb{T}^3}} \asymp \left\|(-\mathcal{H}^\natural)^{\frac{\beta}{2}}u^\natural\right\|_{L^p_{\mathbb{T}^3}} \lesssim \|u^\natural\|_{W^{\beta,p}_{\mathbb{T}^3}}. \qquad (3.9)$$

**Proof.** The first and last bounds follow from the embeddings between fractional Sobolev and Bessel potential spaces (A.1). The main part is of course the norm equivalence.

We first consider $0 \leq \alpha < \frac{3}{2}$ and write the difference using functional calculus

$$(-\mathcal{H}^\natural)^{-\frac{\alpha}{2}} - (1-\Delta)^{-\frac{\alpha}{2}} = \int_0^1 t^{\frac{\alpha}{2}-1}[e^{t\mathcal{H}^\natural} - e^{-t(1-\Delta)}]\mathrm{d}t + \int_1^\infty t^{\frac{\alpha}{2}-1}[e^{t\mathcal{H}^\natural} - e^{-t(1-\Delta)}]\mathrm{d}t \qquad (3.10)$$

where the second term is a bounded operator from $L^2_{\mathbb{T}^3}$ to $W^{2,\infty}_{\mathbb{T}^3}$ by (2.31) for the $e^{t\mathcal{H}^\natural}$ term and by standard tools for the $e^{t\Delta}$ term. For the first term after the equality above, we have according to the heat-kernel estimate of $e^{t\Delta}$ and (3.6) that

$$\left\|\int_0^1 t^{\frac{\alpha}{2}-1}[e^{t\mathcal{H}^\natural} - e^{-t(1-\Delta)}]u\,dt\right\|_{W^{2,p}_{\mathbb{T}^3}}$$
$$\leq \int_0^1 t^{\frac{\alpha}{2}-1}\left\|\int_0^t e^{(t-s)(\Delta-1)}(\mathcal{H}^\natural - \Delta + 1)e^{s\mathcal{H}^\natural}u\,ds\right\|_{W^{2,p}_{\mathbb{T}^3}}dt$$
$$\lesssim \int_0^1 t^{\frac{\alpha}{2}-1}\int_0^t (t-s)^{-1+\varepsilon}ds\|u\|_{W^{\frac{3}{2}+2\varepsilon,p}_{\mathbb{T}^3}}dt \lesssim \|u\|_{W^{\frac{3}{2}+2\varepsilon,p}_{\mathbb{T}^3}}, \quad (3.11)$$

where we have used $\alpha \geq 0$ in the last inequality and (2.32). We now fix some $\beta \in (3/2, 2]$ and derive from the expression (3.10) above and

$$(-\mathcal{H}^\natural)^{\frac{\beta}{2}} - (1-\Delta)^{\frac{\beta}{2}} = (1-\Delta)\left[(-\mathcal{H}^\natural)^{\frac{\beta}{2}-1} - (1-\Delta)^{\frac{\beta}{2}-1}\right]$$
$$+ (\Delta - 1 - \mathcal{H}^\natural)\left[(-\mathcal{H}^\natural)^{\frac{\beta}{2}-1} - (1-\Delta)^{\frac{\beta}{2}-1}\right]$$
$$+ (\Delta - 1 - \mathcal{H}^\natural)(1-\Delta)^{\frac{\beta}{2}-1},$$

that for any $\varepsilon > 0$ small enough,

$$\left\|\left[(-\mathcal{H}^\natural)^{\frac{\beta}{2}} - (1-\Delta)^{\frac{\beta}{2}}\right]u\right\|_{L^p_{\mathbb{T}^3}} \stackrel{(3.11)}{\lesssim} \|u\|_{W^{\frac{3}{2}+\varepsilon,p}_{\mathbb{T}^3}} + \left\|(\Delta - 1 - \mathcal{H}^\natural)(1-\Delta)^{\frac{\beta}{2}-1}u\right\|_{L^p_{\mathbb{T}^3}}$$
$$\lesssim \|u\|_{W^{\frac{3}{2}+\varepsilon,p}_{\mathbb{T}^3}} + \|u\|_{W^{\beta-\frac{1}{2}+\varepsilon,p}_{\mathbb{T}^3}} \lesssim \|u\|_{W^{\frac{3}{2}+\varepsilon,p}_{\mathbb{T}^3}},$$

where we have used again (3.6). To conclude with (3.9), it suffices to choose $\varepsilon < \beta - 3/2$ and bound $\|u\|_{W^{\frac{3}{2}+\varepsilon,p}_{\mathbb{T}^3}} \lesssim \|u\|_{W^{\beta,p}_{\mathbb{T}^3}}$. $\square$

Due to the continuity estimates of $\Theta$ and its inverse established in Lemma 2.8, it is straightforward to translate the Strichartz estimate established in Lemma 3.4 for the group $(e^{-it\mathcal{H}^\natural})_{t\in\mathbb{R}}$ and $(e^{-it\mathcal{H}^\sharp})_{t\in\mathbb{R}}$.

**Corollary 3.6.** *For any small $\varepsilon > 0$, we have the homogeneous Strichartz-type estimate holds for $\mathcal{H}^\sharp$ and $\mathcal{H}^\natural$, that is*

$$\|e^{-it\mathcal{H}^\sharp}u^\sharp\|_{L^{\frac{10}{3}}_{[0,1]}W^{s-\frac{9}{20}-\varepsilon,\frac{10}{3}}_{\mathbb{T}^3}} \lesssim \|u^\sharp\|_{H^s_{\mathbb{T}^3}}, \quad (3.12)$$

*and*

$$\|e^{-it\mathcal{H}^\natural}u^\natural\|_{L^{\frac{10}{3}}_{[0,1]}W^{s-\frac{9}{20}-\varepsilon,\frac{10}{3}}_{\mathbb{T}^3}} \lesssim \|u^\natural\|_{H^s_{\mathbb{T}^3}}, \quad (3.13)$$

*where $s \in \left[\frac{9}{20} + \varepsilon, 2\right]$.*

**Proof.** The bound (3.12) follows from (3.13) and Lemma 2.8, while the latter with $s = 2$ is shown directly as follows

$$\|e^{-it\mathcal{H}^\natural}u^\natural\|_{L^{\frac{10}{3}}_{t;[0,1]}W^{\frac{31}{20}-2\varepsilon,\frac{10}{3}}_{\mathbb{T}^3}} \stackrel{(3.9)}{\lesssim} \left\|(-\mathcal{H}^\natural)^{\frac{31}{40}-\varepsilon}e^{-it\mathcal{H}^\natural}u^\natural\right\|_{L^{\frac{10}{3}}_{t;[0,1]}L^{\frac{10}{3}}_{\mathbb{T}^3}}$$
$$\stackrel{(3.7)}{\lesssim} \left\|(-\mathcal{H}^\natural)^{\frac{31}{40}-\varepsilon}u^\natural\right\|_{H^{\frac{9}{20}+2\varepsilon}_{\mathbb{T}^3}} \stackrel{(3.9)}{\lesssim} \|u^\natural\|_{H^2_{\mathbb{T}^3}}.$$

We can then conclude the proof with an interpolation argument. □

# 4 Local Well-posedness of ANLS

Consider now a symmetric potential $V_\beta \in W^{\beta,1}$ with $\beta \in (0,3)$, one can think of the generalised Coulomb potential, which is

$$V_\beta = \frac{c_\beta}{|\cdot|^{3-\beta}} \quad \text{near the origin}^{4.1}, \tag{4.1}$$

and smooth elsewhere.

In this section, we prove local well-posedness (LWP) of the nonlinear Schrödinger equation with the Anderson Hamiltonian and Hartree nonlinearity

$$i\partial_t u = \mathcal{H} u \pm u(|u|^2 * V_\beta) \quad \text{on } [0,T] \times \mathbb{T}^3, \tag{4.2}$$

for some $T$ depending on the initial data, while the initial data are given either in the domain or the energy domain of $\mathcal{H}$ defined in Definition 2.3 and right above Lemma 2.6 respectively. When the initial data are in the domain of $\mathcal{H}$, we will instead consider the (more general) local-wellposedness of the nonlinear (focusing or defocusing) Schrödinger equation with the Anderson Hamiltonian $\mathcal{H}$ and polynomial nonlinearity on $\mathbb{T}^3$:

$$i\partial_t u = \mathcal{H} u \pm u|u|^{p-1} \quad \text{on } [0,T] \times \mathbb{T}^3, \tag{4.3}$$

where $p > 1$. The LWP of (4.2) in the domain of $\mathcal{H}$ follows from exactly the same argument as that of (4.3), since the convolution with $V_\beta$ only improves the regularity. We thus focus on the latter in this section.

## 4.1 Initial Datum in the Domain

We consider in this section the local well-posedness (LWP) problem of (4.3) in the domain of $\mathcal{D}(\mathcal{H})$ with the mild formulation and Banach's fixed point argument, that is, with

$$u(t) = e^{-it\mathcal{H}} u_0 \mp i \int_0^t e^{-is\mathcal{H}} [u(t-s)|u(t-s)|^{p-1}] ds, \tag{4.4}$$

and $u_0 \in \mathcal{D}(\mathcal{H})$, and equivalently (according to Lemma 2.5)

$$u^\natural(t) = e^{-it\mathcal{H}^\natural} \Lambda^{-1} u_0 \mp i \int_0^t e^{-is\mathcal{H}^\natural} \Lambda^{-1} [\Lambda u^\natural(t-s)|\Lambda u^\natural(t-s)|^{p-1}] ds. \tag{4.5}$$

The result of the LWP of (4.3) reads as follows:

**Proposition 4.1.** *For fixed enhancement $\Xi = \Xi_\omega$, (4.4) is locally well-posedness, in the sense that given any initial data $u_0 \in \mathcal{D}(\mathcal{H})$, there exists*

$$T = T(\Xi, \|\mathcal{H} u_0\|_{L^2}) \in (0, \infty],$$

---

4.1. In particular, we have $V_\beta \in L_w^{\frac{3}{3-\beta}}(\mathbb{T}^3)$.

such that there exists a unique solution $u \in C([0,T], \mathcal{D}(\mathcal{H})) \cap C^1([0,T], L^2_{\mathbb{T}^3})$, depending continuously on the data, and one has the blow-up alternative: there either exists a maximal time $T_{\max} < \infty$ with

$$\lim_{t \to T_{\max}} \|\mathcal{H} u(t)\|_{L^2} = \infty, \tag{4.6}$$

or $T_{\max} = \infty$.

**Proof.** This is quite similar to the proof of Theorem 3.6 in [GUZ20]. We abbreviate $\mathcal{X} := C([0,T], \mathcal{D}(\mathcal{H})) \cap C^1([0,T], L^2)$ with norm

$$\|u\|_{\mathcal{X}} := \|u\|_{C([0,T], \mathcal{D}(\mathcal{H}))} + \|\partial_t u\|_{C([0,T], L^2_{\mathbb{T}^3})} \tag{4.7}$$

and consider initial condition $\|u_0\|_{\mathcal{D}(\mathcal{H})} \leq R$ for some $R < \infty$. According to Banach's fixed point, it suffices to prove that the map $\mathcal{S}$ on the ball $B_{\mathcal{X}}(M_R)$, where

$$(\mathcal{S} u)(t) := e^{-it\mathcal{H}} u_0 \mp i \int_0^t e^{-is\mathcal{H}} (u(t-s)|u(t-s)|^{p-1}) ds,$$

and $M_R < \infty$ is a constant to be decided later depending on $\Xi$ and $R$, is a contraction for $T = T(R, \Sigma_w) > 0$ small enough. Indeed, $\mathcal{S}$ is a well-defined map, since

$$\begin{aligned}
\mathcal{H}\mathcal{S}(u(t)) &= e^{-it\mathcal{H}} \mathcal{H} u_0 \pm \int_0^t \partial_s(e^{-it\mathcal{H}})(u(t-s)|u(t-s)|^{p-1}) ds \\
&= e^{-it\mathcal{H}} \mathcal{H} u_0 \pm e^{-it\mathcal{H}}(u_0|u_0|^{p-1}) \mp u(t)|u(t)|^{p-1} \mp \int_0^t e^{-is\mathcal{H}} \partial_s(u(t-s)|u(t-s)|^{p-1}) ds,
\end{aligned}$$

where the maps $t \mapsto e^{-it\mathcal{H}} \mathcal{H} u_0$, $t \mapsto e^{-it\mathcal{H}}(u_0|u_0|^{p-1})$ and $t \mapsto u(t)(|u(t)|^{p-1})$ are strongly continuous in $L^2_{\mathbb{T}^3}$ by Stone's theorem, using that $\mathcal{H} u_0 \in L^2_{\mathbb{T}^3}$, and $\mathcal{D}(\mathcal{H})$ continuously embeds into $L^\infty_{\mathbb{T}^3}$. Besides,

$$\begin{aligned}
&\int_0^t e^{-is\mathcal{H}} \partial_s(u(t-s)|u(t-s)|^{p-1}) ds \\
&= -\frac{p+1}{2} \int_0^t e^{-is\mathcal{H}} [(\partial_t u(t-s))|u(t-s)|^{p-1}] ds \\
&\quad - \frac{p-1}{2} \int_0^t e^{-is\mathcal{H}} [u^2(t-s)|u(t-s)|^{p-3} \partial_t \bar{u}(t-s)] ds,
\end{aligned} \tag{4.8}$$

which is continuous in $L^2_{\mathbb{T}^3}$, since $t \mapsto \partial_t u(t)$ is continuous in $L^2_{\mathbb{T}^3}$, and $t \mapsto u(t)$ is continuous in $L^\infty_{\mathbb{T}^3}$. Similarly, $\partial_t \mathcal{S}(u) \in C([0,T], L^2)$, since

$$\partial_t \mathcal{S}_t(u) = -i e^{-it\mathcal{H}} \mathcal{H} u_0 \mp i \int_0^t e^{-is\mathcal{H}} \partial_t(u(t-s)|u(t-s)|^{p-1}) ds \mp i e^{-it\mathcal{H}}(u_0|u_0|^{p-1}),$$

where the second term can be treated as in (4.7). Now we can set up the fixed point by noticing

$$\begin{aligned}
\|\mathcal{S}(u)\|_{\mathcal{X}} &= \|\mathcal{H}\mathcal{S}(u)\|_{C_{[0,T]} L^2_{\mathbb{T}^3}} + \|\partial_t \mathcal{S}(u)\|_{C_{[0,T]} L^2_{\mathbb{T}^3}} \\
&\lesssim_{\Xi,p} \|\mathcal{H} u_0\|_{L^2} + \|u_0|u_0|^{p-1}\|_{L^2_{\mathbb{T}^3}} + \|u|u|^{p-1}\|_{C_{[0,T]} L^2_{\mathbb{T}^3}} + \int_0^T \|\partial_t u(t)|u(t)|^{p-1}\|_{L^2_{\mathbb{T}^3}} dt \\
&\quad + \int_0^T \|u^2(t-s)|u(t-s)|^{p-3} \partial_t \bar{u}(t)\|_{L^2_{\mathbb{T}^3}} dt \\
&\lesssim_{\Xi,p} \|\mathcal{H} u_0\|_{L^2} + \|\mathcal{H} u_0\|_{L^2}^p + T\|\partial_t u\|_{C_{[0,T]} L^2_{\mathbb{T}^3}} \|u\|_{C_{[0,T]} \mathcal{D}(\mathcal{H})}^{p-1},
\end{aligned}$$

since

$$\begin{aligned}
\|u|u|^{p-1}\|_{C_{[0,T]} L^2_{\mathbb{T}^3}} &\leq \|u_0|u_0|^{p-1}\|_{L^2_{\mathbb{T}^3}} + \int_0^T \|\partial_t(u|u|^{p-1})\|_{L^2_{\mathbb{T}^3}} dt \\
&\lesssim \|u_0|u_0|^{p-1}\|_{L^2_{\mathbb{T}^3}} + T\|\partial_t u\|_{C_{[0,T]} L^2_{\mathbb{T}^3}} \|u\|_{C_{[0,T]} \mathcal{D}(\mathcal{H})}^{p-1}.
\end{aligned}$$

We have therefore $\|\mathcal{S}(u)\|_{\mathcal{X}} \lesssim_{\Xi,p} R + R^p + TM_R^p \lesssim R^p \vee 1 + TM_R^p$, and we can then first choose $M_R \approx_{\Xi,p} R^p \vee 1$, so that

$$\|\mathcal{S}(u)\|_{\mathcal{X}} \le \frac{1}{2}M_R + C(\Xi,p)TM_R^p,$$

and then $T = T(R,\Xi)$ small enough such that $\|\mathcal{S}(u)\|_{\mathcal{X}} \le M_R$, that is,

$$\|u\|_{\mathcal{X}} \le M_R \Longrightarrow \|\mathcal{S}(u)\|_{\mathcal{X}} \le M_R.$$

To get a contraction, we have by integration by parts in the time-variable again for any $u,v \in \mathcal{X}$ with $u_0 = v_0$ that

$$\mathcal{H}[\mathcal{S}(u(t)) - \mathcal{S}(v(t))] = u(t)|u(t)|^{p-1} - v(t)|v(t)|^{p-1} \tag{4.9}$$
$$\mp \int_0^t e^{-is\mathcal{H}} \partial_s[u(t-s)(|u(t-s)|^{p-1}) - v(t-s)(|v(t-s)|^{p-1})]ds,$$

since the initial data coincide. The first term can then be controlled in $L^2$ by

$$\|u|u|^{p-1}(t) - v|v|^{p-1}(t)\|_{L^2_{\mathbb{T}^3}} \le T\|\partial_t(u|u|^{p-1} - v|v|^{p-1})\|_{C([0,T],L^2_{\mathbb{T}^3})}$$

which we further split into

$$\partial_t(u|u|^{p-1} - v|v|^{p-1}) = \frac{p+1}{2}[\partial_t(u-v)|u|^{p-1} + \partial_t v(|u|^{p-1} - |v|^{p-1})]$$
$$+ \frac{p-1}{2}[(u^2|u|^{p-3} - v^2|v|^{p-3})\partial_t \bar{u} + v^2|v|^{p-3}\partial_t(\overline{u-v})].$$

Therefore, we have

$$\|u|u|^{p-1} - v|v|^{p-1}\|_{C([0,T],L^2_{\mathbb{T}^3})} \lesssim_{\Xi,p} T\|u-v\|_{\mathcal{X}}(\|u\|_{\mathcal{X}}^{p-1} + \|v\|_{\mathcal{X}}^{p-1} + 1),$$

and with the same argument, we can bound the second difference in (4.9) so that

$$\|\mathcal{S}(u) - \mathcal{S}(v)\|_{\mathcal{X}} \lesssim_{\Xi,p} T\|u-v\|_{\mathcal{X}}(\|u\|_{\mathcal{X}}^{p-1} + \|v\|_{\mathcal{X}}^{p-1} + 1), \tag{4.10}$$

so that for $T = T(R,\Xi)$ small enough, the map $\mathcal{S}$ is a contraction on $B_\Xi(M_R)$ as claimed, (4.6) is proved as in [Caz03]. □

## 4.2 Initial Datum in the Energy Domain

Now we consider the LWP of the equation (4.2) when the initial datum is in the energy domain $\mathcal{D}(\mathcal{H})$ of $\mathcal{H}$, using in particular the Strichartz estimates that we have established in Lemma 3.4. We recall that $u \in \mathcal{D}(\mathcal{H})$ if and only if $(e^W\Gamma)^{-1}u \in H^2_{\mathbb{T}^3}$ according to (2.19) and

$$\|u\|_{\mathcal{D}(\sqrt{-\mathcal{H}})} \asymp \|e^{-W}u\|_{H^1_{\mathbb{T}^3}} \asymp \|(e^W\Gamma)^{-1}u\|_{H^1_{\mathbb{T}^3}}, \tag{4.11}$$

where the second equivalence follows from the continuity of $\Gamma^{-1}: H^\beta_{\mathbb{T}^3} \to H^\beta_{\mathbb{T}^3}$ for $\beta < 3/2$, see Lemma 2.5. It is technically more convenient to consider the transformed equation for $u^\sharp = (e^W\Gamma)^{-1}u$ with initial datum $u_0^\sharp := (e^W\Gamma)^{-1}u_0 \in H^1_{\mathbb{T}^3}$:

$$u^\sharp(t) = e^{-it\mathcal{H}^\sharp}u_0^\sharp \mp i\int_0^t e^{-i(t-s)\mathcal{H}^\sharp}(e^W\Gamma)^{-1}[(e^W\Gamma u^\sharp(s))(|e^W\Gamma u^\sharp(s)|^2 * V_\beta)]ds$$
$$= e^{-it\mathcal{H}^\sharp}u_0^\sharp \mp i\int_0^t e^{-i(t-s)\mathcal{H}^\sharp}\Gamma^{-1}[\Gamma u^\sharp(s)(|e^W\Gamma u^\sharp(s)|^2 * V_\beta)]ds, \tag{4.12}$$

where we recall that $e^{-it\mathcal{H}^\sharp} = (e^W\Gamma)^{-1}e^{-it\mathcal{H}}(e^W\Gamma)$ and emphasize that there has been a multiplicative cancellation of $e^{\pm W}$ in the final step.

**Remark 4.2.** Since we solve in the natural domain of the energy of the equation (4.2) given by

$$\mathcal{E}_0(u) := -(u, \mathcal{H}u) \pm \frac{1}{p}\int |u|^2 V_\beta * |u|^2$$

which controls the $L^\infty H^1_{\mathbb{T}^3}$ norm of $u^\sharp$, proving local well-posedness allows us to immediately obtain global well-posedness as well as long as we choose the nonlinearity to be *defocussing*.

The well-posedness result in the energy domain reads as follows:

**Proposition 4.3.** *For $\beta > \frac{17}{20}$, assume that the interaction potential V is symmetric and*

$$V_\beta \in W^{\beta,1}(\mathbb{T}^3). \tag{4.13}$$

*For fixed enhancement $\Xi = \Xi_\omega$, (4.2) is locally well-posedness in the energy domain, in the sense that given any initial data $u_0 \in \mathcal{D}(\sqrt{-\mathcal{H}})$, there exists*

$$T = T(\Xi, \|\sqrt{-\mathcal{H}}u_0\|_{L^2}, \beta) \leq 1,$$

*such that there exists a unique solution $u \in C([0,T], \mathcal{D}(\sqrt{-\mathcal{H}})) \cap L^{\frac{10}{3}}([0,T], e^W W^{\frac{11}{20}-,\frac{10}{3}}_{\mathbb{T}^3})$ with the solution u depending continuously on the initial data.*

**Proof.** We consider the mild formulation (4.12) of the transformed equation (4.2) as the fixed point of the map

$$\mathcal{T}: B_\mathcal{Y}(M_R) \longrightarrow B_\mathcal{Y}(M_R), u^\sharp \mapsto \mathcal{T}(u^\sharp), \tag{4.14}$$

for any initial data $\|u_0^\sharp\|_{H^1_{\mathbb{T}^3}} \leq R$ for some given $R > 0$,

$$\mathcal{T}u^\sharp(t) := e^{-it\mathcal{H}^\sharp}u_0^\sharp \mp i\int_0^t e^{-i(t-s)\mathcal{H}^\sharp}\Gamma^{-1}[\Gamma u^\sharp(s)(|e^W\Gamma u^\sharp(s)|^2 * V_\beta)]ds$$

on the space

$$\mathcal{Y} := C_{[0,T]}H^1_{\mathbb{T}^3} \cap L^{\frac{10}{3}}_{[0,T]}W^{\frac{11}{20}-,\frac{10}{3}}_{\mathbb{T}^3},$$

where the maximal time of existence $T$ and the radius $M_R > 0$ will be determined depending on $R$ later. The continuity in time of $\mathcal{T}(u^\sharp)$ follows by Stone's theorem as in the proof of Proposition 4.1.

We now assume without loss of generality that $\beta = \frac{17}{20}+$, since it is easier to construct with larger $\beta \leq 3$ such that the convolution with $V_\beta$ gives more smoothing effect. By continuity of $\Gamma^{-1}$ and $e^{-it\mathcal{H}^\sharp}$ on $H^1_{\mathbb{T}^3}$, we have that

$$\|\mathcal{T}u^\sharp\|_{L^\infty_{[0,T]}H^1_{\mathbb{T}^3}} \lesssim \|u_0^\sharp\|_{H^1_{\mathbb{T}^3}} + \int_0^T \|\Gamma u^\sharp(s)\|_{H^1_{\mathbb{T}^3}}\||e^W\Gamma u^\sharp(s)|^2\|_{W^{-\frac{17}{20}-,\infty}_{\mathbb{T}^3}} + \|\Gamma u^\sharp(s)\|_{L^6_{\mathbb{T}^3}}\||e^W\Gamma u^\sharp(s)|^2\|_{W^{\frac{3}{20}-,3}_{\mathbb{T}^3}}ds$$

$$\lesssim \|u_0^\sharp\|_{H^1_{\mathbb{T}^3}} + \|\Gamma u^\sharp\|_{L^\infty_{[0,T]}H^1_{\mathbb{T}^3}}\int_0^T \left(\||e^W\Gamma u^\sharp(s)|^2\|_{W^{-\frac{1}{2}+,\frac{60}{7}-}_{\mathbb{T}^3}} + \||e^W\Gamma u^\sharp(s)|^2\|_{W^{\frac{3}{20}-,3}_{\mathbb{T}^3}}\right)ds$$

$$\lesssim \|u_0^\sharp\|_{H^1_{\mathbb{T}^3}} + \|\Gamma u^\sharp\|_{L^\infty_{[0,T]}H^1_{\mathbb{T}^3}}\int_0^T \left(\||\Gamma u_s^\sharp|^2(s)\|_{W^{\frac{11}{20}-,\frac{15}{7}}_{\mathbb{T}^3}} + \||\Gamma u_s^\sharp|^2(s)\|_{W^{\frac{3}{20}-,3}_{\mathbb{T}^3}}\right)ds$$

$$\lesssim \|u_0^\sharp\|_{H^1_{\mathbb{T}^3}} + \|\Gamma u^\sharp\|_{L^\infty_{[0,T]}H^1_{\mathbb{T}^3}}^2 T^{\frac{7}{10}}\|\Gamma u^\sharp\|_{L^{\frac{10}{3}}_{[0,T]}W^{\frac{11}{20}-,\frac{10}{3}}_{\mathbb{T}^3}}$$

where we used the fractional Leibniz rule in Lemma A.6 with the coefficients $\frac{1}{6} + \frac{1}{3} = \frac{1}{2} + \frac{1}{\infty} = \frac{1}{2}$ in the first inequality, the bound $\|e^{2W}\|_{W_{\mathbb{T}^3}^{1/2-,\infty}} < \infty$ and Sobolev embeddings

$$W_{\mathbb{T}^3}^{\frac{11}{20}-,\frac{15}{7}-} \hookrightarrow W_{\mathbb{T}^3}^{-\frac{1}{2}+,\frac{60}{7}-} \hookrightarrow W_{\mathbb{T}^3}^{\frac{17}{20}-,\infty}, \quad H_{\mathbb{T}^3}^1 \hookrightarrow L_{\mathbb{T}^3}^6$$

in the second and third inequalities, and the fractional Leibniz rule again with coefficients $\frac{7}{15} = \frac{1}{6} + \frac{3}{10}$ for the first term and $\frac{1}{6} + \frac{1}{6} = \frac{1}{3}$ for the second, the Sobolev embeddings

$$W_{\mathbb{T}^3}^{\frac{11}{20}-,\frac{10}{3}} \hookrightarrow W_{\mathbb{T}^3}^{\frac{3}{20}-,6};$$

while by the Strichartz-type of estimates for the Anderson Hamiltonian and the continuity of $\Gamma^{-1}$ on $H_{\mathbb{T}^3}^1$ again, we have as above that

$$\|\mathcal{T}u^\sharp\|_{L_{[0,T]}^{\frac{10}{3}} W_{\mathbb{T}^3}^{\frac{11}{20}-,\frac{10}{3}}} \lesssim \|u_0^\sharp\|_{H_{\mathbb{T}^3}^1} + \|\Gamma u^\sharp\|_{L_{[0,T]}^\infty H_{\mathbb{T}^3}^1}^2 T^{\frac{7}{10}} \|\Gamma u^\sharp\|_{L_{[0,T]}^{\frac{10}{3}} W_{\mathbb{T}^3}^{\frac{11}{20}-,\frac{10}{3}}}.$$

With the same argument as in the proof of Proposition 4.1, we can find $M_R > 0$ large enough and then $T > 0$ small enough such that (4.14) is a well-defined map. On the other hand, for any $u^{\sharp,1}, u^{\sharp,2} \in B_{\mathcal{Y}}(M_R)$ with the same inital datum, we have

$$\|\mathcal{T}u^{\sharp,1} - \mathcal{T}u^{\sharp,2}\|_{L_{[0,T]}^\infty H_{\mathbb{T}^3}^1} \lesssim \int_0^T \|\Gamma(u^{\sharp,1}(s) - u^{\sharp,2}(s))(|e^W \Gamma u^{\sharp,1}(s)|^2 * V_\beta)\|_{H_{\mathbb{T}^3}^1} ds$$

$$+ \int_0^T \|\Gamma u^{\sharp,2}(s)(|e^W \Gamma u^{\sharp,1}(s)|^2 * V_\beta - |e^W \Gamma u^{\sharp,2}(s)|^2 * V_\beta)\|_{H_{\mathbb{T}^3}^1} ds := (I) + (II),$$

where

$$(I) \lesssim \|\Gamma(u^{\sharp,1}(s) - u^{\sharp,2}(s))\|_{L_{[0,T]}^\infty H_{\mathbb{T}^3}^1} \int_0^T \left( \||\Gamma u^{\sharp,1}(s)|^2\|_{W_{\mathbb{T}^3}^{\frac{11}{20}-,\frac{15}{7}}} + \||\Gamma u^{\sharp,1}(s)|^2\|_{W_{\mathbb{T}^3}^{\frac{3}{20}-,3}} \right) ds$$

$$\lesssim T^{\frac{10}{3}} M_R^2 \|u^{\sharp,1} - u^{\sharp,2}\|_{L_{[0,T]}^\infty H_{\mathbb{T}^3}^1},$$

and

$$(II) \lesssim \int_0^T \|\Gamma u^{\sharp,2}(s)(|e^W \Gamma u^{\sharp,1}(s)|^2 * V_\beta - |e^W \Gamma u^{\sharp,2}(s)|^2 * V_\beta)\|_{H_{\mathbb{T}^3}^1} ds$$

$$\lesssim \|\Gamma u^{\sharp,1}\|_{L_{[0,T]}^\infty H_{\mathbb{T}^3}^1} \int_0^T \left( \||\Gamma(u^{\sharp,1}(s) - u^{\sharp,2}(s))|^2\|_{W_{\mathbb{T}^3}^{\frac{11}{20}-,\frac{15}{7}}} + \||\Gamma(u^{\sharp,1}(s) - u^{\sharp,2}(s))|^2\|_{W_{\mathbb{T}^3}^{\frac{3}{20}-,3}} \right) ds$$

$$\lesssim T^{\frac{10}{3}} M_R^2 \|u^{\sharp,1} - u^{\sharp,2}\|_{L_{[0,T]}^{\frac{10}{3}} W_{\mathbb{T}^3}^{\frac{11}{20}-,\frac{10}{3}}},$$

since

$$\int_0^T \||\Gamma(u^{\sharp,1}(s) - u^{\sharp,2}(s))|^2\|_{W_{\mathbb{T}^3}^{\frac{11}{20}-,\frac{15}{7}}} ds \lesssim \int_0^T \|\Gamma(u^{\sharp,1}(s) - u^{\sharp,2}(s))\|_{W_{\mathbb{T}^3}^{\frac{11}{20}-,\frac{10}{3}}} \|\Gamma(u^{\sharp,1}(s) - u^{\sharp,2}(s))\|_{L_{\mathbb{T}^3}^6} ds$$

$$\lesssim T^{\frac{7}{10}} \|u^{\sharp,1} - u^{\sharp,2}\|_{L_{[0,T]}^{\frac{10}{3}} W_{\mathbb{T}^3}^{\frac{11}{20}-,\frac{10}{3}}} \left( \|u^{\sharp,1}\|_{L_{[0,T]}^\infty L_{\mathbb{T}^3}^6} + \|u^{\sharp,2}\|_{L_{[0,T]}^\infty L_{\mathbb{T}^3}^6} \right).$$

The control of the difference in the space of $L_{[0,T]}^{\frac{10}{3}} W_{\mathbb{T}^3}^{\frac{11}{20}-,\frac{10}{3}}$ is identical, which indicates that for $T > 0$ small enough, $\mathcal{T}$ is contraction on $B_{\mathcal{Y}}(M_R)$. The local well-posedness of (4.2) follows then from the Banach fixed point. We can then translate the result on $u^\sharp$ to $u$ via Lemma 2.5. □

# 5 Global Well-posedness of stochastic Hartree NLS

## 5.1 Global Well-posedness in the Energy-Domain

We define the mass and energy functionals to be

$$\mathcal{M}[u] := \int_{\mathbb{T}^3} |u|^2, \quad \mathcal{E}_0[u] := \int_{\mathbb{T}^3} (-\mathcal{H}u)\bar{u} + \frac{1}{2}\int_{\mathbb{T}^3} |u|^2 |u|^2 * V_\beta.$$

From now on, we only consider the *defocusing* case of (4.2), that is

$$i\partial_t u = \mathcal{H}u - u(|u|^2 * V_\beta), \quad \text{on } [0,\infty) \times \mathbb{T}^3, \tag{5.1}$$

for which we have the conservation of mass

$$\frac{d}{dt}\mathcal{M}[u] = 2\Re\left(-i\int_{\mathbb{T}^3} u\overline{[\mathcal{H}u - u(|u|^2 * V_\beta)]}\right) = 0, \tag{5.2}$$

and (formally) conservation of energy

$$\begin{aligned}\frac{d}{dt}\mathcal{E}_0[u] &= 2\Re\int_{\mathbb{T}^3}(-\mathcal{H}u)\overline{\partial_t u} + \int_{\mathbb{T}^3}(\partial_t|u|^2)|u|^2 * V_\beta \\ &= 2\Re\left[i\int_{\mathbb{T}^3}(\mathcal{H}u)\overline{u(|u|^2 * V_\beta)}\right] - 2\Re\left[i\int_{\mathbb{T}^3}(\mathcal{H}u)\bar{u}|u|^2 * V_\beta\right] = 0.\end{aligned} \tag{5.3}$$

It then follows directly from Proposition 4.3 and the energy conservation (5.2) that we have the following well-posedness result in the energy domain $\mathcal{D}(\sqrt{-\mathcal{H}})$. We also claim that (5.1) is globally well-posed in $\mathcal{D}(\mathcal{H})$ under the assumption of Proposition 4.3, the discussion of which is postponed to the next section.

**Corollary 5.1.** *Under the assumption of Proposition 4.3, assuming in addition that the equation is defocusing, i.e. (5.1). Also we assume $V_\beta \geq 0$. Then (5.1) is globally well-posed in $\mathcal{D}(\sqrt{-\mathcal{H}})$.*

**Proof.** This is quite classical, as long as one can justify that the energy is conserved up to time $T$ given by Proposition 4.3, i.e. (5.2) holds since the time of existence therein can be chosen as $T \sim \|u_0\|_{\mathcal{D}(\sqrt{-\mathcal{H}})}^{-k}$ for some $k > 0$, and one has

$$\|u(T)\|_{\mathcal{D}(\sqrt{-\mathcal{H}})}^2 \lesssim \mathcal{E}_0[u(T)] = \mathcal{E}_0[u_0],$$

so one may restart the flow from $T$ to $2T$ and etc.

In order to justify the conservation of energy in the energy space (5.2), we can proceed as in Section 3.2.2 of [GUZ20]. Clearly (5.2) is rigorous in the "more regular" case $u_0 \in \mathcal{D}(\mathcal{H})$ and, as in [GUZ20], we can approximate the solution $u$ to (5.1) by the (global in time by Theorem 5.3 that we prove later) solution $u^\lambda$ to (5.1) with initial data

$$u_0^\lambda := (1 + \lambda\sqrt{-\mathcal{H}})^{-1} u_0,$$

for which we have

$$\mathcal{E}_0[u_0^\lambda] \to \mathcal{E}_0[u_0] \quad \text{as } \lambda \to 0,$$

and by the local well-posedness Theorem 4.3, especially the continuity of the solution map with respect to the initial conditions, the solutions converge in the sense that

$$u^\lambda \to u \quad \text{in } C([0,T]; \mathcal{D}(\sqrt{-\mathcal{H}})),$$

which further implies using the $\mathcal{D}(\sqrt{-\mathcal{H}}) \hookrightarrow L^6_{\mathbb{T}^3}$ embedding from Lemma 2.6 that the interaction part in $\mathcal{E}_0$ converges, and therefore

$$\mathcal{E}_0[u_0^\lambda] = \mathcal{E}_0[u^\lambda(t)] \to \mathcal{E}_0[u(t)] \quad \text{as } \lambda \to 0^+,$$

for all $t \leq T$. □

## 5.2 Global well-posedness in the Domain of Anderson Hamiltonian

We now consider the global well-posedness of (5.1) in the domain $\mathcal{D}(\mathcal{H})$ using the method of *modified energy* introduced in [PTV17] and applied for the multiplicative stochastic NLS on $\mathbb{T}^2$ and $\mathbb{R}^2$ in [TV23a],[TV23b] and [DLTV24] respectively. We first notice that it follows immediately from the Sobolev embedding, the energy conservation, Lemma 2.5 and the boundedness of $e^W$ that

$$\|u\|_{L^\infty_{[0,T]}L^6_{\mathbb{T}^3}} \lesssim \|u^\sharp\|_{L^\infty_{[0,T]}L^6_{\mathbb{T}^3}} \lesssim \|u^\sharp\|_{L^\infty_{[0,T]}H^1_{\mathbb{T}^3}} \lesssim \mathcal{E}_0^{\frac{1}{2}}[u_0], \tag{5.4}$$

which is uniform in time, while with the help of Strichartz estimates in Corollary 3.6, we expect to *bound higher spatial-integrability of $u^\sharp$ at the cost of some time integrability*. This is indeed possible thanks to the following *a priori estimates*.

**Lemma 5.2. (the "Nonlinear Strichartz" Estimate)** *For any solution $u^\sharp$ to (4.12) in the space of $H^2_{\mathbb{T}^3}$, we have for any parameter $\beta$ of the Hartree interaction such that*

$$\frac{7}{10} < \beta < 1,$$

*the bounds*

$$\|u^\sharp\|_{L^{\frac{10}{3}}_{[0,T]}W^{\frac{11}{20}-,\frac{10}{3}}_{\mathbb{T}^3}} \lesssim \|u^\sharp\|_{L^\infty_{[0,T]}H^2_{\mathbb{T}^3}}^{\frac{3}{7}\left(\frac{17}{20}-\beta\right)_+}, \tag{5.5}$$

*and*[5.1]

$$\|u^\sharp\|_{L^{\frac{10}{3}}_{[0,T]}W^{\frac{31}{20}-,\frac{10}{3}}_{\mathbb{T}^3}} \lesssim \|u^\sharp\|_{L^\infty_{[0,T]}H^2_{\mathbb{T}^3}}^{1+\frac{20}{7}\left(\frac{17}{20}-\beta\right)_+(1-\beta)}, \tag{5.6}$$

*where the constant only depends on the terminal time $T > 0$, $\|u_0^\sharp\|_{H^1_{\mathbb{T}^3}}$, the realization of the noise $\Xi$ and the norms of the potential $V_\beta$.*

Heuristically, (5.5) in the case $\beta \geq \frac{17}{20}$ is the a priori bounds applying Strichartz estimates in the mild formulation of energy solutions. In the case that $\beta < \frac{17}{20}$, it says how much of the $L^\infty_{[0,T]}H^2_{\mathbb{T}^3}$ one needs to bound the nonlinear term in the mild formulation. We postpone the proof to the end of the chapter in 5.2.2.

Let us firstly give a "naive" result coming from a Gronwall-type bound, i.e. how low one can choose the parameter $\beta$ if one simply wants a linear bound in the mild formulation leading to an exponential bound in time.

**Theorem 5.3. (Naive GWP in domain)** *Let $u_0 \in \mathcal{D}(\mathcal{H})$ and $T > 0$ and $V_\beta$ as above with $\beta > \frac{17}{20}$. Then there exists a unique mild solution to (5.1)*

$$u \in C([0,T]; \mathcal{D}(\mathcal{H})) \cap C^1([0,T]; L^2_{\mathbb{T}^3})$$

---

[5.1]. We use the notation $a_+ := \max\{a, 0\}$.

*with the solution depending continuously on the initial data.*

*Moreover, there are some constants $c(\mathcal{E}_0(u_0)), c'(\mathcal{E}_0(u_0)) > 0$ depending on the energy of the initial data $\mathcal{E}_0(u_0)$ s.t. we have the bound*

$$\|u\|_{C([0,T];\mathcal{D}(\mathcal{H}))\cap C^1([0,T];L^2_{\mathbb{T}^3})} \leq c(\mathcal{E}_0(u_0)) e^{c'(\mathcal{E}_0(u_0))T} \|u_0\|_{\mathcal{D}(\mathcal{H})}.$$

**Proof.** The main ingredient is an apriori estimate on the mild formulation (using the same integration by parts in time trick as in the proof of Proposition 4.1)

$$\mathcal{H}u(t) = e^{-it\mathcal{H}}\mathcal{H}u_0 - e^{-it\mathcal{H}}(u_0 V_\beta * |u_0|^2) + u(t) V_\beta * |u(t)|^2 - \int_0^t e^{-is\mathcal{H}} \partial_s(u(t-s) V_\beta * |u(t-s)|^2) ds,$$

combined with the local well-posedness from Section 4.1 gives GWP in a straightforward manner. We bound

$$\|\mathcal{H}u(t)\|_{L^2} \lesssim \|\mathcal{H}u_0\|_{L^2_{\mathbb{T}^3}}(1 + \mathcal{E}_0(u_0)) + \int_0^t \|\partial_s u\|_{L^\infty_T L^2_{\mathbb{T}^3}} \|V_\beta * |u(s)|^2\|_{L^\infty_{\mathbb{T}^3}} + \|u(s) V_\beta * [\partial_s u(s) \bar{u}(s)]\|_{L^2_{\mathbb{T}^3}} ds,$$

where we further bound the integral in time on the right by

$$\|\partial_s u\|_{L^\infty_{[0,t]} L^2_{\mathbb{T}^3}} \int_0^t \|u(s)\|^2_{L^{\frac{6}{\beta}+}_{\mathbb{T}^3}} \|V_\beta\|_{L^{\frac{3}{3-\beta}-}_{\mathbb{T}^3}} ds \lesssim \|V_\beta\|_{L^{\frac{3}{3-\beta}}_{\mathbb{T}^3}} \|\partial_s u\|_{L^\infty_{[0,t]} L^2_{\mathbb{T}^3}} t^{\frac{2}{5}} \|u\|^2_{L^{\frac{10}{3}}_{[0,T]} L^{\frac{6}{\beta}+}_{\mathbb{T}^3}}$$

$$\lesssim \|V_\beta\|_{L^{\frac{3}{3-\beta}}_{\mathbb{T}^3}} \|\partial_s u\|_{L^\infty_{[0,T]} L^2_{\mathbb{T}^3}} T^{\frac{2}{5}} \|u\|^2_{L^{\frac{10}{3}}_{[0,T]} W^{\frac{9}{10}-\frac{\beta}{2}+}_{\mathbb{T}^3}},$$

having used Hölder's/Young's inequality and the conservation of energy and the embeddings from Lemma 2.6. Using the equation, see also (5.8), we have

$$\|\partial_s u\|_{L^\infty_{[0,T]} L^2_{\mathbb{T}^3}} \lesssim \|\mathcal{H}u\|_{L^\infty_{[0,T]} L^2_{\mathbb{T}^3}} + C\mathcal{E}_0^{\frac{3}{2}}(u_0) \leq \|\mathcal{H}u\|_{L^\infty_{[0,T]} L^2_{\mathbb{T}^3}} + C\|\mathcal{H}u_0\|_{L^2_{\mathbb{T}^3}} \mathcal{E}_0(u_0)$$

and for the final term we use (5.5) which, since we take $\beta > \frac{17}{20}$ implying $\frac{9}{10} - \frac{\beta}{2} < \frac{19}{40}$, takes the form

$$\|u\|_{L^{\frac{10}{3}}_{[0,T]} W^{\frac{9}{10}-\frac{\beta}{2}+,\frac{10}{3}}_{\mathbb{T}^3}} \overset{\text{Lemma 2.6}}{\lesssim} \|u^\sharp\|_{L^{\frac{10}{3}}_{[0,T]} W^{\frac{9}{10}-\frac{\beta}{2}+,\frac{10}{3}}_{\mathbb{T}^3}} \overset{(5.5)}{\lesssim} \|u^\sharp\|_{L^{\frac{10}{3}}_{[0,T]} W^{\frac{11}{20}-,\frac{10}{3}}_{\mathbb{T}^3}} \lesssim 1. \tag{5.7}$$

This means in total we have

$$\|\mathcal{H}u\|_{L^\infty_{[0,t]} L^2_{\mathbb{T}^3}} \leq C(\mathcal{E}_0(u_0)) \|\mathcal{H}u_0\|_{L^2_{\mathbb{T}^3}} + t^{\frac{2}{5}} C'(\mathcal{E}_0(u_0)) \|\mathcal{H}u\|_{L^\infty_{[0,t]} L^2_{\mathbb{T}^3}},$$

and analogously for $0 < s < t < T$,

$$\|\mathcal{H}u\|_{L^\infty_{[0,t]} L^2_{\mathbb{T}^3}} \leq C(\mathcal{E}_0(u_0)) \|\mathcal{H}u\|_{L^\infty_{[0,s]} L^2_{\mathbb{T}^3}} + (t-s)^{\frac{2}{5}} C'(\mathcal{E}_0(u_0)) \|\mathcal{H}u\|_{L^\infty_{[0,t]} L^2_{\mathbb{T}^3}}.$$

We now set $f(t) := \|\mathcal{H}u\|_{L^\infty_{[0,t]} L^2_{\mathbb{T}^3}}$ and taking intervals $[s,t]$ of length $|t-s| < (C'(\mathcal{E}_0(u_0)))^{-\frac{5}{2}}$, this leads to

$$f(t) \leq [1 + C(\mathcal{E}_0(u_0))] f(s) \leq [1 + C(\mathcal{E}_0(u_0))]^{T(2C'(\mathcal{E}_0(u_0)))^{\frac{5}{2}}} f(0) \leq c(\mathcal{E}_0(u_0)) e^{c'(\mathcal{E}_0(u_0))T} f(0)$$

which is what we tried to prove, noting that because of the continuity in time of the solution in the domain, $f(0) = \lim_{t \to 0^+} \|\mathcal{H}u\|_{L^\infty_{[0,t]} L^2_{\mathbb{T}^3}} = \|\mathcal{H}u_0\|_{L^2_{\mathbb{T}^3}}$. □

**Remark 5.4.** Note that (5.7) is slightly improvable, since $\frac{9}{10} - \frac{\beta}{2} < \frac{19}{40} < \frac{11}{20}$, which is the regularity exponent appearing in the nonlinear Strichartz bound (5.5). This seems to imply that one can go further but choosing $\beta < \frac{17}{20}$ leads to a bound which is super-linear in $\|\mathcal{H}u\|_{L^\infty_{[0,t]}L^2_{\mathbb{T}^3}}$ hence the iteration argument does not work.

### 5.2.1 Second-Order Energy

It turns out that we can still push the $\beta$ slightly below the threshold of $17/20$ from Theorem 5.3 if we employ the method of second order energies, inspired by the series of works [PTV17],[TV23b],[TV23a] and [DLTV24]. We fix therefore in this section a solution $u \in C([0,T], \mathcal{D}(\mathcal{H}))$ to (4.2) with initial datum $u_0 \in \mathcal{D}(\mathcal{H})$ and denote

$$u_0^\sharp = (e^W \Gamma)^{-1} u_0 \in H^2_{\mathbb{T}^3} \quad \text{and} \quad u_0^\natural = \Lambda^{-1} u_0 \in H^2_{\mathbb{T}^3}.$$

As a first step to warm up, we show that for any $t \leq T$, the $L^2_{\mathbb{T}^3}$-norm of the time-derivative of $u$ is comparable to that of $\mathcal{H}u$ up to a lower-order term, that is

$$|\|\partial_t u(t)\|_{L^2_{\mathbb{T}^3}} - \|\mathcal{H}u(t)\|_{L^2_{\mathbb{T}^3}}| \lesssim_{V_\beta} (\mathcal{E}_0[u_0])^{\frac{3}{2}}. \tag{5.8}$$

Indeed, the bound (5.8) follows from (5.4) and the equation (5.1):

$$\|u(|u|^2 \ast V_\beta)\|_{L^2_{\mathbb{T}^3}} \lesssim_{V_\beta} \|u\|^3_{L^6_{\mathbb{T}^3}} \lesssim_\Xi (\mathcal{E}_0[u_0])^{\frac{3}{2}}. \tag{5.9}$$

Before stating the main theorem of this section, we derive the *second-order energy* we will consider. Following [PTV17], the preliminary ansatz is

$$\tilde{\mathcal{E}}_1[u(t)] = \|\partial_t u(t)\|^2_{L^2_{\mathbb{T}^3}}. \tag{5.10}$$

By (5.8) and Lemma 2.6, control of $\tilde{\mathcal{E}}_1$ gives a bound for $\|u^\sharp\|_{L^\infty H^2_{\mathbb{T}^3}}$. We also give useful but basic bounds in negative-regularity spaces, which will be frequently used to interpolate between the first- and second-order energies.

**Lemma 5.5. (Time derivative and interpolation)** *One has the following uniform in time bounds*

$$\|\partial_t u^\sharp(t)\|_{H^{-1}_{\mathbb{T}^3}} \asymp \|\partial_t u^\natural(t)\|_{H^{-1}_{\mathbb{T}^3}} \asymp \|e^{-W}\partial_t u(t)\|_{H^{-1}_{\mathbb{T}^3}} \asymp \|\partial_t u(t)\|_{(\mathcal{D}(\sqrt{-\mathcal{H}}))^*} \lesssim C(\mathcal{E}(u_0)), \tag{5.11}$$

*and*

$$\|\partial_t u^\natural\|_{H^{-\alpha}_{\mathbb{T}^3}} \lesssim \|\partial_t u^\natural\|^\alpha_{H^{-1}_{\mathbb{T}^3}} \|\partial_t u^\natural\|^{1-\alpha}_{L^2_{\mathbb{T}^3}} \lesssim C(\mathcal{E}(u_0)) (\tilde{\mathcal{E}}_1[u(t)])^{\frac{1-\alpha}{2}}, \quad \text{for } \alpha \in [0,1]. \tag{5.12}$$

**Proof.** The main point in (5.11) is that since we have identified[5.2] $\mathcal{D}(\sqrt{-\mathcal{H}})$ with $e^W H^1_{\mathbb{T}^3}$ one canonically has $\mathcal{D}(\sqrt{-\mathcal{H}}) \approx e^{-W} H^{-1}_{\mathbb{T}^3}$ which is precisely the final equivalence. This can be seen simply by noting

$$\|\varphi\|_{(\mathcal{D}(\sqrt{-\mathcal{H}}))^*} \asymp \sup_{\psi \in \mathcal{D}(\sqrt{-\mathcal{H}})} \langle \psi, \varphi \rangle_{\mathcal{D}(\sqrt{-\mathcal{H}}),(\mathcal{D}(\sqrt{-\mathcal{H}}))^*}$$

$$= \sup_{\psi' \in H^1_{\mathbb{T}^3}} \langle e^W \psi', \varphi \rangle_{\mathcal{D}(\sqrt{-\mathcal{H}}),(\mathcal{D}(\sqrt{-\mathcal{H}}))^*} \asymp \|e^W \varphi\|_{H^{-1}_{\mathbb{T}^3}}$$

---

[5.2]. c.f. Proposition 2.53 in [GUZ20].

The other equivalences follow from the boundedness of the maps $\Gamma$ and $\Theta$. The bound follows by using the equation, i.e.

$$\|\partial_t u(t)\|_{(\mathcal{D}(\sqrt{-\mathcal{H}}))^*} \leq \|\mathcal{H}u(t)\|_{(\mathcal{D}(\sqrt{-\mathcal{H}}))^*} + \|u(|u|^2 * V_\beta)\|_{(\mathcal{D}(\sqrt{-\mathcal{H}}))^*}$$
$$\overset{(5.9)}{\lesssim} \|u(t)\|_{\mathcal{D}(\sqrt{-\mathcal{H}})} + C(\mathcal{E}(u_0)) \lesssim C(\mathcal{E}(u_0)).$$

Finally, (5.11) follows by interpolation and (5.11). □

The strategy is to use the equation and some additional manipulations to write

$$\frac{d}{dt}\tilde{\mathcal{E}}_1 = \frac{d}{dt}A + O((\tilde{\mathcal{E}}_1)^{1-\kappa}) \text{ for some } \kappa > 0, \tag{5.13}$$

an explicit $A$ which will be included in the energy and a sublinear remainder, schematically denoted by $O((\tilde{\mathcal{E}}_1)^{1-\varepsilon})$. This will lead even to a polynomially growing bound in time for the second order energy.

In order to achieve this, we compute

$$\frac{d}{dt}\int_{\mathbb{T}^3}|\partial_t u|^2 = 2\Re\int_{\mathbb{T}^3}\overline{\partial_t u}\partial_t^2 u \overset{(5.1)}{=} 2\Re\left[-i\int_{\mathbb{T}^3}\overline{\partial_t u}\partial_t(\mathcal{H}u \pm u(|u|^2 * V_\beta))\right]$$
$$= \mp 2\Re\left[i\int_{\mathbb{T}^3}\overline{\partial_t u}u(\partial_t|u|^2 * V_\beta)\right] \overset{(5.1)}{=} \pm 2\Re\left[\int_{\mathbb{T}^3}\overline{\mathcal{H}u \pm u(|u|^2 * V_\beta)}u(\partial_t|u|^2 * V_\beta)\right]$$
$$= \pm 2\Re\left[\int_{\mathbb{T}^3}\overline{\mathcal{H}u}u(\partial_t|u|^2 * V_\beta)\right] + \frac{d}{dt}\int_{\mathbb{T}^3}|u|^2(|u|^2 * V_\beta)^2 - \int_{\mathbb{T}^3}(\partial_t|u|^2)(|u|^2 * V_\beta)^2,$$

having used the self-adjointness of $\mathcal{H}$ and the product rule to show that

$$\Re\left[-i\int_{\mathbb{T}^3}\overline{\partial_t u}\mathcal{H}\partial_t u\right] = \Re\left[-i\int_{\mathbb{T}^3}\overline{\partial_t u}\partial_t u|u|^2 * V_\beta\right] = 0$$

in the third equality.

The first term is still very badly behaved, and we want to manipulate it as in [PTV17], but since the operator $\mathcal{H}$ does not obey Leibnitz rule, one can not apply their method directly. We can, however, use the idea from [Zac19], which is to apply a suitable transformation to $\mathcal{H}u$ which makes it a perturbation of the Laplacian. In this setting this is precisely the $\Lambda$ transformation from Lemma 2.9. In fact, this yields

$$\Re\left[\int_{\mathbb{T}^3}\overline{\mathcal{H}u}u(\partial_t|u|^2 * V_\beta)\right]$$
$$= \Re\left[\int_{\mathbb{T}^3}\overline{\Delta u^\natural}u^\natural(\partial_t|u^\natural|^2 * V_\beta)\right] + \Re\left[\int_{\mathbb{T}^3}\overline{(\mathcal{H}\Lambda - \Delta)u^\natural}u^\natural(\partial_t|u^\natural|^2 * V_\beta)\right]$$
$$+ \Re\left[\int_{\mathbb{T}^3}\overline{\mathcal{H}u}(\Lambda - \mathrm{id})u^\natural(\partial_t|u^\natural|^2 * V_\beta)\right] + \Re\left[\int_{\mathbb{T}^3}\overline{\mathcal{H}u}u(\partial_t(|u|^2 - |u^\natural|^2) * V_\beta)\right]$$
$$=: \text{(I)} + \text{(II)} + \text{(III)} + \text{(IV)}, \tag{5.14}$$

where the first term can be further written as (again, inspired by [PTV17])

$$\text{(I)} = \int_{\mathbb{T}^3}\Delta|u^\natural|^2(\partial_t|u^\natural|^2 * V_\beta) - 2\int_{\mathbb{T}^3}|\nabla u^\natural|^2(\partial_t|u^\natural|^2 * V_\beta)$$
$$= -\frac{1}{2}\frac{d}{dt}\int_{\mathbb{T}^3}\nabla|u^\natural|^2(\nabla|u^\natural|^2 * V_\beta) - 2\int_{\mathbb{T}^3}|\nabla u^\natural|^2(\partial_t|u^\natural|^2 * V_\beta). \tag{5.15}$$

**Remark 5.6.** We notice that in the case of $3 \geq \beta \geq 1$, the transformation (5.14) is no longer necessary, since we can instead bound

$$\left|\Re\left[\int_{\mathbb{T}^3}\overline{\mathcal{H}u}u(\partial_t|u|^2 \ast V_\beta)\right]\right|(t) \lesssim \|V_\beta\|_{L^{\frac{3}{2}}_{\mathbb{T}^3}} \|\partial_t|u|^2\|_{L^{\frac{3}{2}}_{\mathbb{T}^3}} \|\mathcal{H}u\|_{L^2_{\mathbb{T}^3}} \|u\|_{L^6_{\mathbb{T}^3}}$$

$$\lesssim \|V_\beta\|_{L^{\frac{3}{3-\beta}}_{\mathbb{T}^3}} \left(1 + \|\mathcal{H}u\|^2_{L^2_{\mathbb{T}^3}}\right)\|u\|^2_{L^6_{\mathbb{T}^3}}$$

$$\lesssim C(\mathcal{E}_0[u_0])\|V_\beta\|_{L^{\frac{3}{3-\beta}}_{\mathbb{T}^3}}(1 + \tilde{\mathcal{E}}_1[u(t)]),$$

due to (5.8). Such a bound suffices to close the estimate in Proposition 5.11 with Grönwall inequality and obtain an exponential bound, similar to Theorem 5.3. We will therefore only focus on the case of $0 \leq \beta < 1$ from now on.

**Remark 5.7.** If one did not have the map $\Theta$, but simply the exponential and paracontrolled maps, the differences $\Delta - e^{-W}\mathcal{H}e^W\Gamma$ and $\Gamma - \mathrm{id}$ would appear. Indeed, one would get

$$\int_{\mathbb{T}^3}\overline{\mathcal{H}u}u(\partial_t|u|^2 \ast V_\beta) = \int_{\mathbb{T}^3}\overline{e^{-W}\mathcal{H}e^W\Gamma u^\sharp}\Gamma u^\sharp(\partial_t|u|^2 \ast V_\beta)e^{2W},$$

but then one would instead have the leading term in (5.14) after the replacement of $e^{-W}\mathcal{H}e^W\Gamma$ with $\Delta$ and $\Gamma$ with id respectively as

$$\int_{\mathbb{T}^3}\Delta|u^\sharp|^2(\partial_t|u^\sharp|^2 \ast V_\beta)e^{2W} = \frac{1}{2}\frac{\mathrm{d}}{\mathrm{d}t}\int_{\mathbb{T}^3}\nabla|u^\sharp|^2 \cdot (\nabla|u^\sharp|^2 \ast V_\beta)e^{2W} - 2\int_{\mathbb{T}^3}\nabla|u^\sharp|^2 \cdot \nabla W(\partial_t|u^\sharp|^2 \ast V_\beta)e^{2W},$$

which is much worse, since $\nabla W \in C^{-1/2-}_{\mathbb{T}^3}$ only. We have therefore introduced the "natural" transformation $\Theta$, which absorbs the exponential ansatz see Lemma 2.8.

We can thus set the second order energy functional of (4.2) to be

$$\mathcal{E}_1[u] := \int_{\mathbb{T}^3}|\partial_t u|^2 - \int_{\mathbb{T}^3}|u|^2(|u|^2 \ast V_\beta)^2 - \int_{\mathbb{T}^3}\nabla|u^\sharp|^2(\nabla|u^\sharp|^2 \ast V_\beta). \tag{5.16}$$

Before we prove that $\frac{\mathrm{d}}{\mathrm{d}t}\mathcal{E}_1[u] \lesssim 1 + O(\mathcal{E}_1[u])^{1-\kappa}$ for some $0 < \kappa < 1$ (after shifting the $\mathcal{E}_1$ to make it non-negative), we indicate that the other terms in $\mathcal{E}_1$ are lower-order with respect to the first term.

**Lemma 5.8.** *One has the bound for any $\beta > 0$ that*

$$\left|\int_{\mathbb{T}^3}\nabla|u^\sharp|^2(\nabla|u^\sharp|^2 \ast V_\beta)\right| + \int_{\mathbb{T}^3}|u|^2(|u|^2 \ast V_\beta)^2 \leq \frac{1}{4}\int_{\mathbb{T}^3}|\partial_t u|^2 + C(\mathcal{E}_0(u_0)),$$

*which implies*

$$|\mathcal{E}_1[u](t) - \tilde{\mathcal{E}}_1[u](t)| \leq \frac{1}{2}\tilde{\mathcal{E}}_1[u](t) + C(\mathcal{E}_0(u_0)),$$

*and*

$$\tilde{\mathcal{E}}_1[u](t) \leq 2\mathcal{E}_1[u](t) + C(\mathcal{E}_0(u_0)),$$
$$|\mathcal{E}_1[u](t)| \leq 2\tilde{\mathcal{E}}_1[u](t) + C(\mathcal{E}_0(u_0)).$$

**Proof.** This follows quite straightforwardly by Lemma 2.6 and (5.8) and the standard combination of Sobolev interpolation and Young's inequality. Indeed,

$$\|\nabla|u^\natural|^2\|_{L^p_{\mathbb{T}^3}} \lesssim \|u^\natural\|_{H^1_{\mathbb{T}^3}}\|u^\natural\|_{L^{\frac{2p}{2-p}}_{\mathbb{T}^3}} \lesssim \|u^\natural\|_{H^1_{\mathbb{T}^3}}^{\frac{3}{p}}\|u^\natural\|_{H^2_{\mathbb{T}^3}}^{2-\frac{3}{p}}, \quad \forall p \in \left[\frac{3}{2}, 2\right],$$

and therefore for $\beta \in \left(0, \frac{1}{2}\right]$,

$$\left|\int_{\mathbb{T}^3} \nabla|u^\natural|^2 (\nabla|u^\natural|^2 \star V_\beta)\right| \lesssim_{V_\beta} \|\nabla|u^\natural|^2\|_{L^{\frac{6}{3+2\beta}}_{\mathbb{T}^3}} \|\nabla|u^\natural|^2\|_{L^2_{\mathbb{T}^3}} \lesssim \|u^\natural\|_{H^1_{\mathbb{T}^3}}^{3+\beta}\|u^\natural\|_{H^2_{\mathbb{T}^3}}^{1-\beta},$$

which gives the desired bound by Young's inequality since $1-\beta < 1$. The bound for $\beta \geq \frac{1}{2}$ for $\int_{\mathbb{T}^3} |u|^2 (|u|^2 \star V_\beta)^2$ is easier.

$\square$

Now we show that the terms on the right hand side of (5.14) are indeed of the form schematically given by (5.13).

**Proposition 5.9.** *We have the following bounds for the replacement in (5.14) in terms of the higher-order energy, that is*

$$|(\mathrm{II})| + |(\mathrm{III})| + |(\mathrm{IV})| \leq C(\mathcal{E}_0(u_0))(\mathcal{E}_1[u])^{1-\kappa}. \tag{5.17}$$

*for a suitable $\kappa > 0$.*

*Furthermore, the second term in* (I) *can be bounded by*

$$\left|2\int_0^T \int_{\mathbb{T}^3} |\nabla u^\natural|^2 (\partial_t |u^\natural|^2 \star V_\beta)\,\mathrm{d}t\right| \lesssim \|u^\sharp\|_{L^\infty_{[0,T]}H^2_{\mathbb{T}^3}}^{s_\beta}$$

*where*

$$s_\beta = -\frac{20}{7}\beta^3 + \frac{65}{7}\beta^2 - \frac{1663}{140}\beta + \frac{19159}{2800}. \tag{5.18}$$

**Remark 5.10.** It suffices to have $s_\beta < 2$ to close the bound for the higher-order energy $\mathcal{E}_1[u]$, which gives approximately

$$0.733 < \beta < \frac{17}{20}. \tag{5.19}$$

**Proof.** For the first inequality we bound using Hölder's and Young's convolution inequalities using that $\beta > \frac{1}{2}$ implies $V_\beta \in L^{\frac{6}{5}+}$ and the usual embedding results from Lemma 2.6

$$\begin{aligned}
|(\mathrm{II})| &\leq \|(\mathcal{H}\exp(W^M)\Gamma\Theta^{-1} - \Delta)u^\natural\|_{L^{\frac{3}{1+\varepsilon}}_{\mathbb{T}^3}} \|u^\natural\|_{L^6_{\mathbb{T}^3}} \|\partial_t|u^\natural|^2\|_{L^{\frac{3}{2}}_{\mathbb{T}^3}} \|V_\beta\|_{L^{\frac{6}{5-\varepsilon}}_{\mathbb{T}^3}} \\
&\stackrel{(2.27)}{\lesssim} \mathcal{E}_0^{\frac{1}{2}}[u_0] \|u^\natural\|_{W^{2,\frac{3}{1+\varepsilon}}_{\mathbb{T}^3}^{\frac{3}{2}+\delta}} \|u^\natural\|_{L^6} \|\partial_t u^\natural\|_{L^2_{\mathbb{T}^3}} \|V_\beta\|_{L^{\frac{6}{5-\varepsilon}}_{\mathbb{T}^3}} \\
&\lesssim \mathcal{E}_1^{\frac{1}{2}}[u]\mathcal{E}_0[u_0]\|u^\natural\|_{H^{2+\delta-\varepsilon}_{\mathbb{T}^3}}\|V_\beta\|_{L^{\frac{6}{5-\varepsilon}}_{\mathbb{T}^3}} \\
&\lesssim C(\mathcal{E}_0[u_0])\|u^\natural\|_{H^2_{\mathbb{T}^3}}^{1-2\kappa}\mathcal{E}_1^{\frac{1}{2}}[u] \lesssim C(\mathcal{E}_0[u_0])\mathcal{E}_1^{1-\kappa}[u]
\end{aligned}$$

for $\varepsilon > \delta > 0$ small enough and $\kappa = \frac{1}{2}(\varepsilon - \delta) > 0$; while for the second one, we have similarly

$$\begin{aligned}
|(\text{III})| &\leq \|\mathcal{H}u\|_{L^2}\|(\Lambda - \text{id})u^\natural\|_{L^{\frac{2}{\varepsilon}}_{\mathbb{T}^3}}\|V_\beta * \partial_t |u^\natural|^2\|_{L^{\frac{2}{1-\varepsilon}}_{\mathbb{T}^3}} \\
&\lesssim \mathcal{E}_1[u]^{\frac{1}{2}}\|(\Lambda - \text{id})u^\natural\|_{W^{\frac{1-3\varepsilon}{2},6}_{\mathbb{T}^3}}\|\partial_t|u^\natural|^2\|_{W^{-\delta,\frac{3}{2+\varepsilon}}_{\mathbb{T}^3}}\|V_\beta\|_{W^{\delta,\frac{6}{5-5\varepsilon}}} \\
&\stackrel{(2.26)}{\lesssim} \mathcal{E}_1[u]^{\frac{1}{2}}\mathcal{E}_0^{\frac{1}{2}}[u_0]\|\partial_t u^\natural\|_{H^{-\delta}_{\mathbb{T}^3}}\|u^\natural\|_{W^{2\delta,\frac{6}{1+2\varepsilon}}_{\mathbb{T}^3}} \\
&\stackrel{(5.11)}{\lesssim} \mathcal{E}_1[u]^{1-\kappa}C(\mathcal{E}_0[u_0])
\end{aligned}$$

where $\varepsilon, \delta > 0$ are small enough and $\kappa = \frac{\delta}{2}$; moreover we have

$$\begin{aligned}
|(\text{IV})| &\leq \|\mathcal{H}u\|_{L^2_{\mathbb{T}^3}}\|u\|_{L^6_{\mathbb{T}^3}}\|V_\beta * \partial_t(|u|^2 - |u^\natural|^2)\|_{L^3_{\mathbb{T}^3}} \\
&\lesssim \mathcal{E}_0^{\frac{1}{2}}[u_0]\mathcal{E}_1^{\frac{1}{2}}[u]\|V_\beta\|_{W^{\delta,\frac{6}{5-3\varepsilon}}_{\mathbb{T}^3}}\|\partial_t(|u|^2 - |u^\natural|^2)\|_{W^{-\delta,\frac{2}{1+\varepsilon}}_{\mathbb{T}^3}} \\
&\lesssim \mathcal{E}_0^{\frac{1}{2}}[u_0]\mathcal{E}_1^{\frac{1}{2}}[u]\left(\|(1-\Lambda)\partial_t u\|_{W^{-\delta,3}_{\mathbb{T}^3}}\|u\|_{W^{2\delta,\frac{6}{1+3\varepsilon}}_{\mathbb{T}^3}} + \|(1-\Lambda)u\|_{W^{\delta,\frac{2}{\varepsilon}}_{\mathbb{T}^3}}\|\partial_t u\|_{H^{-\delta}_{\mathbb{T}^3}}\right) \\
&\lesssim \mathcal{E}_0^{\frac{1}{2}}[u_0]\mathcal{E}_1^{\frac{1}{2}}[u]\left(\|(1-\Lambda)\partial_t u\|_{H^{\frac{1}{2}-\delta}_{\mathbb{T}^3}}\|u\|_{W^{2\delta,\frac{6}{1+3\varepsilon}}_{\mathbb{T}^3}} + \|(1-\Lambda)u\|_{W^{\frac{1}{2}-,6}_{\mathbb{T}^3}}\|\partial_t u\|_{H^{-\delta}_{\mathbb{T}^3}}\right) \\
&\stackrel{(2.26)}{\lesssim} \mathcal{E}_0[u_0]\mathcal{E}_1^{\frac{1}{2}}[u]\|\partial_t u\|_{H^{-\frac{\delta}{2}}_{\mathbb{T}^3}} \stackrel{(5.11)}{\lesssim} \mathcal{E}_1[u]^{1-\kappa}C(\mathcal{E}_0[u_0]),
\end{aligned}$$

for $\varepsilon, \delta > 0$ sufficiently small and $\kappa = \frac{\delta}{4}$. The $\kappa > 0$ in (5.17) is then chosen as the minimum of the $\kappa$ parameters in the previous bounds.

Lastly, we have for the final term which is the most irregular,

$$\begin{aligned}
\left|\int_0^T \int_{\mathbb{T}^3} |\nabla u^\natural|^2 (\partial_t |u^\natural|^2 * V_\beta) \, dt\right| &\lesssim T^{\frac{1}{10}}\|V_\beta\|_{L^{\frac{3}{3-\beta}}_{\mathbb{T}^3}}\|\partial_t|u^\natural|^2\|_{L^{\frac{10}{3}}_{[0,T]}L^{\frac{60}{37}}_{\mathbb{T}^3}}\|\,|\nabla u^\natural|^2\|_{L^{\frac{5}{3}}_{[0,T]}L^{\frac{60}{23+20\beta}+}_{\mathbb{T}^3}} \\
&\lesssim_{V_\beta} T^{\frac{1}{10}}\|u^\natural\|_{L^{\frac{10}{3}}_{[0,T]}L^{\frac{60}{7}}_{\mathbb{T}^3}}\|\partial_t u^\natural\|_{L^\infty_{[0,T]}L^2_{\mathbb{T}^3}}\|\nabla u^\natural\|^2_{L^{\frac{10}{3}}_{[0,T]}L^{\frac{120}{23+20\beta}+}_{\mathbb{T}^3}}
\end{aligned}$$

having used Hölder's inequality in both space and time.

We assume from now on $\frac{17}{20} > \beta > \frac{13}{20}$ (the lower bound is equivalent to $p_\beta := \frac{120}{23+20\beta} + < \frac{10}{3}$) and bound

$$\|u^\natural(s)\|_{W^{1,p_\beta}_{\mathbb{T}^3}} \lesssim \|u^\natural(s)\|^{\frac{9}{40}+\frac{\beta}{2}}_{W^{\frac{9}{40}+\frac{\beta}{2},p_\beta}_{\mathbb{T}^3}}\|u^\natural(s)\|^{\frac{31}{40}-\frac{\beta}{2}}_{W^{\frac{49}{40}+\frac{\beta}{2},p_\beta}_{\mathbb{T}^3}} \lesssim \|u^\natural(s)\|^{\frac{9}{40}+\frac{\beta}{2}}_{W^{\frac{11}{20},\frac{10}{3}}_{\mathbb{T}^3}}\|u^\natural(s)\|^{\frac{31}{40}-\frac{\beta}{2}}_{W^{\frac{31}{20},\frac{10}{3}}_{\mathbb{T}^3}},$$

and applying Lemma 5.2, this indicates that

$$\begin{aligned}
\left|\int_0^T \int_{\mathbb{T}^3} |\nabla u^\natural|^2 (\partial_t |u^\natural|^2 * V_\beta) \, dt\right| &\lesssim T^{\frac{1}{10}}\|u^\natural\|_{L^{\frac{10}{3}}_{[0,T]}L^{\frac{60}{7}-}_{\mathbb{T}^3}}\|\partial_t u^\natural\|_{L^\infty_{[0,T]}L^2_{\mathbb{T}^3}}\|u^\natural\|^{\frac{9}{20}+\beta}_{L^{\frac{10}{3}}_{[0,T]}W^{\frac{11}{20}-,\frac{10}{3}}_{\mathbb{T}^3}}\|u^\natural\|^{\frac{31}{20}-\beta}_{L^{\frac{10}{3}}_{[0,T]}W^{\frac{31}{20}-,\frac{10}{3}}_{\mathbb{T}^3}} \\
&\lesssim T^{\frac{1}{10}}\|u^\natural\|^{\frac{29}{20}+\beta}_{L^{\frac{10}{3}}_{[0,T]}W^{\frac{11}{20}-,\frac{10}{3}}_{\mathbb{T}^3}}\|\partial_t u^\natural\|_{L^\infty_{[0,T]}L^2_{\mathbb{T}^3}}\|u^\natural\|^{\frac{31}{20}-\beta}_{L^{\frac{10}{3}}_{[0,T]}W^{\frac{31}{20}-,\frac{10}{3}}_{\mathbb{T}^3}} \lesssim T^{\frac{1}{10}}\|u^\natural\|^{s_\beta}_{L^\infty_{[0,T]}H^2_{\mathbb{T}^3}},
\end{aligned}$$

where $s_\beta$ is defined by (5.18). □

**Proposition 5.11.** *For the second-order energy defined in (5.16), we have for any solution $u \in C([0, T], \mathcal{D}(\mathcal{H}))$ to (4.2) that*

$$\mathcal{E}_1[u](t) - \mathcal{E}_1[u](s) \leq C(\mathcal{E}_0(u_0))(t-s)^{\frac{1}{10}}(1 + (\mathcal{E}_1[u](t))^{1-\kappa}) \tag{5.20}$$

*for $|t - s| < 1$ and some $0 < \kappa < 1$.*

*This leads to an explicit growth for $\mathcal{E}_1$ namely*

$$\sup_{0 \leq \tau \leq T} (\mathcal{E}_1[u])(\tau) \lesssim (\mathcal{E}_1[u])(0) + C(\mathcal{E}_0(u_0)) + T^{\frac{1}{\kappa}} \text{ for any } T \geq 1 \tag{5.21}$$

*i.e. we have a polynomial growth in time of the $\mathcal{D}(\mathcal{H})$ norm of the solution u.*

**Proof.** By an approximation argument, we can make the previous computations rigorous for the second order energy $\mathcal{E}_1$ of the solution of the 3d stochastic Hartree NLS (5.1) in an integrated form which by Proposition 5.9 is precisely (5.20), choosing the small parameter $\kappa$ s.t. $2 - 2\kappa \leq s_\beta$, see (5.18). Note that the restriction $|t - s| < 1$ comes from the fact that the Strichartz estimates hold only on intervals of length $\sim 1$.

The norm growth bound (5.21) follows in a straightforward way from (5.20). Indeed, fixing $T > 0$, we simply estimate

$$\mathcal{E}_1[u](\tau) = \mathcal{E}_1[u](0) + \sum_{i=1}^{\lceil \tau \rceil}\left(\mathcal{E}_1[u]\left(\tau\frac{i}{\lceil \tau \rceil}\right) - \mathcal{E}_1[u]\left(\tau\frac{i-1}{\lceil \tau \rceil}\right)\right)$$

$$\Rightarrow$$

$$|\mathcal{E}_1[u](\tau)| \stackrel{(5.20)}{\leq} |\mathcal{E}_1[u](0)| + C(\mathcal{E}_0(u_0))\left(\lceil \tau \rceil + \sum_{i=1}^{\lceil \tau \rceil}\left(\mathcal{E}_1[u]\left(\tau\frac{i}{\lceil \tau \rceil}\right)\right)^{1-\kappa}\right)$$

$$\leq |\mathcal{E}_1[u](0)| + C(\mathcal{E}_0(u_0))\left(\lceil T \rceil + \sup_{0 \leq \tau \leq T}(\mathcal{E}_1[u](\tau))^{1-\kappa}T\right)$$

which implies

$$\sup_{0 \leq \tau \leq T}|\mathcal{E}_1[u](\tau)| \leq |\mathcal{E}_1[u](0)| + C(\mathcal{E}_0(u_0))\left(\lceil T \rceil + T^{\frac{1-\kappa}{\kappa}}\right) + \frac{1}{2}\sup_{0 \leq \tau \leq T}(\mathcal{E}_1[u](\tau))$$

$$\Rightarrow$$

$$\sup_{0 \leq \tau \leq T}|\mathcal{E}_1[u](\tau)| \leq 2|\mathcal{E}_1[u](0)| + C(\mathcal{E}_0(u_0))\left(\lceil T \rceil + T^{\frac{1}{\kappa}}\right).$$

this finishes the proof. □

The main theorem of this chapter follows from the above considerations and (4.6).

**Theorem 5.12.** *Let $u_0 \in \mathcal{D}(\mathcal{H})$, $T > 0$ and $V_\beta \in W^{\beta,1}_{\mathbb{T}^3}$ non-negative and symmetric with $\beta < \frac{17}{20}$ be such that $s_\beta < 2$, where $s_\beta$ is defined in (5.18). Then there exists a unique mild solution to (5.1)*

$$u \in C([0, T]; \mathcal{D}(\mathcal{H})) \cap C^1([0, T]; L^2_{\mathbb{T}^3})$$

*with the solution depending continuously on the initial data and its norm growing algebraically in $T$.*

### 5.2.2 Proof of Lemma 5.2

In the end of the chapter, we present the proof of the bounds in Lemma 5.2, which have been essential for the improvement of the choice of smaller $\beta > 0$ for the well-posedness of (5.1).

**Proof. (Lemma 5.2)** We have from the mild formulation and the homogeneous and inhomogenuous Strichartz-type estimates from Corollary 3.6 and Lemma 2.5 that for the time interval $[0, T]$ with $T \sim 1$ and any $0 < T_0 \leq T$ one has

$$\|u^\#\|_{L^{\frac{10}{3}}_{[0,T_0]} W^{\frac{11}{20},\frac{10}{3}}_{\mathbb{T}^3}} \lesssim \|u_0^\#\|_{H^1_{\mathbb{T}^3}} + \left\|\int_0^t e^{-i(t-s)\mathcal{H}^\#}\Gamma^{-1}[\Gamma u^\#(s)(|e^W \Gamma u^\#(s)|^2 * V_\beta)]ds\right\|_{L^{\frac{10}{3}}_{t;[0,T_0]} W^{\frac{11}{20},\frac{10}{3}}_{\mathbb{T}^3}}$$

$$\lesssim \|u_0^\#\|_{H^1_{\mathbb{T}^3}} + \int_0^{T_0} \|\Gamma u^\#(s)\|_{H^1_{\mathbb{T}^3}} \||e^W \Gamma u^\#(s)|^2 * V_\beta\|_{L^\infty_{\mathbb{T}^3}} ds$$

$$+ \int_0^{T_0} \|\Gamma u^\#(s)\|_{L^{\frac{60}{7}-}_{\mathbb{T}^3}} \||e^W \Gamma u^\#(s)|^2 * V_\beta\|_{W^{1,\frac{60}{23}+}_{\mathbb{T}^3}} ds,$$

where we have also applied the fractional Leibniz rule in the last inequality. We can further bound

$$\||e^W \Gamma u^\#(s)|^2 * V_\beta\|_{L^\infty_{\mathbb{T}^3}} \lesssim \||\Gamma u^\#(s)|^2\|_{L^{\frac{3}{\beta}+}_{\mathbb{T}^3}} \|V_\beta\|_{L^{\frac{3}{3-\beta}-}_{\mathbb{T}^3}}$$

$$\lesssim_{V_\beta} \|u^\#(s)\|_{L^{\frac{60}{7}}_{\mathbb{T}^3}} \|u^\#(s)\|_{L^{\frac{60}{20\beta-7}+}_{\mathbb{T}^3}}$$

$$\lesssim_{V_\beta} \|u^\#(s)\|_{W^{\frac{11}{20},\frac{10}{3}}_{\mathbb{T}^3}} \|u^\#(s)\|_{H^1_{\mathbb{T}^3}}^{\frac{3}{20}+\beta-} \|u(s)\|_{H^2_{\mathbb{T}^3}}^{\frac{17}{20}-\beta+},$$

where we have used the choice that $\beta \leq \frac{17}{20}$ and

$$\|u^\#(s)\|_{L^{\frac{60}{20\beta-7}+}_{\mathbb{T}^3}} \lesssim \|u^\#(s)\|_{H^{\frac{37}{20}-\beta+}_{\mathbb{T}^3}} \lesssim \|u^\#(s)\|_{H^1_{\mathbb{T}^3}}^{\frac{3}{20}+\beta-} \|u^\#(s)\|_{H^2_{\mathbb{T}^3}}^{\frac{17}{20}-\beta+},$$

by interpolation. Besides, we have by fractional Leibniz again and Sobolev embedding that

$$\||e^W \Gamma u^\#(s)|^2 * V_\beta\|_{W^{1,\frac{60}{23}+}_{\mathbb{T}^3}} \lesssim \||e^W \Gamma u^\#(s)|^2\|_{W^{\frac{1}{2},\frac{60}{20\beta+13}+}_{\mathbb{T}^3}} \|V_\beta\|_{W^{\frac{1}{2}+,\frac{6}{7-2\beta}-}_{\mathbb{T}^3}}$$

$$\lesssim \||\Gamma u^\#(s)|^2\|_{W^{\frac{1}{2}-,\frac{60}{20\beta+13}}_{\mathbb{T}^3}} \|V_\beta\|_{W^{\frac{1}{2}+,\frac{6}{7-2\beta}-}_{\mathbb{T}^3}}$$

$$\lesssim \|u^\#(s)\|_{L^6_{\mathbb{T}^3}} \|u^\#(s)\|_{W^{\frac{1}{2}-,\frac{60}{20\beta+3}}_{\mathbb{T}^3}} \|V_\beta\|_{W^{\frac{1}{2}+,\frac{6}{7-2\beta}-}_{\mathbb{T}^3}}$$

$$\lesssim \|u(s)\|_{H^1_{\mathbb{T}^3}}^{\frac{3}{20}+\beta-} \|u(s)\|_{H^2_{\mathbb{T}^3}}^{\frac{17}{20}-\beta+} \|V_\beta\|_{W^{\frac{1}{2}+,\frac{6}{7-2\beta}-}_{\mathbb{T}^3}},$$

since $\frac{1}{2} < \beta < \frac{17}{20}$, and also $\|u^\#(s)\|_{W^{\frac{1}{2}-,\frac{60}{20\beta+3}}_{\mathbb{T}^3}} \lesssim \|u^\#(s)\|_{H^{\frac{37}{20}-\beta+}_{\mathbb{T}^3}}$ as above.

Combining the bounds, we have derived that for any $0 < T_0 \leq T$ that

$$\|u^\#\|_{L^{\frac{10}{3}}_{[0,T_0]} W^{\frac{11}{20},\frac{10}{3}}_{\mathbb{T}^3}} \lesssim \|u_0^\#\|_{H^1_{\mathbb{T}^3}} + T_0^{\frac{7}{10}} \|u^\#\|_{L^{\frac{10}{3}}_{[0,T_0]} W^{\frac{11}{20},\frac{10}{3}}_{\mathbb{T}^3}} \|u^\#\|_{L^\infty_{[0,T]} H^2_{\mathbb{T}^3}}^{\frac{17}{20}-\beta+},$$

where the constant only depends on the constant in the linear Strichartz estimate (3.12), the norms of $V_\beta$ and $\|u^\sharp\|_{L^\infty_{[0,T]}H^1_{\mathbb{T}^3}} \lesssim \|u^\sharp_0\|_{H^1_{\mathbb{T}^3}}$ but not on $T_0$. We can therefore assume without loss of generality that $\|u^\sharp\|_{L^\infty_{[0,T]}H^2_{\mathbb{T}^3}} > 0$ and choose

$$T_0^{-\frac{7}{10}} \sim \|u^\sharp\|_{L^\infty_{[0,T]}H^2_{\mathbb{T}^3}}^{\frac{17}{20} - \beta+} \tag{5.22}$$

small enough which gives rise to a partition $\{I\}$ of $[0,T]$ such that the mesh of $|I| \sim T_0$ such that

$$\|u^\sharp\|_{L^{\frac{10}{3}}_{[0,T]} W^{\frac{11}{20},\frac{10}{3}}_{\mathbb{T}^3}} = \left(\sum_I \|u^\sharp\|_{L^{\frac{10}{3}}_I W^{\frac{11}{20},\frac{10}{3}}_{\mathbb{T}^3}}^{\frac{10}{3}}\right)^{\frac{3}{10}} \lesssim |I|^{-\frac{3}{10}} \overset{(5.22)}{\asymp} \|u^\sharp\|_{L^\infty_{[0,T]}H^2_{\mathbb{T}^3}}^{\frac{3}{7}\left(\frac{17}{20}-\beta\right)+},$$

thanks to the flow property of the equation and energy conservation.

Furthermore, if we bound the right hand side of mild formulation

$$u^\sharp(t) = e^{-it\mathcal{H}^\sharp} u_0^\sharp - i\int_0^t e^{-i(t-s)\mathcal{H}^\sharp} \Lambda^{-1}[\Lambda u^\sharp(s)(|\Lambda u^\sharp(s)|^2 * V_\beta)]ds$$

in $H^2_{\mathbb{T}^3}$, we can bound for $0 < T \sim 1$ that

$$\|u^\sharp\|_{L^{\frac{10}{3}}_{[0,T]} W^{\frac{31}{20},\frac{10}{3}}_{\mathbb{T}^3}} \lesssim \|u_0^\sharp\|_{H^2_{\mathbb{T}^3}} + \left\|\int_0^t e^{-i(t-s)\mathcal{H}^\sharp} \Lambda^{-1}[\Lambda u^\sharp(s)(|\Lambda u^\sharp(s)|^2 * V_\beta)]ds\right\|_{L^{\frac{10}{3}}_{t;[0,T]} W^{\frac{31}{20},\frac{10}{3}}_{\mathbb{T}^3}}$$

$$\lesssim \|u_0^\sharp\|_{H^2_{\mathbb{T}^3}} + \left\|\int_0^t (-\mathcal{H}^\sharp)^{\frac{31}{40}-} e^{-i(t-s)\mathcal{H}^\sharp} \Lambda^{-1}[\Lambda u^\sharp(s)(|\Lambda u^\sharp(s)|^2 * V_\beta)]ds\right\|_{L^{\frac{10}{3}}_{[0,T]\times\mathbb{T}^3}},$$

having used Lemma 3.5. Now we can perform again the integration by parts with respect to time as in the proof of Proposition 4.1 to get

$$\int_0^t (-\mathcal{H}^\sharp)^{\frac{31}{40}-} e^{-i(t-s)\mathcal{H}^\sharp} \Lambda^{-1}[\Lambda u^\sharp(s)(|\Lambda u^\sharp(s)|^2 * V_\beta)]ds$$
$$= \int_0^t -\mathcal{H}^\sharp e^{-i(t-s)\mathcal{H}^\sharp} (-\mathcal{H}^\sharp)^{-\frac{9}{40}-} \Lambda^{-1}[\Lambda u^\sharp(s)(|\Lambda u^\sharp(s)|^2 * V_\beta)]ds$$
$$= i\int_0^t e^{-i(t-s)\mathcal{H}^\sharp} (-\mathcal{H}^\sharp)^{-\frac{9}{40}-} \Lambda^{-1} \partial_s[\Lambda u^\sharp(s)(|\Lambda u^\sharp(s)|^2 * V_\beta)]ds + \mathrm{Bd}[u^\sharp](t),$$

where the boundary terms are given by

$$\mathrm{Bd}[u^\sharp](t) = i(-\mathcal{H}^\sharp)^{-\frac{9}{40}-} \Lambda^{-1}[\Lambda u^\sharp(t)(|\Lambda u^\sharp(t)|^2 * V_\beta)]$$
$$+ ie^{-it\mathcal{H}^\sharp} (-\mathcal{H}^\sharp)^{-\frac{9}{40}-} \Lambda^{-1}[\Lambda u^\sharp(0)(|\Lambda u^\sharp(0)|^2 * V_\beta)],$$

and they can be bounded in $L^{\frac{10}{3}}_{[0,T]\times\mathbb{T}^3}$ by $\|\mathrm{Bd}[u^\sharp]\|_{L^{\frac{10}{3}}_{[0,T]\times\mathbb{T}^3}} \lesssim \|u^\sharp\|^3_{L^\infty_{[0,T]}H^1_{\mathbb{T}^3}} \leq C(\mathcal{E}_0(u_0))$. We postpone the proof of Lemma B.2 to the appendix and are therefore left to control

$$\left\|\int_0^t e^{-i(t-s)\mathcal{H}^\sharp} (-\mathcal{H}^\sharp)^{-\frac{9}{40}-} \Lambda^{-1} \partial_s[\Lambda u^\sharp(s)(|\Lambda u^\sharp(s)|^2 * V_\beta)]ds\right\|_{L^{\frac{10}{3}}_{[0,T]\times\mathbb{T}^3}}$$

$$\lesssim \int_0^T \left\|(-\mathcal{H}^\sharp)^{-\frac{9}{40}-} \Lambda^{-1} \partial_s[\Lambda u^\sharp(s)(|\Lambda u^\sharp(s)|^2 * V_\beta)]\right\|_{H^{\frac{9}{20}+}_{\mathbb{T}^3}} ds$$

$$\lesssim \int_0^T \|[\partial_s \Lambda u^\sharp(s)](|\Lambda u^\sharp(s)|^2 * V_\beta)\|_{L^2_{\mathbb{T}^3}} + \|\Lambda u^\sharp(s)(\partial_s|\Lambda u^\sharp(s)|^2 * V_\beta)\|_{L^2_{\mathbb{T}^3}} ds = \mathrm{(I)} + \mathrm{(II)}.$$

Having used the Strichartz estimate Lemma 3.4 and Lemma 3.5.

We can now exploit the bound (5.5) we have derived using Strichartz estimates in $H^1_{\mathbb{T}^3}$ so that the first term above can be controlled. Since the bound for (II) is identical, we only present the bound of the first term as follows

$$(\mathrm{I}) \lesssim \int_0^T \|\partial_s u^\natural(s)\|_{L^2_{\mathbb{T}^3}} \||\Lambda u^\natural(s)|^2 * V_\beta\|_{L^\infty_{\mathbb{T}^3}} ds$$

$$\lesssim \|V_\beta\|_{L^{\frac{3}{3-\beta}-}_{\mathbb{T}^3}} \int_0^T \|\partial_s u^\natural(s)\|_{L^2_{\mathbb{T}^3}} \||u^\natural(s)|^2\|_{L^{\frac{3}{\beta}+}_{\mathbb{T}^3}} ds,$$

where we notice that if $\frac{7}{10} < \beta < 1$ then $6 < \frac{6}{\beta} < \frac{60}{7}$ and $W^{\frac{11}{20}-,\frac{10}{3}}_{\mathbb{T}^3} \hookrightarrow L^{\frac{60}{7}-}_{\mathbb{T}^3}$ and $W^{\frac{9}{10}-\frac{\beta}{2}+,\frac{10}{3}}_{\mathbb{T}^3} \hookrightarrow L^{\frac{6}{\beta}+}_{\mathbb{T}^3}$ hence

$$\||u^\natural(s)|^2\|_{L^{\frac{3}{\beta}+}_{\mathbb{T}^3}} = \|u^\natural(s)\|^2_{L^{\frac{6}{\beta}+}_{\mathbb{T}^3}} \lesssim \|u^\sharp(s)\|^{\frac{20}{3}(1-\beta)+}_{L^{\frac{60}{7}-}_{\mathbb{T}^3}} \|u^\sharp(s)\|^{\frac{20}{3}(\beta-\frac{7}{10})+}_{L^6_{\mathbb{T}^3}} \lesssim \|u^\sharp(s)\|^{\frac{20}{3}(1-\beta)+}_{W^{\frac{11}{20}-,\frac{10}{3}}_{\mathbb{T}^3}} \|u^\sharp(s)\|^{\frac{20}{3}(\beta-\frac{7}{10})-}_{H^1_{\mathbb{T}^3}},$$

by interpolation and noting that the last term is controlled again simply by the energy.

This gives

$$\||u^\sharp|^2\|_{L^1_{[0,T]}L^{\frac{3}{\beta}}_{\mathbb{T}^3}} \lesssim T^{1-2(1-\beta)+} \|u^\sharp\|^{\frac{20}{3}(1-\beta)+}_{L^{\frac{10}{3}}_{[0,T]}W^{\frac{11}{20},\frac{10}{3}}_{\mathbb{T}^3}}. \tag{5.23}$$

We can therefore conclude with (5.5) that

$$\|u^\sharp\|_{L^{\frac{10}{3}}_{[0,T]}W^{\frac{31}{20}-,\frac{10}{3}}_{\mathbb{T}^3}} \lesssim \|u^\sharp\|^{1+\frac{20}{7}\left(\frac{17}{20}-\beta\right)(1-\beta)+}_{L^\infty_{[0,T]}H^2_{\mathbb{T}^3}}, \quad \beta \in \left(\frac{7}{10}, \frac{17}{20}\right). \tag{5.24} \quad \square$$

## 5.3 The case of bounded interaction

In this chapter we discuss how to obtain global well-posedness of

$$i\partial_t u = \mathcal{H}u - u(|u|^2 * V), \quad \text{on } [0,T] \times \mathbb{T}^3 \tag{5.25}$$

in the case of bounded potential, i.e. $V = L^\infty_{\mathbb{T}^3}$ which is in many ways simpler, since the nonlinear term $V * |u|^2$ is controlled in $L^\infty_{[0,T]} L^\infty_{\mathbb{T}^3}$ by the conservation of mass through

$$\|V * |u|^2(t)\|_{L^\infty_{\mathbb{T}^3}} \lesssim \|V\|_{L^\infty_{\mathbb{T}^3}} \||u|^2(t)\|_{L^1_{\mathbb{T}^3}} = \|V\|_{L^\infty_{\mathbb{T}^3}} \|u(t)\|^2_{L^2_{\mathbb{T}^3}} = \|V\|_{L^\infty_{\mathbb{T}^3}} \|u_0\|^2_{L^2_{\mathbb{T}^3}}. \tag{5.26}$$

This shows that one could obtain global well-posedness in the space $L^2_{\mathbb{T}^3}$ using the mass conservation (5.2), since in that space, the nonlinearity is Lipschitz.

**Proposition 5.13.** *Let $u_0 \in L^2_{\mathbb{T}^3}$, $V \in L^\infty_{\mathbb{T}^3}$, and $T > 0$. Then there exists a unique mild solution to (5.25) in the space $C([0,T]; L^2_{\mathbb{T}^3})$ which depends continuously on the data.*

**Proof.** This is quite standard using the mild formulation

$$u(t) = e^{-it\mathcal{H}} u_0 - i \int_0^t e^{-i(t-s)\mathcal{H}} [u(s) V * |u(s)|^2] ds \tag{5.27}$$

and the bound (5.26) to get first LWP by a contraction argument and GWP by the conservation of mass, see (5.2) which can be justified again by an approximation argument, we omit the details. □

Similarly, in the strong setting i.e. $u_0 \in \mathcal{D}(\mathcal{H})$, one gets local well-posedness in the same way as in Section 4.1 i.e. by integrating by parts in time in the mild fomulation. In order to get global well-posedness, one simply applies Gronwall inequality to the mild formulation of the time derivative. If we assume that $V$ has some regularity e.g. $W_{\mathbb{T}^3}^{\beta,1}$ then this was already proven, but the global well-posedness still holds if only $V \in L_{\mathbb{T}^3}^\infty$ with a Gronwall-type argument.

**Theorem 5.14.** *Let $u_0 \in \mathcal{D}(\mathcal{H})$ and $T > 0$ and $V \in L_{\mathbb{T}^3}^\infty$ and non-negative. Then there exists a unique mild solution to (5.1)*

$$u \in C([0,T]; \mathcal{D}(\mathcal{H})) \cap C^1([0,T]; L_{\mathbb{T}^3}^2)$$

*with the solution depending continuously on the initial data.*

*Moreover, there are some constants $c(\mathcal{E}_0(u_0)), c'(\mathcal{E}_0(u_0)) > 0$ depending on the parameters and the energy of the initial data $\mathcal{E}_0(u_0)$ s.t. we have the bound*

$$\|u\|_{C([0,T];\mathcal{D}(\mathcal{H})) \cap C^1([0,T];L_{\mathbb{T}^3}^2)} \lesssim c(\mathcal{E}_0(u_0)) \, e^{c'(\mathcal{E}_0(u_0))T} \|u_0\|_{\mathcal{D}(\mathcal{H})}. \tag{5.28}$$

**Proof.** This is proved basically the same way as Theorem 5.3, integrating by parts in time, using (5.26) and applying a Gronwall-type argument to the time derivative $\partial_t u$. □

In this setting there is an interesting subtlety, namely that it is not obvious how to obtain energy solutions[5.3] from the mild formulation, basically because the term $V * |u|^2$ can not have more than $\frac{1}{2}-$ spatial derivatives so even in the transformed setting, one can not close the fixed point. It is even less evident how in (5.27) one would directly bound

$$\left\| \int_0^t e^{-i(t-s)\mathcal{H}} u(s) V * |u(s)|^2 \right\|_{\mathcal{D}(\sqrt{-\mathcal{H}})}.$$

For a similar issue to what happens in the two-dimensional NLS case, see [GUZ20] Theorem 3.11 and [Zac19] Theorem 1.3.

**Proposition 5.15.** *Let $u_0 \in \mathcal{D}(\sqrt{-\mathcal{H}})$ and $T > 0$ and $V \in L_{\mathbb{T}^3}^\infty$ and non-negative. Then there exists a unique mild solution to (5.1)*

$$u \in C([0,T], \mathcal{D}(\sqrt{-\mathcal{H}})),$$

*with the solution depending continuously on the initial data.*

**Proof.** One obtains a solution by approximating the initial data as in the proof of Corollary 5.1, which exists globally by Theorem 5.14 and using the conserved energy as an apriori estimate as in [GUZ20]. The solution is also unique by Theorem 5.13.

---

5.3. c.f. [GUZ20] 3.2.2.

To get the continuity of the solution map in $\mathcal{D}(\sqrt{-\mathcal{H}})$ as well as continuity in time, it suffices to show that the solution map of the mild formulation of the sharpened equation (4.12), where we replace $V_\beta$ by $V \in L^\infty$ as in the statement, verifies the conditions of nonlinear interpolation (see [ABI+24] Remark 4 after Theorem 1 or Theorem 18 for the full generality).

Namely, we have to check that given $R, T > 0$, $u_0^\# \in H^1_{\mathbb{T}^3}$ and the solution map

$$\mathcal{S}: B_{H^1_{\mathbb{T}^3}}(u_0^\#, R) \to L^\infty([0,T], L^2_{\mathbb{T}^3}), \quad u_0^\# \mapsto u^\# \text{ as the RHS of (4.12) with } V_\beta = V, \tag{5.29}$$

the following conditions hold:

i. (Weak Lipschitz) there exists $C_0 = C_0(R, T, u_0^\#) < \infty$, such that

$$v_0^\#, w_0^\# \in B_{H^1_{\mathbb{T}^3}}(u_0, R) \Longrightarrow \|\mathcal{S}(v_0^\#) - \mathcal{S}(w_0^\#)\|_{L^\infty([0,T], L^2_{\mathbb{T}^3})} \leq C_0 \|v_0^\# - w_0^\#\|_{L^2_{\mathbb{T}^3}}, \tag{5.30}$$

ii. (Tame estimate) there exists $C_1 = C_1(R, T, u_0^\#) < \infty$, such that

$$v_0^\# \in B_{H^1_{\mathbb{T}^3}}(u_0^\#, R) \cap C^\infty_{\mathbb{T}^3} \Longrightarrow \|\mathcal{S}(v_0^\#)\|_{L^\infty([0,T], H^2_{\mathbb{T}^3})} \leq C_0 \|v_0^\#\|_{H^2_{\mathbb{T}^3}}, \tag{5.31}$$

iii. (Time Continuity)

$$v_0^\# \in B_{H^1_{\mathbb{T}^3}}(u_0^\#, R) \Longrightarrow \mathcal{S}(v_0^\#) \in C([0,T], L^2_{\mathbb{T}^3}). \tag{5.32}$$

Since (5.32) follows directly from Proposition 5.13, and (5.31) from Theorem 5.14, more precisely (5.28), and the equivalence of norms (2.19), we are left to verify the weak Lipschitz bound (5.30). Abbreviating $v^\# := \mathcal{S}(v_0^\#)$ and $w^\#$ likewise, we have from the mild formulation

$$\begin{aligned}
\|v^\#(t) - w^\#(t)\|_{L^2_{\mathbb{T}^3}} &\lesssim \|e^{-it\mathcal{H}^\#}(v_0^\# - w_0^\#)\|_{L^2_{\mathbb{T}^3}} \\
&\quad + \int_0^t \|\Gamma v^\#(s)(|e^W \Gamma v^\#(s)|^2 * V) - \Gamma w^\#(s)(|e^W \Gamma w^\#(s)|^2 * V)\|_{L^2_{\mathbb{T}^3}} ds \\
&\lesssim \|v_0^\# - w_0^\#\|_{L^2_{\mathbb{T}^3}} + \int_0^t \|\Gamma[v^\#(s) - w^\#(s)](|e^W \Gamma v^\#(s)|^2 * V)\|_{L^2_{\mathbb{T}^3}} ds \\
&\quad + \int_0^t \|\Gamma w^\#(s)([|e^W \Gamma v^\#(s)|^2 - |e^W \Gamma w^\#(s)|^2] * V)\|_{L^2_{\mathbb{T}^3}} ds \\
&\lesssim \|v_0^\# - w_0^\#\|_{L^2_{\mathbb{T}^3}} + \|V\|_{L^\infty} \int_0^t \|v^\#(s) - w^\#(s)\|_{L^2_{\mathbb{T}^3}} \|e^W \Gamma v^\#(s)\|^2_{L^2_{\mathbb{T}^3}} ds \\
&\quad + \|V\|_{L^\infty} \int_0^t \|\Gamma w^\#(s)\|_{L^2_{\mathbb{T}^3}} \||e^W \Gamma v^\#(s)|^2 - |e^W \Gamma w^\#(s)|^2\|_{L^1_{\mathbb{T}^3}} ds, \\
&\lesssim_R \|v_0^\# - w_0^\#\|_{L^2_{\mathbb{T}^3}} + \int_0^t \|v^\#(s) - w^\#(s)\|_{L^2_{\mathbb{T}^3}} ds,
\end{aligned}$$

where we used the conservation of mass

$$\|e^W \Gamma w^\#(s)\|^2_{L^2_{\mathbb{T}^3}} + \|e^W \Gamma v^\#(s)\|^2_{L^2_{\mathbb{T}^3}} = \|w_0\|^2_{L^2} + \|v_0\|^2_{L^2} \lesssim R^2$$

as in Proposition 5.13. Moreover, using again the conservation of mass

$$\begin{aligned}
\||e^W \Gamma v^\#(s)|^2 - |e^W \Gamma w^\#(s)|^2\|_{L^1_{\mathbb{T}^3}} &\lesssim (\|e^W \Gamma v^\#(s)\|_{L^2_{\mathbb{T}^3}} + \|e^W \Gamma w^\#(s)\|_{L^2_{\mathbb{T}^3}}) \|e^W \Gamma[v^\#(s) - w^\#(s)]\|_{L^2_{\mathbb{T}^3}} \\
&\lesssim_R \|v^\#(s) - w^\#(s)\|_{L^2_{\mathbb{T}^3}}.
\end{aligned}$$

(5.30) then follows from Grönwall's inequality. □

# 6 Mean-Field Limit of N-particle Systems

In this chapter we discuss the convergence of the solutions to the $N$-particle linear Schrödinger equations (1.2) whose linear part is given by the mollified Anderson Hamiltonian to the solution to the Hartree NLS (1.1) in a suitable sense. More precisely, we recall from (1.2) the Hamiltonian of an $N$-particle Bosonic system with mean field interaction in a random external potential as

$$H_{N,\delta} := -\frac{1}{2}\sum_{k=1}^{N} \mathcal{H}_{\delta,x_k} + \frac{1}{N}\sum_{1 \leq j < k \leq N} V_{j,k},$$

on $L_s^2((\mathbb{T}^3)^N)$ with interaction potential $V$ and the self-adjoint operator $\mathcal{H}_{\delta,x_k}$ on $L^2(\mathbb{T}^3)$ is defined as the operator

$$\mathcal{H}_\delta = \Delta + \xi_\delta(\cdot) + c_\delta, \tag{6.1}$$

acting on the $k$-th variable, where $\xi_\delta$ is Gaussian with smooth spatial correlations for any $\delta > 0$, which tends to spatial white noise as $\delta \to 0^+$, while $(c_\delta)_{\delta > 0}$ is the family of renormalization constants discussed in Section 2.1. We also recall that $V_{j,k}$ is used to denote the multiplication operator with $V(x_j - x_k)$. The Heisenberg picture of the dynamic is given by

$$i\partial_t \rho_{N,\delta}(t) = [H_{N,\delta}, \rho_{N,\delta}(t)], \tag{6.2}$$

where $\rho_{N,\delta}$ is the density matrix of the system and $[\cdot,\cdot]$ is the commutator. We assume that the mollification vanishes as the number of particles increases, that is we have some sequences $\delta \to 0^+$ and $N \to \infty$ simultaneously, so that in the limit we have a genuine white-noise potential but all the approximate objects are smooth and classically defined.

Notice that $\rho_{N,\delta}$ is not necessarily the projection onto wave functions, but in case it is i.e.

$$\rho_{N,\delta} = \Pi_{\Psi_{N,\delta}}$$

for an $N$-particle wave function $\Psi_{N,\delta}$, or equivalently in the form of kernels, we can define the path of density matrices $\rho_{N,\delta} \in L^\infty(\mathbb{R}_+, \mathcal{L}_N^1)$ as

$$\rho_{N,\delta}(t, x_{1:N}, y_{1:N}) = \Psi_{N,\delta}(t, x_{1:N}) \overline{\Psi_{N,\delta}(t, y_{1:N})},$$

and the Heisenberg picture (6.2) is equivalent to the linear Schrödinger equation

$$i\partial_t \Psi_{N,\delta} = H_{N,\delta}\Psi_{N,\delta}, \tag{6.3}$$

which is

$$i\partial_t \Psi_{N,\delta} = -\frac{1}{2}\sum_{k=1}^{N}(\Delta_{x_k} + \xi_\delta(x_k) + c_\delta)\Psi_{N,\delta} + \frac{1}{N}\sum_{1 \leq j < k \leq N} V(x_j - x_k)\Psi_{N,\delta}, \tag{6.4}$$

where $\Psi_{N,\delta} \in L^\infty(\mathbb{R}_+, L_s^2((\mathbb{T}^3)^N))$ is the unique solution to the equation in $L_s^2((\mathbb{T}^3)^N)$ with initial condition $\|\Psi_{N,\delta}(0)\|_{L^2} = 1$ for $\delta > 0$.

We recall that $\xi_\delta = a_\delta * \xi$ with $a_\delta$ being a smooth mollifier and the external potential is in the Kato class (c.f. [RS80] X.2). In that case, it follows from (6.3) that the energy and mass conservations hold, that is

$$\frac{d}{dt}\langle H_{N,\delta}\Psi_{N,\delta}(t), \Psi_{N,\delta}(t)\rangle = 0, \quad \frac{d}{dt}\|\Psi_{N,\delta}(t)\|^2_{L^2(\mathbb{T}^{3N})} = 0, \tag{6.5}$$

whenever $|\langle H_{N,\delta}\Psi_{N,\delta}(0), \Psi_{N,\delta}(0)\rangle| < +\infty$. Now, we can define its $n$-point density matrix, which should be interpreted as the $n$-particle density marginal of the system for the "first" $n$-particles $1 \leq n \leq N$, as $\rho_{N,\delta}^n := \text{Tr}_{n+1:N}(\rho_{N,\delta})$[6.1], or in kernel forms: $\rho_{N,\delta}^n = \rho_{N,\delta}$ when $n = N$, and

$$\rho_{N,\delta}^n(t, x_{1:n}, y_{1:n}) = \int_{\mathbb{T}^{3(N-n)}} \Psi_{N,\delta}(x_{1:n}, z_{n+1:N})\overline{\Psi_{N,\delta}(t, y_{1:n}, z_{n+1:N})}dz_{n+1:N}, \quad n < N.$$

It is possible to derive equations for the trace class operators $(\rho_{N,\delta}^n)_{n=1,\ldots,N}$. Indeed, we have

$$\begin{aligned} i\partial_t \rho_{N,\delta}^n &= \text{Tr}_{n+1:N}([H_{N,\delta}^n, \rho_{N,\delta}]) + \text{Tr}_{n+1:N}([H_{N,\delta} - H_{N,\delta}^n, \rho_{N,\delta}]) \\ &= [H_{N,\delta}^n, \text{Tr}_{n+1:N}(\rho_{N,\delta})] + \text{Tr}_{n+1:N}([H_{N,\delta} - H_{N,\delta}^n, \rho_{N,\delta}]), \end{aligned}$$

and we abbreviate $H_{N,\delta}^n := -\frac{1}{2}\sum_{k=1}^n \mathcal{H}_{\delta,x_k} + \frac{1}{N}\sum_{1 \leq j < k \leq n} V_{j,k}$, and the second term on the R.H.S. can be further expanded by

$$\begin{aligned} \text{Tr}_{n+1:N}([H_{N,\delta} - H_{N,\delta}^n, \rho_{N,\delta}]) &= -\frac{1}{2}\sum_{j=n+1}^N \text{Tr}_{n+1:N}([\mathcal{H}_{\delta,x_j}, \rho_{N,\delta}]) \\ &\quad + \frac{1}{N}\sum_{j<n<k \text{ or } n<j<k} \text{Tr}_{n+1:N}([V_{j,k}, \rho_{N,\delta}]). \end{aligned} \tag{6.6}$$

We further notice that for $n + 1 \leq l < j \leq N$, we have

$$\text{Tr}_{n+1:N}([V_{j,k}, \rho_{N,\delta}]) = \text{Tr}_{j,k}([V_{j,k}, \text{Tr}_{n+1:\hat{j}:\hat{k}:N}\rho_{N,\delta}]) = 0$$

due to the cyclic permutation invariance of the trace operation, where $\text{Tr}_{n+1:\hat{j}:\hat{k}:N}$ is the partial trace with respect to particles $\{n+1,\ldots,N\}\setminus\{j,k\}$. Plugging it back into (6.6) and using the symmetry of the underlying spaces, we have

$$i\partial_t \rho_{N,\delta}^n = [H_{N,\delta}^n, \rho_{N,\delta}^n] + \frac{N-n}{N}\sum_{k \leq n} \text{Tr}_{n+1}([V_{k,n+1}, \rho_{N,\delta}^{n+1}]), \tag{6.7}$$

or equivalently in kernel form,

$$\begin{aligned} i\partial_t \rho_{N,\delta}^n &= -\frac{1}{2}\left(\sum_{k=1}^n \mathcal{H}_{\delta,x_k} - \sum_{k=1}^n \mathcal{H}_{\delta,y_k}\right)\rho_{N,\delta}^n + \frac{1}{N}\sum_{j<k \leq n}[V(x_j - x_k) - V(y_j - y_k)]\rho_{N,\delta}^n \\ &\quad + \frac{N-n}{N}\sum_{k=1}^n \int_{\mathbb{T}^3}[V(x_k - z) - V(y_k - z)]\rho_{N,\delta}^{n+1}((x_{1:n}, z), (y_{1:n}, z))dz, \end{aligned} \tag{6.8}$$

---

6.1. Note that $\rho_{N,\delta}^n$ defined above are always density matrices, that is they are trace-class non-negative self-adjoint operators such that $\text{Tr}(\rho_{N,\delta}^n) = 1$.

By sending $N\to\infty$ and $\delta\to 0^+$ simultaneously– since we expect $\mathcal{H}_\delta\to\mathcal{H}$ in a suitable sense from Section 2.1– we have then derived formally the limiting equation, i.e. the BBGKY hierarchy which is

$$i\partial_t\rho_\infty^n = -\frac{1}{2}\sum_{k=1}^n [\mathcal{H}_{x_k},\rho_\infty^n] + \sum_{1\leq k\leq n}\mathrm{Tr}_{n+1}([V_{k,n+1},\rho_\infty^{n+1}]) \tag{6.9}$$

where $(\rho_\infty^n)_{n\in\mathbb{N}}$ is a suitbale sequence of trace class operators. The corresponding Schrödinger equation with Hartree-type nonlinearity, which is just (1.1), reads

$$i\partial_t\Psi = -\frac{1}{2}\mathcal{H}\Psi + \Psi(V*|\Psi|^2) \tag{6.10}$$

and to witness the relation between (6.9) and (6.10), we formulate (6.10) equivalently in the Heisenberg picture, that is if we let $\rho:=\Pi_\Psi$, and therefore $|\Psi|^2(z)=\rho(z,z)$ for $z\in\mathbb{T}^3$ in the language of kernels (1.2), we have for any test function $\Phi$ in the domain of $\mathcal{H}$ that

$$\begin{aligned}i\partial_t\rho_t(\Phi) &= -\frac{1}{2}[\mathcal{H},\Pi_\Psi]\Phi - \langle(V*|\Psi|^2)\Phi,\Psi\rangle\Psi + (V*|\Psi|^2)\langle\Phi,\Psi\rangle\Psi \\ &= \left[-\frac{1}{2}\mathcal{H} + \int V(\cdot-z)\rho_t(z,z)\mathrm{d}z,\rho_t\right]\Phi.\end{aligned} \tag{6.11}$$

The consideration of the ANLS (6.10) is *formally* justified, since the family of tensorizations $(\rho^{\otimes n})_n$ of the solutions to the Hartree equation (6.10) satisfies the BBGKY hierarchy (6.9). Indeed, according to the derivation above, it suffices to notice that

$$\left[\int V(\cdot-z)\rho_t(z,z)\mathrm{d}z,\rho_t\right](x_1,y_1) = \mathrm{Tr}_{x_2}[V(\cdot-x_2),\rho_t^{\otimes 2}](x_1,y_1), \tag{6.12}$$

for any solution $\Psi$ to the nonlinear Schrödinger equation (6.10) and $\rho_t^{\otimes 2}=\Pi_{\Psi(t,\cdot)^{\otimes 2}}$ for $t\geq 0$ acts on two particles. The assertion for $\rho^{\otimes n}$ for $n\geq 2$ follows by induction.

We have essentially used formal arguments in two steps above, namely in the convergence of the many-body density evolution to the BBGKY hierarchy and in the identification of the two solutions of the BBGKY hierarchy. In the rest of this chapter, we will make those two steps rigorous in the following Sections 6.2 and 6.3 respectively under different conditions of the interaction potential $V$.

## 6.1 Precise Formulation

In order to give a precise definition of solution to the BBGKY hierarchy (6.9), we denote the unitary flow on the bounded operators on $L_s^2((\mathbb{T}^3)^n)$ as

$$\mathcal{U}_t^n(\rho) := \exp\left(\frac{it}{2}\sum_{k=1}^n \mathcal{H}_{x_k}\right)\rho\exp\left(-\frac{it}{2}\sum_{k=1}^n \mathcal{H}_{x_k}\right),$$

or equivalently in kernel form by

$$\mathcal{U}_t^n(\rho)(x_{1:n},y_{1:n}) := \exp\left(\frac{it}{2}\sum_{k=1}^n \mathcal{H}_{x_k} - \frac{it}{2}\sum_{k=1}^n \mathcal{H}_{y_k}\right)\rho(x_{1:n},y_{1:n}), \tag{6.13}$$

for any $t\in\mathbb{R}$, and $t\mapsto\mathcal{U}_{\delta,t}^n(\rho)$ is defined similarly by replacing $\mathcal{H}$ with $\mathcal{H}_\delta$. Since $V_{j,k}$ is the multiplication operator with $V(x_j-x_k)$, we have also in kernel form that

$$[V_{j,k},\rho^n](x_{1:n},y_{1:n}) = [V(x_j-x_k) - V(y_j-y_k)]\rho^n(x_{1:n},y_{1:n}).$$

We recall from (1.5) that $\rho^n \in \mathcal{L}_n^p$ belongs to the space $\mathcal{W}_{n,\delta}^{p,\alpha}$ if

$$\|\rho^n\|_{\mathcal{W}_{n,\delta}^{p,\alpha}} = n\|(-\mathcal{H}_{\delta,x_1})^{\alpha/2}\rho^n(-\mathcal{H}_{\delta,x_1})^{\alpha/2}\|_{\mathcal{L}_n^p} < +\infty, \tag{6.14}$$

using the symmetry of the underlying space. We denote by $\mathcal{W}_n^{p,\alpha}$ the subspace of $\mathcal{L}_n^p$ such that the norm $\|\cdot\|_{\mathcal{W}_n^{p,\alpha}}$, which is defined as in (6.14) where the operator $\mathcal{H}_\delta$ is replaced by the limiting Anderson Hamiltonian $\mathcal{H}$, is finite. The notion of solutions to (6.9) can be formulated as follows:

**Definition 6.1.** *Consider a sequence of time-dependent trace-class non-negative self-adjoint operators $(\rho_\infty^n)_{n\in\mathbb{N}}$, such that $\rho_\infty^n \in L^\infty(\mathbb{R}_+, \mathcal{L}_n^1)$. We say that $(\rho_\infty^n)_{n\in\mathbb{N}}$ is a mild solution to the hierarchy equation (6.7) with initial condition $(\rho_{\infty,0}^n)_{n\in\mathbb{N}}$ if, for any $n\in\mathbb{N}$,*

$$\rho_\infty^n(t) = \mathcal{U}_t^n(\rho_{\infty,0}^n) + \sum_{k=1}^n \int_0^t \mathcal{U}_{t-s}^n(\mathrm{Tr}_{n+1}[V_{k,n+1}, \rho_\infty^{n+1}(s)])\mathrm{d}s, \tag{6.15}$$

*and furthermore, $\rho_\infty^n \in L_{\mathbb{R}_+}^\infty \mathcal{W}_n^{1,1}$, i.e.*

$$\sup_{t\geq 0}\|\rho_\infty^n(t)\|_{\mathcal{W}_n^{1,1}} = n\sup_{t\geq 0}\|(-\mathcal{H}_{x_1})^{1/2}\rho_\infty^n(t)(-\mathcal{H}_{x_1})^{1/2}\|_{\mathcal{L}^1(L_s^2((\mathbb{T}^3)^n))} < +\infty. \tag{6.16}$$

**Notation 6.2.** *We recall that a sequence converges weakly\* in $\mathcal{L}_n^1 = \mathcal{L}^1(L_s^2((\mathbb{T}^3)^n))$ in the sense that it converges weakly with respect to the duality $(\mathcal{K}_n)^* \approx \mathcal{L}_n^1$, where $\mathcal{K}_n = \mathcal{K}(L_s^2((\mathbb{T}^3)^n))$ is the space of compact operators defined on the Hilbert space $L_s^2((\mathbb{T}^3)^n)$.*

*Since $\mathcal{K}_n$ is separable, the bounded sets of $\mathcal{L}_n^1$ are metrizable if equipped with the weak\* topology. Thus when saying that a sequence, depending on a parameter (usually the time or mollifiction parameters), converges uniformly in $\mathcal{L}_n^1$ w.r.t. the weak\* topology we mean that the sequence is bounded in the trace norm and converges uniformly with respect that parameter, when $\mathcal{L}_n^1$ is equipped with such a metric inducing the weak\* topology.*

*We also say that a sequence converges weakly\* in $L_{\mathrm{loc}}^\infty(\mathbb{R}_+, \mathcal{L}^1)$ in the sense of weak\* convergence if it converges weakly\* with respect to the duality*

$$(L^1([0,T],\mathcal{K}_n), L^\infty([0,T],\mathcal{L}_n^1)) \text{ for all } T>0. \tag{6.17}$$

We are ready to state the main convergence theorem of this section.

**Theorem 6.3.** *Suppose that the potential $V:\mathbb{T}^3 \to \mathbb{R}$ satisfies*

$$V = \frac{V_1}{|\cdot|^{3-\beta}} + V_2, \quad V_1 \in L_{\mathbb{T}^3}^\infty, V_2 \in L_{\mathbb{T}^3}^{3+}, \tag{6.18}$$

*where $3 > \beta \geq 2$. Consider $\Psi_{N,\delta}(0) \in L_s^2((\mathbb{T}^3)^N)$ such that*

$$\|\Psi_{N,\delta}(0)\|_{L_s^2((\mathbb{T}^3)^N)} = 1, \quad \sup_{N\in\mathbb{N}, 0<\delta\leq 1}\|(-\mathcal{H}_{\delta,x_1})^{1/2}\Psi_{N,\delta}(0)\|_{L_s^2((\mathbb{T}^3)^N)} < +\infty. \tag{6.19}$$

*Then there exists subsequences $N_k \to +\infty$ and $\delta_k \to 0^+$ together with a family of density matrices $(\rho_\infty^n)_{n\in\mathbb{N}}$ such that, for each $n\geq 0$, $\rho_\infty^n \in L^\infty(\mathbb{R}_+, \mathcal{W}_n^{1,1})$ and $\rho_{N_k,\delta_k}^n \to \rho_\infty^n$ weakly\* in $L_{\mathrm{loc}}^\infty(\mathbb{R}_+, \mathcal{L}_n^1)$.*

*Furthermore, the limit $(\rho_\infty^n)_{n\in\mathbb{N}}$ is a mild $\mathcal{W}^{1,1}$ solution to the hierarchy equation (6.1) in the sense of Definition 6.1.*

To justify the formal discussion under (6.10), we suppose that the potential $V$ as in (6.18) belongs to $W_{\mathbb{T}^3}^{17/20+,1}$ or $L_{\mathbb{T}^3}^\infty$ so that there is a unique continuous global (mild) solution

$$u_0 \in \mathcal{D}(\sqrt{-\mathcal{H}}) \Rightarrow t \mapsto u(t) \in C([0,T]; \mathcal{D}(\sqrt{-\mathcal{H}})) \, \forall T > 0$$

to equation (5.1) according to Section 4 and 5, which implies that the sequence of compatible operators $(\rho_\infty^n)_{n \in \mathbb{N}}$ defined by

$$\rho_\infty^n(t) = \Pi_{(u(t,\cdot))^{\otimes n}},$$

is a mild solution to the hierarchy equation (6.15).

Consider now the setup in Theorem 6.3 where we have taken sequences $\delta_k \to 0^+$ and $N_k \to \infty$ as $k \to \infty$, together with a sequence of initial conditions $\Psi_{N_k, \delta_k}(0)$ for the Schrödinger equation (6.4) satisfying the hypothesis of Theorem 6.3 and furthermore assume that

$$\rho_{N_k,\delta_k}^n(0, x_{1:n}, y_{1:n}) \to \Pi_{u_0}(x_{1:n}, y_{1:n}) = u_0^{\otimes n}(x_{1:n})\overline{u_0^{\otimes n}(y_{1:n})}, \quad k \to \infty,$$

weakly* in $\mathcal{W}_n^{1,1}$. As a consequence of Theorem 6.3 and the discussion above, we would expect that $\rho_{N_k,\delta_k}^n(t, x_{1:n}, y_{1:n})$ converges to $u^{\otimes n}(t, x_{1:n})\overline{u^{\otimes n}(t, y_{1:n})}$ weakly* in $L^\infty(\mathbb{R}_+, \mathcal{L}_n^1)$.

In general, this could be however not true, since it is not clear how "singular" the potential $V$ or the initial condition $(\rho_{\infty,0}^n)_{n \in \mathbb{N}}$ are allowed to be, such that the hierarchy equation (6.15) admits a unique mild solution. In the following theorem we give a *sufficient* condition under which the hierarchy equation (6.15) admits a unique mild solution, which then implies the convergence of the solutions to the $N$-particle linear Schrödinger equation.

**Theorem 6.4.** *Under the hypothesis of Theorem 6.3, we further assume that $V \in L_{\mathbb{T}^3}^\infty$ and the family of initial conditions $(\rho_{\infty,0}^n)_{n \in \mathbb{N}}$ satisfies $\rho_{\infty,0}^n \in \mathcal{L}_n^1$. The BBGKY hierarchy equation (6.15) admits a unique mild solution in the sense of Definition 6.1.*

As a consequence, solutions $(\Psi_{N,\delta})$ of the linear Schrödinger equation (6.4) with bounded interactions to the solution converge to the NLS (1.1).

**Corollary 6.5.** *Under the hypotheses of Theorem 6.4, we consider $\delta \to 0^+$ and $N \to \infty$, such that $\Psi_{N,\delta,0}$ satisfies (6.19), and $u_0 \in \mathcal{D}(\sqrt{-\mathcal{H}})$ verifies*

$$\mathrm{Tr}_{n+1:N}(\Pi_{\Psi_{N,\delta,0}}) = \rho_{N,\delta}^n(0, \cdot) \to \Pi_{u_0^{\otimes n}}, \quad \text{weakly* in } \mathcal{L}^1(L_s^2((\mathbb{T}^3)^n)), \quad N \to \infty, \delta \to 0^+.$$

*Then we have for any $n \in \mathbb{N}$*

$$\rho_{N,\delta}^n \to \Pi_{u^{\otimes n}}, \quad \text{weakly* in } L_{\mathrm{loc}}^\infty(\mathbb{R}_+, \mathcal{L}_n^1), \quad N \to \infty, \delta \to 0^+,$$

*where $u$ is the unique solution to equation (5.1) with initial condition $u(0, \cdot) = u_0$.*

## 6.2 Proof of the convergence, Theorem 6.3

We start with the following version of Hardy's inequality for the Anderson Hamiltonian, sometimes we write for brevity $\mathcal{H}_0 = \mathcal{H}$.

**Lemma 6.6. (Hardy's inequality)** *Let $f \in \mathcal{D}((-\mathcal{H}_\delta)^{1/2})$ for $\delta \in [0,1]$, then we have*

$$\int_{\mathbb{T}^3} \frac{|f(y)|^2}{|x-y|^2} dy \lesssim \|f\|^2_{\mathcal{D}((-\mathcal{H}_\delta)^{1/2})}, \tag{6.20}$$

*where the implicit multiplicative constant is uniform in $x \in \mathbb{T}^3$ and $0 \leq \delta \leq 1$. We have further the bound*

$$\int_{\mathbb{T}^3} \frac{|f(y)|^2}{|x-y|^{2s}} dy \leq \varepsilon \|f\|^2_{\mathcal{D}((-\mathcal{H}_\delta)^{1/2})} + C_\varepsilon \|f\|^2_{L^2_{\mathbb{T}^3}}, \tag{6.21}$$

*for any $\varepsilon > 0$, $0 \leq s < 1$ and $C_\varepsilon$ is uniform in $x \in \mathbb{T}^3$ and $0 \leq \delta \leq 1$.*

**Proof.** Following Theorem 2.57 in [BCD11][6.2], we have Hardy's inequality

$$\int_{\mathbb{T}^3} |x|^{-2s} |f(x)|^2 dx \lesssim \|f\|^2_{H^s_{\mathbb{T}^3}} \qquad 0 < s \leq 1$$

and we recall the norm equivalence (2.18) that $u \in \mathcal{D}(\sqrt{-\mathcal{H}})$ iff $e^{-W}u \in H^1_{\mathbb{T}^3}$, and $\|e^{-W}u\|_{H^1_{\mathbb{T}^3}} \asymp (-u, \mathcal{H}u)$. Thus we have

$$\int_{\mathbb{T}^3} |x|^{-2} |f(x)|^2 dx \leq \|e^{2W}\|_{L^\infty_{\mathbb{T}^3}} \int_{\mathbb{T}^3} |x|^{-2} |e^{-W}f(x)|^2 dx \lesssim \|e^{2W}\|_{L^\infty_{\mathbb{T}^3}} \|e^{-W}f\|^2_{H^1_{\mathbb{T}^3}} \lesssim \|f\|^2_{\mathcal{D}(\sqrt{-\mathcal{H}})},$$

with a multiplicative constant depending on the noise. For $s < 1$ we have

$$\int_{\mathbb{T}^3} |x|^{-2s} |f(x)|^2 dx \lesssim \|e^{-W}f\|^2_{H^s_{\mathbb{T}^3}} \leq \varepsilon \|f\|^2_{\mathcal{D}(\sqrt{-\mathcal{H}})} + C_\varepsilon \|f\|^2_{L^2_{\mathbb{T}^3}},$$

for arbitrarily small $\varepsilon > 0$ where we have used in the last inequality simply the interpolation

$$\|e^{-W}f\|_{H^s_{\mathbb{T}^3}} \lesssim \|e^{-W}f\|^s_{H^1_{\mathbb{T}^3}} \|e^{-W}f\|^{1-s}_{L^2_{\mathbb{T}^3}} \lesssim \tilde{\varepsilon} \|e^{-W}f\|_{H^1_{\mathbb{T}^3}} + C_{\tilde{\varepsilon}} \|e^{-W}f\|_{L^2_{\mathbb{T}^3}},$$

for any $\tilde{\varepsilon} > 0$. Clearly the same computation is true for $\mathcal{H}$ replaced by $\mathcal{H}_\delta$. □

Thanks to the previous result we obtain the following lemma.

**Lemma 6.7.** *For the interaction potential $V$ as in Theorem 6.3, for any $\ell \neq k \leq N$, we have*

$$|V_{\ell,k}|^2 \lesssim -\mathcal{H}_{\delta,x_\ell}, \tag{6.22}$$

*where the constant is uniform in $0 \leq \delta \leq 1$.*

**Proof.** We have that

$$|V_{\ell,k}|^2 = |V(x_\ell - x_k)|^2 \leq 2 \frac{\|V_1\|^2_{L^\infty_{\mathbb{T}^3}}}{|x_\ell - x_k|^2} + 2|V_2(x_\ell - x_k)|^2.$$

For the Hardy inequality in Lemma 6.6, we have

$$2 \frac{\|V_1\|^2_{L^\infty_{\mathbb{T}^3}}}{|x_\ell - x_k|^2} \lesssim -\mathcal{H}_{\delta,x_\ell}$$

---

6.2. Observe that in [BCD11] Hardy's inequality is proven on the whole space but the proof on $\mathbb{T}^n$ is identical.

uniformly in $0 \leq \delta \leq 1$. Using the fact that

$$\int |V_2(x_\ell - x_k)\varphi(x_k)|^2 dx_k \lesssim \|V_2\|_{L^3_{\mathbb{T}^3}}^2 \|\varphi\|_{L^6_{\mathbb{T}^3}}^2 \lesssim \|V_2\|_{L^3_{\mathbb{T}^3}}^2 \|\sqrt{-\mathcal{H}_\delta}\varphi\|_{L^2_{\mathbb{T}^3}}^2$$

where we have exploited inequality (2.20) (in the case $s=0$) in Lemma 2.6, and Corollary 2.7 for $0 \leq \delta \leq 1$. □

We now introduce the following operators. For any $\rho^{n+1} \in \mathcal{L}^1_{n+1}$, the following expressions are well-defined:

$$\mathcal{I}_{n,\ell,V}(\rho^{n+1}) := \mathrm{Tr}_{n+1}(V_{\ell,n+1}\rho^{n+1}), \tag{6.23}$$
$$\mathcal{J}_{n,\ell,V}(\rho^{n+1}) := \mathrm{Tr}_{n+1}(\rho^{n+1}V_{\ell,n+1}), \tag{6.24}$$

for $n \in \mathbb{N}$, $l \leq n+1$, with which we can write

$$\mathrm{Tr}_{n+1}([V_{\ell,n+1},\rho^{n+1}]) = \mathcal{I}_{n,\ell,V}(\rho^{n+1}) - \mathcal{J}_{n,\ell,V}(\rho^{n+1}).$$

**Lemma 6.8.** *The expressions (6.23) and (6.24) are linear bounded operators from $\mathcal{W}^{1,1}_{n+1,\delta}$ to $\mathcal{L}^1_n$ uniformly in $\delta \in [0,1]$ for V of the form (6.18).*

**Proof.** The result follows from (6.22). We give here the details on for the operator $\mathcal{I}_{n,\ell,V}$, the proof for the operator $\mathcal{J}_{n,\ell,V}$ being completely analogous. For any $\rho^{n+1} \in \mathcal{W}^{1,1}_{n+1,\delta}$, we have

$$\mathrm{Tr}(|V_{\ell,n+1}\rho^{n+1}|) = \mathrm{Tr}\left(\sqrt{V_{\ell,n+1}(\rho^{n+1})^2 V_{\ell,n+1}}\right),$$

and thus

$$\begin{aligned}\|\mathcal{I}_{n,\ell,V}(\rho^{n+1})\|_{\mathcal{L}^1_n} &\leq \mathrm{Tr}|V_{\ell,n+1}\rho^{n+1}| \leq \mathrm{Tr}\sqrt{\rho^{n+1}V^2_{\ell,n+1}\rho^{n+1}} \lesssim \mathrm{Tr}\sqrt{\rho^{n+1}(-\mathcal{H}_{\delta,x_l})\rho^{n+1}} \\ &\lesssim \mathrm{Tr}\sqrt{(-\mathcal{H}_{\delta,x_l})^{1/2}\rho^{n+1}(-\mathcal{H}_{\delta,x_l})\rho^{n+1}(-\mathcal{H}_{\delta,x_l})^{1/2}} \\ &= \mathrm{Tr}[(-\mathcal{H}_{\delta,x_l})^{1/2}\rho^{n+1}(-\mathcal{H}_{\delta,x_l})^{1/2}] = \|\rho^{n+1}\|_{\mathcal{W}^{1,1}_{n+1}},\end{aligned}$$

where we use Remark 6.7 and (2.12) in the third and fourth inequalities. The implicit constants are all uniform in $0 \leq \delta \leq 1$. □

We study now the properties of the solutions $(\Psi_{N,\delta})$ of the linear Schrödinger equations (6.4) with finite initial energy:

$$\sup_{k \leq N \in \mathbb{N}, 0 < \delta \leq 1} \int_{\mathbb{T}^{3N}} \Psi_{N,\delta,0} \overline{(-\mathcal{H}_{\delta,x_k})\Psi_{N,\delta,0}} = \sup_{0 < \delta \leq 1} \int_{\mathbb{T}^{3N}} \Psi_{N,\delta,0} \overline{(-\mathcal{H}_{\delta,x_1})\Psi_{N,\delta,0}} < \infty, \tag{6.25}$$

and notice that for $\delta > 0$ the conservation laws (6.5) rigorously hold since the random potential is smooth on $\mathbb{T}^3$.

**Lemma 6.9.** *Under the hypotheses of Theorem 6.3 and for any $n \in \mathbb{N}$ we have*

$$\sup_{N \in \mathbb{N}, 0 < \delta \leq 1} \|\rho^n_{N,\delta}\|_{L^\infty(\mathbb{R}_+, \mathcal{W}^{1,1}_{\delta,n})} < +\infty. \tag{6.26}$$

**Proof.** By inequality (6.22) and mass conservation in (6.5), we have for any $N \in \mathbb{N}$, $\varepsilon$, $\delta > 0$, there is a constant $C_\varepsilon$ (depending only on $\varepsilon$) such that

$$\left|\frac{1}{N}\sum_{j<k\leq N}\int_{(\mathbb{T}^3)^N} V(x_j - x_k)\Psi_{N,\delta}(t)\overline{\Psi_{N,\delta}(t)}dx_{1:N}\right| \leq C_\varepsilon + \varepsilon\int_{(\mathbb{T}^3)^N}\Psi_{N,\delta}(t)\overline{(-\mathcal{H}_{\delta,x_1})\Psi_{N,\delta}(t)} \tag{6.27}$$

Now by symmetry, we have

$$\int_{(\mathbb{T}^3)^N} \Psi_{N,\delta}(t) \overline{(-\mathcal{H}_{\delta,x_1})\Psi_{N,\delta}(t)} dx_{1:N} = \frac{1}{N} \sum_{j \leq N} \int_{(\mathbb{T}^3)^N} \Psi_{N,\delta}(t) \overline{(-\mathcal{H}_{\delta,x_j})\Psi_{N,\delta}(t)} dx_{1:N}.$$

Using the conservation laws (6.5) and triangle inequality, we obtain

$$\frac{1}{N} \sum_{j \leq N} \int_{(\mathbb{T}^3)^N} \Psi_{N,\delta}(t) \overline{(-\mathcal{H}_{\delta,x_j})\Psi_{N,\delta}(t)} dx_{1:N}$$

$$\stackrel{(6.5)}{\lesssim} \left| \frac{1}{N} \sum_{j \leq N} \int_{(\mathbb{T}^3)^N} \Psi_{N,\delta,0} \overline{(-\mathcal{H}_{\delta,x_j})\Psi_{N,\delta,0}} dx_{1:N} + \frac{1}{N^2} \sum_{j < k \leq N} \int_{(\mathbb{T}^3)^N} V(x_j - x_k) \Psi_{N,\delta,0} \overline{\Psi_{N,\delta,0}} dx_{1:N} \right|$$

$$+ \left| \frac{1}{N^2} \sum_{j < k \leq N} \int_{(\mathbb{T}^3)^N} V(x_j - x_k) \Psi_{N,\delta}(t) \overline{\Psi_{N,\delta}(t)} dx_{1:N} \right|$$

$$\stackrel{(6.27)}{\lesssim} 2 \left| \int_{(\mathbb{T}^3)^N} \Psi_{N,\delta,0} \overline{(-\mathcal{H}_{\delta,x_1})\Psi_{N,\delta,0}} dx_{1:N} \right| + \varepsilon \int_{(\mathbb{T}^3)^N} \Psi_{N,\delta}(t) \overline{(-\mathcal{H}_{\delta,x_1})\Psi_{N,\delta}(t)} dx_{1:N} + C_\varepsilon.$$

By choosing $\varepsilon$ small enough, and using the symmetry of $\Psi_{N,\delta}(t)$ for $t \geq 0$ again and assumption on the initial condition, this implies that

$$\sup_{N \in \mathbb{N}, 0 < \delta \leq 1} \int_{(\mathbb{T}^3)^N} \Psi_{N,\delta}(t) \overline{(-\mathcal{H}_{\delta,x_1})\Psi_{N,\delta}(t)} dx_{1:N} \lesssim \sup_{N \in \mathbb{N}, 0 < \delta \leq 1} \int_{(\mathbb{T}^3)^N} \Psi_{N,\delta,0} \overline{(-\mathcal{H}_{\delta,x_1})\Psi_{N,\delta,0}} dx_{1:N}$$
$$< +\infty \qquad (6.28)$$

Thus

$$\begin{aligned}
\|\rho_{N,\delta}^n(t)\|_{\mathcal{W}_{\delta,n}^{1,1}} &= n\text{Tr}[(-\mathcal{H}_{\delta,x_1})^{1/2} \rho_{N,\delta}^n (-\mathcal{H}_{\delta,x_1})^{1/2}] \\
&= n\text{Tr}[(-\mathcal{H}_{\delta,x_1})^{1/2} \text{Tr}_{n+1:N}(\Pi_{\Psi_{N,\delta}}) (-\mathcal{H}_{\delta,x_1})^{1/2}] \\
&= n\text{Tr}[(-\mathcal{H}_{\delta,x_1})^{1/2} \Pi_{\Psi_{N,\delta}} (-\mathcal{H}_{\delta,x_1})^{1/2}] \\
&= n \int_{(\mathbb{T}^3)^N} \Psi_{N,\delta}(t) \overline{(-\mathcal{H}_{\delta,x_1})\Psi_{N,\delta}(t)} dx_{1:N} < \infty,
\end{aligned}$$

which by inequality (6.28) is uniformly bounded in $N \in \mathbb{N}$, $\delta > 0$ and $t \in \mathbb{R}_+$. □

**Remark 6.10.** Thanks to Lemma 6.9 we can extract subsequences of converging density matrices with $\mathcal{W}^{1,1}$ regularity. We also notice that the weak* limit $\rho$ of density matrices $(\rho_n)_{n \geq 0}$ on $L_s^2((\mathbb{T}^3)^n)$ is always self-adjoint and non-negative. Indeed, self-adjointness follows immediately from testing with projection $f \otimes \bar{g}$ for arbitrary $f, g \in L_s^2((\mathbb{T}^3)^n)$:

$$(f, \rho_n g) = \text{Tr}(\rho_n \cdot f \otimes \bar{g}) \to \text{Tr}(\rho \cdot f \otimes \bar{g}) = (\rho f, g),$$

and similarly $(\rho_n f, g) \to (\rho f, g)$. Non-negativity of $\rho$ follows similarly by testing against $f \otimes \bar{f}$ for arbitrary $f \in L_s^2((\mathbb{T}^3)^n)$. On the other hand, by the cyclicity of trace, any solution to (6.15) satisfies

$$\text{Tr}[\rho_\infty^n(t)] = \text{Tr}[\rho_\infty^n(0)] = 1.$$

**Lemma 6.11.** *Under the hypotheses of Theorem 6.3, there are subsequences $N_k \to +\infty$ and $\delta_k \to 0^+$ together with a sequence $(\rho_\infty^n)_{n \in \mathbb{N}}$ such that, for any $n \in \mathbb{N}$, it holds that*

$$\rho_\infty^n \in L^\infty(\mathbb{R}_+, \mathcal{W}_n^{1,1})$$

*and*

$$\rho_{N_k,\delta_k}^n \to \rho_\infty^n, \quad (-\mathcal{H}_{\delta_k,x_l})^{1/2} \rho_{N_k,\delta_k}^n (-\mathcal{H}_{\delta_k,x_l})^{1/2} \to (-\mathcal{H}_{x_l})^{1/2} \rho_\infty^n (-\mathcal{H}_{x_l})^{1/2}, \quad l = 1, \ldots, n,$$

*weakly\* in $L^\infty(\mathbb{R}_+, \mathcal{L}_n^1)$.*

**Proof.** Since

$$\sup_{t \in \mathbb{R}_+} \|\rho_{N,\delta}^n(t)\|_{\mathcal{L}_n^1} = 1, \quad \forall n \leq N, \delta > 0,$$

there exist subsequences $N_k \to +\infty$, $\delta_k \to 0^+$ and a family of operators $\rho_\infty^n \in L^\infty(\mathbb{R}_+, \mathcal{L}_n^1)$ indexed by $n \in \mathbb{N}$, s.t.

$$\rho_k^n := \rho_{N_k, \delta_k}^n \to \rho_\infty^n, \quad \text{weakly* in } L^\infty(\mathbb{R}_+, \mathcal{L}_n^1), \quad \forall n \in \mathbb{N},$$

as $k \to \infty$. By Lemma 6.9, we have by defining for $l \leq n \in \mathbb{N}$

$$\tilde{\rho}_{k,\ell}^n := (-\mathcal{H}_{\delta_k, x_l})^{1/2} \rho_k^n (-\mathcal{H}_{\delta_k, x_l})^{1/2} \Longrightarrow \sup_{k \in \mathbb{N}, \ell \leq n} \|\tilde{\rho}_{k,l}^n\|_{\mathcal{L}_n^1} < +\infty.$$

Thus up to further taking a subsequence, there exists a family of operators $\tilde{\rho}_{\infty,l}^n \in L^\infty(\mathbb{R}_+, \mathcal{L}_n^1)$ in $n \in \mathbb{N}$, such that

$$(-\mathcal{H}_{\delta_k, x_l})^{1/2} \rho_k^n (-\mathcal{H}_{\delta_k, x_l})^{1/2} = \tilde{\rho}_{k,l}^n \to \tilde{\rho}_{\infty,l}^n \tag{6.29}$$

weakly\* in $L^\infty(\mathbb{R}_+, \mathcal{L}_n^1)$ as $k \to +\infty$, which implies according to Lemma C.1 and the convergence of $(-\mathcal{H}_\delta)^{-1/2} \to (-\mathcal{H})^{-1/2}$ in $\mathcal{L}(L_{\mathbb{T}^3}^2; L_{\mathbb{T}^3}^2)$ (see Lemma 2.13) that

$$\rho_k^n = (-\mathcal{H}_{\delta_k, x_l})^{-1/2} \tilde{\rho}_{k,l}^n (-\mathcal{H}_{\delta_k, x_l})^{-1/2} \to (-\mathcal{H}_{x_l})^{-1/2} \tilde{\rho}_{\infty,l}^n (-\mathcal{H}_{x_l})^{-1/2}$$

weakly\* in $L^\infty(\mathbb{R}_+, \mathcal{L}_n^1)$ as $k \to +\infty$. We have therefore identified

$$(-\mathcal{H}_{x_l})^{-1/2} \tilde{\rho}_{\infty,l}^n (-\mathcal{H}_{x_l})^{-1/2} = \rho_\infty^n, \quad \forall l \leq n,$$

that is $\rho_\infty^n \in L^\infty(\mathbb{R}_+, \mathcal{W}_n^{1,1})$ and

$$(-\mathcal{H}_{\delta_k, x_l})^{1/2} \rho_k^n (-\mathcal{H}_{\delta_k, x_l})^{1/2} = \tilde{\rho}_k^n \to \mathcal{H}_{x_l}^{1/2} \rho_\infty^n \mathcal{H}_{x_l}^{1/2}$$

weakly\* in $L^\infty(\mathbb{R}_+, \mathcal{L}_n^1)$ as $k \to \infty$. □

Lemma 6.11 ensures the existence of subsequences $N_k, \delta_k$ such that $\rho_{N_k, \delta_k}^n \to \rho_\infty^n$. What remains to prove, in order to obtain Theorem 6.3, is that $(\rho_\infty^n)_{n \in \mathbb{N}}$ satisfies the hierarchy equation (6.15).

Instead of proving directly that $(\rho_\infty^n)_{n \in \mathbb{N}}$ given by Lemma 6.11 satisfies equation (6.15), we prove that there exists an orthonormal basis $\{e_n\}_{n \in \mathbb{N}}$ of $L^2(\mathbb{T}^3)$, such that for any $n \in \mathbb{N}$, $0 \leq t_1 \leq t_2$, $I_n = (i_1, \ldots, i_n)$ and $\mathcal{J}_n = (j_1, \ldots, j_n) \in \mathbb{N}^n$, we have

$$\int_{t_1}^{t_2} \text{Tr}(O^{I_n, \mathcal{J}_n} \rho_\infty^n(t)) dt$$

$$= \int_{t_1}^{t_2} \text{Tr}\left( O^{I_n, \mathcal{J}_n} \left( \prod_{k=1}^n \exp\left( \frac{it}{2}(\mathcal{H}_{x_k} - \mathcal{H}_{y_k}) \right) \rho_\infty^n(0) \right) \right) dt$$

$$+ \int_{t_1}^{t_2} \text{Tr}\left( O^{I_n, \mathcal{J}_n} \int_0^t \left[ \prod_{k=1}^n \exp\left( \frac{i(t-s)}{2}(\mathcal{H}_{x_k} - \mathcal{H}_{y_k}) \right) \left( \sum_{\ell=1}^n \mathcal{I}_{n,\ell,V}(\rho_\infty^{n+1}(s)) \right) ds \right] \right) dt$$

$$- \int_{t_1}^{t_2} \text{Tr}\left( O^{I_n, \mathcal{J}_n} \int_0^t \left[ \prod_{k=1}^n \exp\left( \frac{i(t-s)}{2}(\mathcal{H}_{x_k} - \mathcal{H}_{y_k}) \right) \left( \sum_{\ell=1}^n \mathcal{J}_{n,\ell,V}(\rho_\infty^{n+1}(s)) \right) ds \right] \right) dt, \tag{6.30}$$

where $O^{I_n, \mathcal{J}_n}$ is the bounded range operator associated with the kernel

$$O^{I_n, \mathcal{J}_n}(x_{1:n}, y_{1:n}) = e_{I_n}(x_{1:n}) e_{\mathcal{J}_n}(y_{1:n}), \tag{6.31}$$

where $e_{I_n}$ is the *symmetric tensor product* of $(e_i)_{i \in I_n}$.

**Notation 6.12.** *If a sequence $(\rho_\infty^n)_{n \in \mathbb{N}}$ satisfies equations (6.30) for any $t_1, t_2$ and $O^{I_n, \mathcal{J}_n}$ described above, we say that $(\rho_\infty^n)_{n \in \mathbb{N}}$ is a weak mild solution to the hierarchy equation (6.15) with respect to the basis $\{e_n\}_{n \in \mathbb{N}}$.*

**Lemma 6.13.** *Suppose that the sequence $(\rho_\infty^n)_{n \in \mathbb{N}}$ is a weak mild solution to the hierarchy (6.15), and $\rho_\infty^n \in L^\infty(\mathbb{R}_+, \mathcal{W}_n^{1,1})$, then $(\rho_\infty^n)_{n \in \mathbb{N}}$ is a $\mathcal{W}^{1,1}$ solution to the to the hierarchy equation.*

**Proof.** Let $(e_i)_{i \geq 1}$ be a fixed basis of $L^2(\mathbb{T}^3)$, and let $S \subset L^2(\mathbb{T}^3)$ be the set of (finite) linear combination of the elements in $(e_i)_{i \geq 1}$, we define

$$\mathcal{D}_n := \{ f \Pi_{g_1 \otimes \cdots \otimes g_n} | f \in L^1(\mathbb{R}_+) \text{ simple functions}, g_1, \ldots, g_n \in S \} \subset L^1(\mathbb{R}_+, \mathcal{K}_n),$$

that is the family of functions from $\mathbb{R}_+$ to $\mathcal{K}_n$ of the form

$$\mathbb{R}_+ \ni t \mapsto A(t) = \sum_{k=1}^{\ell} f_k \mathbb{I}_{[t_k, t_k')}(t) \Pi_{g_1^k \otimes \cdots \otimes g_n^k}, \quad \text{where} \quad g_j^k = \sum_{m \in \mathbb{N}} h_{j,m}^k e_m,$$

for some $h_{j,m}^k \in \mathbb{R}$ ($h_{j,m}^k \neq 0$ for finite $m$, with any fixed $k, j$). It is straightforward to observe that $\mathcal{D}_n$ *separates points of* $L^\infty(\mathbb{R}_+, \mathcal{L}_n^1)$ (since $\mathcal{D}_n$ is dense in the strong pre-dual $L^1(\mathbb{R}_+, \mathcal{K}_n)$ of $L^\infty(\mathbb{R}_+, \mathcal{L}_n^1)$), and we can write explicitly for any $B \in L^\infty(\mathbb{R}_+, \mathcal{L}_n^1)$ that

$$\int_0^{+\infty} \text{Tr}(A(t) B(t)) dt = \sum_{k=1}^{\ell} \sum_{I_n, \mathcal{J}_n} f_k h_{1,i_1}^k \cdots h_{n,i_n}^k h_{1,j_1}^k \cdots h_{n,j_n}^k \int_{t_k}^{t_k'} \text{Tr}(O^{I_n, \mathcal{J}_n} B(s)) ds,$$

where $O^{I_n, \mathcal{J}_n} = e_{I_N}(x_{1:n}) e_{\mathcal{J}_N}(y_{1:n})$ as in (6.31). Therefore, we have the equivalence

$$L^\infty(\mathbb{R}_+, \mathcal{L}_n^1) \ni B = 0 \iff \int_0^{+\infty} \text{Tr}(A(t) B(t)) dt = 0, \forall A \in \mathcal{D}_n$$

$$\iff \int_{t_k}^{t_k'} \text{Tr}(O^{I_n, \mathcal{J}_n} B(s)) ds = 0, \forall t_k \leq t_{k'}, i_{1:n}, j_{1:n} \in \mathbb{N}^n, \tag{6.32}$$

since $\mathcal{D}_n \subset L^1(\mathbb{R}_+, \mathcal{K}_n)$ separates points of $L^\infty(\mathbb{R}_+, \mathcal{L}_n^1)$. If we take

$$B(t) = \rho_\infty^n(t) - \mathcal{U}_t^n(\rho_{\infty,0}^n) - \sum_{k=1}^n \int_0^t \mathcal{U}_{t-s}^n(\text{Tr}_{n+1}[V_{k,n+1}, \rho_\infty^{n+1}(s)]) ds, \tag{6.33}$$

which is an element of $L_{\text{loc}}^\infty(\mathbb{R}_+, \mathcal{L}^1)$ due to the fact that, by Lemma 6.8, $s \mapsto \text{Tr}_{n+1}[V_{k,n+1}, \rho_\infty^{n+1}(s)] \in L^\infty(\mathbb{R}_+, \mathcal{L}_n^1)$ whenever $\rho_\infty^{n+1} \in L^\infty(\mathbb{R}_+, \mathcal{W}_n^{1,1})$, (6.32) indicates that conditions (6.30) is equivalent to the statement that $(\rho_\infty^n)_{n \in \mathbb{N}}$ is a mild solution to the hierarchy equation (6.15). $\square$

**Proof of Theorem 6.3.** We recall from the discussions above, it suffices to show that $(\rho_\infty^n)_{n \in \mathbb{N}}$ satisfies the hierarchy equation. With the help of Lemma 6.13, this reduces to show that $(\rho_\infty^n)_{n \in \mathbb{N}}$ is a *weak mild solution* to the hierarchy equation (6.15) with the $L_{\mathbb{T}^3}^2$-basis $\{f_n\}_{n \in \mathbb{N}}$ given by the *orthonormal eigenfunctions of the Anderson Hamiltonian $\mathcal{H}$*.

We fix $I = (i_1, \ldots, i_n)$, $\mathcal{J} = (j_1, \ldots, j_n) \in \mathbb{N}^n$, $0 \leq t_1 < t_2$, and consider the corresponding compact operator $O^{I,\mathcal{J}}$ given by formula (6.31) where we chose specifically $(e_i)_{i \geq 0}$ to be the eigenfunctions $(f_i)_{i \geq 0}$ of $\mathcal{H}$. We denote as in the proof of Lemma 6.11

$$\rho_k^n(t) := \rho_{k,t}^n := \rho_{N_k,\delta_k}^n(t), \quad \mathcal{U}_{k,t}^n := \mathcal{U}_k^n(t) := \mathcal{U}_{N_k,\delta_k,t}^n. \tag{6.34}$$

If we write the finite hierarchy equation (6.7) for $\rho_k^n$ in the weak mild formulation, we obtain the equation

$$\begin{aligned}
\int_{t_1}^{t_2} \text{Tr}(O^{I,\mathcal{J}} \rho_{k,t}^n) \, dt &= \int_{t_1}^{t_2} \text{Tr}(O^{I,\mathcal{J}} \mathcal{U}_{k,t}^n(\rho_{k,0}^n)) \, dt \\
&\quad + \frac{1}{N} \sum_{j < \ell \leq n} \int_{t_1}^{t_2} \text{Tr}\left(O^{I,\mathcal{J}} \int_0^t \mathcal{U}_{k,t-s}^n([V_{j,\ell}, \rho_{k,s}^{n+1}]) \, ds\right) dt \\
&\quad + \frac{N-n}{N} \sum_{\ell=1}^n \int_{t_1}^{t_2} \text{Tr}\left(O^{I,\mathcal{J}} \int_0^t \mathcal{U}_{k,t-s}^n(\text{Tr}_{n+1}[V_{\ell,n+1}, \rho_{k,s}^{n+1}]) \, ds\right) dt.
\end{aligned} \tag{6.35}$$

We show in the rest of the proof that each term in (6.35) converges to the corresponding term in (6.30), proving that $(\rho_\infty^n)_{n \in \mathbb{N}}$ is a *weak mild solution* to the BBGKY hierarchy. The L.H.S. of (6.35)

$$\int_{t_1}^{t_2} \text{Tr}(O^{I,\mathcal{J}} \rho_k^n(t)) \, dt \to \int_{t_1}^{t_2} \text{Tr}(O^{I,\mathcal{J}} \rho_\infty^n(t)) \, dt,$$

by definition of weak* convergence in $L^\infty(\mathbb{R}_+, \mathcal{L}_n^1)$. For the first term on the right of (6.35), we note that

$$\|O^{I,\mathcal{J}} \mathcal{U}_{k,t}^n(\rho_{k,0}^n)\|_{\mathcal{L}_n^1} \leq \|O^{I,\mathcal{J}}\|_{\mathcal{L}_n^\infty} \|\mathcal{U}_{k,t}^n(\rho_{k,0}^n)\|_{\mathcal{L}_n^1} \leq \|O^{I,\mathcal{J}}\|_{\mathcal{L}_n^\infty},$$

where we used that $\|\mathcal{U}_{k,t}^n(\rho_{k,0}^n)\|_{\mathcal{L}_n^1} \leq \|\rho_{k,0}^n\|_{\mathcal{L}_n^1} = \text{Tr}(\rho_{k,0}^n) = 1$, since (by Lemma 6.15 below) $\mathcal{U}_k^n$ is unitary. Thus if we are able to prove the convergence of $\text{Tr}(O^{I,\mathcal{J}} \mathcal{U}_{k,t}^n(\rho_{k,0}^n))$ for every fixed $t \in [t_1, t_2]$, the convergence

$$\int_{t_1}^{t_2} \text{Tr}(O^{I,\mathcal{J}} \mathcal{U}_{k,t}^n(\rho_{k,0}^n)) \, dt \to \int_{t_1}^{t_2} \text{Tr}(O^{I,\mathcal{J}} \mathcal{U}_t^n(\rho_{\infty,0}^n)) \, dt \tag{6.36}$$

follows from Lebesgue dominated convergence. Since in kernel form

$$\mathcal{U}_{k,t}^n(\rho_{k,0}^n) = \exp\left(\frac{it}{2} \sum_{h=1}^n \mathcal{H}_{\delta,x_h} - \frac{it}{2} \sum_{h=1}^n \mathcal{H}_{\delta,y_h}\right) \rho_{k,0}^n(x_{1:h}, y_{1:h}), \tag{6.37}$$

and

$$\exp\left(\frac{it}{2} \sum_{h=1}^n \mathcal{H}_{\delta,x_h}\right) \to \exp\left(\frac{it}{2} \sum_{h=1}^n \mathcal{H}_{x_h}\right), \tag{6.38}$$

strongly on $L_s^2((\mathbb{T}^3)^n)$, as $\mathcal{H}_\delta$ converges to $\mathcal{H}$ in the norm resolvent topology as $\delta \to 0^+$, we derive the desired convergence (6.36) by Lemma C.1 and the compactness of $O^{I,\mathcal{J}}$.

We consider now the second term on the right of (6.35). Using Lebesgue dominated convergence again and the fact that, by Lemma 6.8, for every $j < n+1$,

$$\sup_{0 \leq s \leq t \leq t_2} \|O^{I,\mathcal{J}} \mathcal{U}_{k,t-s}^n([V_{j,n+1}, \rho_k^{n+1}(s)])\|_{\mathcal{L}_n^1} \lesssim 1,$$

it is enough to prove the convergence of $\frac{1}{N}\mathrm{Tr}(O^{I,\mathcal{J}}\mathcal{U}^n_{k,t-s}([V_{j,n+1},\rho^{n+1}_{N_k,\delta_k}(s)])) \to 0$ for any $s < t$ fixed. On the other hand, using again Lemma 6.8,

$$\left\| \frac{1}{N_k} O^{I,\mathcal{J}} \mathcal{U}^n_{k,t-s}[V_{i,j},\rho^n_k(s)] \right\|_{\mathcal{L}^1_n} \lesssim \frac{1}{N_k}\|[V_{i,j},\rho^n_k(s)]\|_{\mathcal{L}^1_n}ds \lesssim \frac{1}{N_k}\|\rho^n_k\|_{L^\infty(\mathbb{R}_+,\mathcal{W}^{1,1}_{n,\delta_k})} \to 0,$$

as $k \to +\infty$, since, by Lemma 6.11, $\|\rho^n_k\|_{L^\infty(\mathbb{R}_+,\mathcal{W}^{1,1}_{n,\delta_k})}$ is uniformly bounded in $k \in \mathbb{N}$.

We turn to the last term on the right of (6.35). Exploiting the argument based on the uniform boundedness in $\mathcal{L}^1_n$ of the involved operators and the Lebesgue dominated convergence theorem as above, it is enough to prove that for any $0 \leq s \leq t \leq t_2$ and $\ell = 1, \ldots, n$, we have

$$|\mathrm{Tr}(O^{I,\mathcal{J}}\mathcal{U}^n_{k,t-s}(\mathrm{Tr}_{n+1}[V_{\ell,n+1},\rho^{n+1}_k(s)])) - \mathrm{Tr}(O^{I,\mathcal{J}}\mathcal{U}^n_{t-s}(\mathrm{Tr}_{n+1}[V_{\ell,n+1},\rho^{n+1}_\infty(s)]))| \to 0, \quad (6.39)$$

as $k \to \infty$. First, we can replace $\mathcal{U}^n_{k,t-s}$ in the first term in (6.39) by $\mathcal{U}^n_{t-s}$, since

$$|\mathrm{Tr}((O^{I,\mathcal{J}}\mathcal{U}^n_{k,t-s}(\mathrm{Tr}_{n+1}[V_{\ell,n+1},\rho^{n+1}_k(s)]))) - \mathrm{Tr}(O^{I,\mathcal{J}}\mathcal{U}^n_{t-s}(\mathrm{Tr}_{n+1}[V_{\ell,n+1},\rho^{n+1}_k(s)]))| \to 0, \quad (6.40)$$

as $k \to \infty$. Indeed the convergence (6.40) follows from: the cyclicity of trace; the fact that

$$\left\| \exp\left(-\frac{it}{2}\sum_{h=1}^n \mathcal{H}_{\delta,x_h}\right) O^{I,\mathcal{J}} \exp\left(\frac{it}{2}\sum_{h=1}^n \mathcal{H}_{\delta,x_h}\right) - \exp\left(-\frac{it}{2}\sum_{h=1}^n \mathcal{H}_{x_h}\right) O^{I,\mathcal{J}} \exp\left(\frac{it}{2}\sum_{h=1}^n \mathcal{H}_{x_h}\right) \right\|_{\mathcal{L}^\infty_n} \to 0,$$

as $\delta \to 0^+$ as a consequence of Lemma C.1 and the fact that $e^{\pm it\sum_{h=1}^n \mathcal{H}_{\delta,x_h}} - e^{\pm it\sum_{h=1}^n \mathcal{H}_{x_h}} \to 0$ strongly and $O^{I,\mathcal{J}}$ is a fixed compact operator of the form (6.31); and the fact that, by Lemma 6.8 and Lemma 6.9, for any $\ell = 1, \ldots, n$,

$$\sup_{k\in\mathbb{N},s\in\mathbb{R}_+} \|\mathrm{Tr}_{n+1}[V_{\ell,n+1},\rho^{n+1}_k(s)]\|_{\mathcal{L}^1_n} \leq \sup_{k\in\mathbb{N},s\in\mathbb{R}_+} \|[V_{\ell,n+1},\rho^{n+1}_k(s)]\|_{\mathcal{L}^1_{n+1}} \lesssim \sup_{k\in\mathbb{N}} \|\rho^{n+1}_k\|_{\mathcal{W}^{1,1}_{\delta_k,n+1}}$$

By using that $O^{I,\mathcal{J}}$ is a projection onto the tensor products of eigenfunctions of $\mathcal{H}$, which commutes with the unitary groups $\exp(\pm it\sum_{h=1}^n \mathcal{H}_{x_h})$ in view of (6.40), (6.39) can be reduced to showing

$$\lim_{k\to\infty} |\mathrm{Tr}(O^{I,\mathcal{J}}\mathrm{Tr}_{n+1}[V_{\ell,n+1},\rho^{n+1}_k(s) - \rho^{n+1}_\infty(s)])| = 0 \quad (6.41)$$

for each $s \in \mathbb{R}_+$. In order to prove the limit (6.41), we split the potential $V_{\ell,n+1}$ into three parts, namely

$$V_{\ell,n+1} = \tilde{V}^1_{\ell,n+1}\eta^\varepsilon_{\ell,n+1} + \tilde{V}^1_{\ell,n+1}(1-\eta^\varepsilon_{\ell,n+1}) + V^2_{\ell,n+1} \quad (6.42)$$

where

$$\tilde{V}^1_{\ell,n+1} := \frac{V^1(x_\ell - x_{n+1})}{|x_\ell - x_{n+1}|}, \quad V^2_{\ell,n+1} := V^2(x_\ell - x_{n+1}),$$

where $V^1, V^2$ are defined in (6.18), and $\eta^\varepsilon_{\ell,n+1} := \eta^\varepsilon(x_\ell - x_{n+1})$ with some smooth and radially symmetric bump-function $\eta^\varepsilon$ of radius $\varepsilon > 0$ such that

$$\eta^\varepsilon(x) = 1 \text{ if } |x| \leq \frac{\varepsilon}{4}, \quad \eta^\varepsilon(x) = 0 \text{ if } |x| \geq \frac{\varepsilon}{2}.$$

Since $\tilde{V}^1_{\ell,n+1}\eta^\varepsilon_{\ell,n+1} = (\tilde{V}^1\eta^\varepsilon)_{\ell,n+1}$ keeps the singularity of $\tilde{V}^1$ at the origin, while we can make it arbitrarily small as $\varepsilon \to 0^+$ in the sense that

$$\lim_{\varepsilon\to 0^+} |\mathrm{Tr}(O^{I,\mathcal{J}}\mathrm{Tr}_{n+1}[(\tilde{V}^1\eta^\varepsilon)_{\ell,n+1},\rho^{n+1}_\infty])| + \sup_{k\in\mathbb{N}} |\mathrm{Tr}(O^{I,\mathcal{J}}\mathrm{Tr}_{n+1}[(\tilde{V}^1\eta^\varepsilon)_{\ell,n+1},\rho^{n+1}_k])| = 0. \quad (6.43)$$

Since $\text{Tr}_{n+1}[(\tilde{V}^1\eta^\varepsilon)_{\ell,n+1}, \rho_k^{n+1}] = \mathcal{I}_{n,\ell,\tilde{V}^1\eta^\varepsilon}(\rho_k^{n+1}) - \mathcal{J}_{n,\ell,\tilde{V}^1\eta^\varepsilon}(\rho_k^{n+1})$, and anologously for $\rho_\infty^{n+1}$, we give the proof of the convergence (6.43) when $\text{Tr}_{n+1}[(\tilde{V}^1\eta^\varepsilon)_{\ell,n+1}, \rho_k^{n+1}]$ is replaced by $\mathcal{I}_{n,\ell,\tilde{V}^1\eta^\varepsilon}(\rho_k^{n+1})$, the proof of the term involving the expression $-\mathcal{J}_{n,\ell,\tilde{V}^1\eta^\varepsilon}(\rho_k^{n+1})$ being completely analogous.

By cyclicity of the trace, we can write

$$\text{Tr}(O^{I,J}\mathcal{I}_{n,\ell,\tilde{V}^1\eta^\varepsilon}(\rho_k^{n+1}))$$
$$= \text{Tr}[(\tilde{V}^1)_{\ell,n+1}(-\mathcal{H}_{\delta_k,x_{n+1}})^{-1/2}\tilde{\rho}_k^{n+1}(-\mathcal{H}_{\delta_k,x_{n+1}})^{-1/2}(O^{I,J}\otimes \text{Id}_{n+1})(\eta^\varepsilon)_{\ell,n+1}]$$
$$= \text{Tr}[(\tilde{V}^1)_{\ell,n+1}(-\mathcal{H}_{\delta_k,x_{n+1}})^{-1/2}\tilde{\rho}_k^{n+1}(O^{I,J}\otimes(-\mathcal{H}_{\delta_k,x_{n+1}})^{-1/2})(\eta^\varepsilon)_{\ell,n+1}],$$

where

$$\tilde{\rho}_k^{n+1} = (-\mathcal{H}_{\delta_k,x_{n+1}})^{1/2}\rho_k^{n+1}(-\mathcal{H}_{\delta_k,x_{n+1}})^{1/2}$$

as in (6.29), and we notice that

$$(-\mathcal{H}_{\delta_k,x_{n+1}})^{-1/2}(O^{I,J}\otimes \text{Id}_{n+1}) = O^{I,J}\otimes(-\mathcal{H}_{\delta_k,x_{n+1}})^{-1/2}.$$

and are uniformly bounded (in $k\in\mathbb{N}$) compact operator in $\mathcal{K}_n(L^2(\mathbb{T}^{3(n+1)}))$.

This means that

$$|\text{Tr}[(\tilde{V}^1)_{\ell,n+1}(-\mathcal{H}_{\delta_k,x_{n+1}})^{-1/2}\tilde{\rho}_k^{n+1}(O^{I,J}\otimes(-\mathcal{H}_{\delta_k,x_{n+1}})^{-1/2})(\eta^\varepsilon)_{\ell,n+1}]|$$
$$\leq \|(\tilde{V}^1)_{\ell,n+1}(-\mathcal{H}_{\delta_k,x_{n+1}})^{-1/2}\|_{\mathcal{L}_{n+1}^\infty}\|\tilde{\rho}_k^{n+1}\|_{\mathcal{L}_{n+1}^1}\|(O^{I,J}\otimes(-\mathcal{H}_{\delta_k,x_{n+1}})^{-1/2})(\eta^\varepsilon)_{\ell,n+1}\|_{\mathcal{L}_{n+1}^\infty}.$$

By Lemma 6.7 and Lemma 6.9, we have

$$\sup_{k\in\mathbb{N}}\|(\tilde{V}^1)_{\ell,n+1}(-\mathcal{H}_{\delta_k,x_{n+1}})^{-1/2}\|_{\mathcal{L}_{n+1}^\infty} < +\infty \quad , \sup_{k\in\mathbb{N}}\|\tilde{\rho}_k^{n+1}\|_{\mathcal{L}_{n+1}^1} < +\infty.$$

Furthermore we notice that

$$\sup_{k\in\mathbb{N}}\|(O^{I,J}\otimes(-\mathcal{H}_{\delta_k,x_{n+1}})^{-1/2})(\eta^\varepsilon)_{\ell,n+1}\|_{\mathcal{L}_{n+1}^\infty}$$
$$\leq \sup_{k\in\mathbb{N}}\|(O^{I,J}\otimes((-\mathcal{H}_{\delta_k,x_{n+1}})^{-1/2})(1-\Delta)^\gamma)\|_{\mathcal{L}_{n+1}^\infty}\|(1-\Delta)^{-\gamma}(\eta^\varepsilon)_{\ell,n+1}\|_{\mathcal{L}_{n+1}^\infty}$$
$$\lesssim \|(1-\Delta)^{-\gamma}(\eta^\varepsilon)_{\ell,n+1}\|_{\mathcal{L}_{n+1}^\infty} \to 0$$

as $\varepsilon\to 0$ since, $\sup_{k\in\mathbb{N}}\|(O^{I,J}\otimes((-\mathcal{H}_{\delta_k,x_{n+1}})^{-1/2})(1-\Delta)^\gamma)\|_{\mathcal{L}_{n+1}^\infty} < +\infty$ for $\gamma>0$ small enough by Corollary 2.7, $(\eta^\varepsilon)_{\ell,n+1}\to 0$ strongly, as $\varepsilon\to 0$, and $(1-\Delta)^{-\gamma}$ is a compact operator, and thus we can apply Lemma C.1 on the product $(1-\Delta)^{-\gamma}(\eta^\varepsilon)_{\ell,n+1}$. Since the same argument applies to the term with $\rho_\infty^{n+1}$, we have derived (6.43). From (6.43) follows that, there is a decreasing function $f: \mathbb{R}_+ \to \mathbb{R}_+$ such that $\lim_{\varepsilon\to 0^+} f(\varepsilon) = 0$ such that

$$|\text{Tr}(O^{I,J}\text{Tr}_{n+1}[(\tilde{V}^1\eta^\varepsilon)_{\ell,n+1}, \rho_\infty^{n+1}])| + \sup_{k\in\mathbb{N}}|\text{Tr}(O^{I,J}\text{Tr}_{n+1}[(\tilde{V}^1\eta^\varepsilon)_{\ell,n+1}, \rho_k^{n+1}])| \leq f(\varepsilon) \qquad (6.44)$$

In order to analyze the other two terms in (6.42), we fix $\varepsilon>0$ and we write

$$W^\varepsilon := \tilde{V}^1(1-\eta^\varepsilon) + V^2.$$

It is important to note that since $\tilde{V}^1(1-\eta^\varepsilon) \in L_{\mathbb{T}^3}^\infty$ for any $\varepsilon>0$, $W^\varepsilon \in L_{\mathbb{T}^3}^{3+}$. We claim that for fixed $\varepsilon>0$,

$$\lim_{k\to+\infty}|\text{Tr}(O^{I,J}\mathcal{I}_{n,\ell,W^\varepsilon}(\rho_k^{n+1}(s) - \rho_\infty^{n+1}(s)))| = 0. \qquad (6.45)$$

and an analogous statement holds for the limit $\lim_{\varepsilon \to 0} \mathcal{J}_{n,\ell,W^\varepsilon}$.

We split the difference in (6.45) by

$$\begin{aligned}
&\text{Tr}(O^{I,\mathcal{J}}\text{Tr}_{n+1}(W^\varepsilon_{\ell,n+1} \cdot (\rho^{n+1}_k - \rho^{n+1}_\infty))) \\
&= \text{Tr}(O^{I,\mathcal{J}} \otimes \text{id}_{n+1} W^\varepsilon_{\ell,n+1}(\rho^{n+1}_k - \rho^{n+1}_\infty)) \\
&= \text{Tr}(O^{I,\mathcal{J}} \otimes \text{id}_{n+1} W^\varepsilon_{\ell,n+1}((-\mathcal{H}_{\delta_k, x_{n+1}})^{-1/2} - (-\mathcal{H}_{x_{n+1}})^{-1/2})\tilde{\rho}^{n+1}_{k,n+1}(-\mathcal{H}_{\delta_k,x_{n+1}})^{-1/2}) \\
&\quad +\text{Tr}(O^{I,\mathcal{J}} \otimes \text{id}_{n+1} W^\varepsilon_{\ell,n+1}(-\mathcal{H}_{x_{n+1}})^{-1/2}(\tilde{\rho}^{n+1}_{k,n+1} - \tilde{\rho}^{n+1}_{\infty,n+1})(-\mathcal{H}_{\delta_k,x_{n+1}})^{-1/2}) \\
&\quad +\text{Tr}(O^{I,\mathcal{J}} \otimes \text{id}_{n+1} W^\varepsilon_{\ell,n+1}(-\mathcal{H}_{x_{n+1}})^{-1/2}\tilde{\rho}^{n+1}_{\infty,n+1}((-\mathcal{H}_{\delta_k,x_{n+1}})^{-1/2} - (-\mathcal{H}_{x_{n+1}})^{-1/2})) \\
&=: (\text{I})_k + (\text{II})_k + (\text{III})_k.
\end{aligned} \tag{6.46}$$

For $(\text{I})_k$, we have

$$\begin{aligned}
&|\text{Tr}(O^{I,\mathcal{J}} \otimes \text{id}_{n+1} W^\varepsilon_{\ell,n+1}((-\mathcal{H}_{\delta_k,x_{n+1}})^{-1/2} - (-\mathcal{H}_{x_{n+1}})^{-1/2})\tilde{\rho}^{n+1}_{k,n+1}(-\mathcal{H}_{\delta_k,x_{n+1}})^{-1/2})| \\
&\leq \|O^{I,\mathcal{J}} \otimes (-\mathcal{H}_{\delta_k,x_{n+1}})^{-1/2} W^\varepsilon_{\ell,n+1}((-\mathcal{H}_{\delta_k,x_{n+1}})^{-1/2} - (-\mathcal{H}_{x_{n+1}})^{-1/2})\|_{\mathcal{L}^\infty_{n+1}}\text{Tr}|\tilde{\rho}^{n+1}_{k,n+1}| \\
&\lesssim \|O^{I,\mathcal{J}} \otimes (-\mathcal{H}_{\delta_k,x_{n+1}})^{-1/2}(-\Delta+1)^\gamma\|_{\mathcal{L}^\infty_{n+1}}\|(-\Delta+1)^{-\gamma}W^\varepsilon_{\ell,n+1}((-\mathcal{H}_{\delta_k,x_{n+1}})^{-1/2} - (-\mathcal{H}_{x_{n+1}})^{-1/2})\|_{\mathcal{L}^\infty_{n+1}} \\
&\lesssim \|(-\Delta+1)^{-\gamma}W^\varepsilon_{\ell,n+1}((-\mathcal{H}_{\delta_k,x_{n+1}})^{-1/2} - (-\mathcal{H}_{x_{n+1}})^{-1/2})\|_{\mathcal{L}^\infty_{n+1}}
\end{aligned}$$

where all the implicit constants do not depend on $k$ and we used the fact that

$$\sup_{k \in \mathbb{N}} \|O^{I,\mathcal{J}} \otimes (-\mathcal{H}_{\delta_k,x_{n+1}})^{-1/2}(-\Delta+1)^\gamma\|_{\mathcal{L}^\infty_{n+1}} < +\infty$$

for $\gamma > 0$ small enough by Corollary 2.7.

The operator $W^\varepsilon_{\ell,n+1}((-\mathcal{H}_{\delta_k,x_{n+1}})^{-1/2} - (-\mathcal{H}_{x_{n+1}})^{-1/2})$ converges strongly to 0 as an operator on $L^2_s(\mathbb{T}^{3(n+1)})$ since the multiplication operator $W^\varepsilon_{\ell,n+1}$ is bounded from $L^{6-}_{\mathbb{T}^3}$ into $L^2_{\mathbb{T}^2}$ (since $W^\varepsilon_{\ell,n+1} \in L^{3+}_{\mathbb{T}^3}$) and the sequence $(-\mathcal{H}_{\delta_k,x_{n+1}})^{-1/2} - (-\mathcal{H}_{x_{n+1}})^{-1/2}$ strongly converge to 0 as sequence of operator from $L^2_{\mathbb{T}^3}$ into $L^{6-}_{\mathbb{T}^3}$ by Lemma 2.13 and Corollary 2.7. Thus, since $(-\Delta+1)^{-\gamma}$ is a compat operator whenever $\gamma > 0$, by Lemma C.1,

$$\lim_{k \to +\infty} \|(-\Delta+1)^{-\gamma}W^\varepsilon_{\ell,n+1}((-\mathcal{H}_{\delta_k,x_{n+1}})^{-1/2} - (-\mathcal{H}_{x_{n+1}})^{-1/2})\|_{\mathcal{L}^\infty_{n+1}} = 0,$$

and $\lim_{k \to +\infty}(\text{I})_k = 0$.

Similarly, we for $(\text{III})_k$, we have that

$$\begin{aligned}
&|\text{Tr}(O^{I,\mathcal{J}} \otimes \text{id}_{n+1} W^\varepsilon_{\ell,n+1}(-\mathcal{H}_{x_{n+1}})^{-1/2}\tilde{\rho}^{n+1}_{\infty,n+1}[(-\mathcal{H}_{\delta_k,x_{n+1}})^{-1/2} - (-\mathcal{H}_{x_{n+1}})^{-1/2}])| \\
&\leq \|O^{I,\mathcal{J}} \otimes [(-\mathcal{H}_{\delta_k,x_{n+1}})^{-1/2} - (-\mathcal{H}_{x_{n+1}})^{-1/2}]\|_{\mathcal{L}^\infty_{n+1}}\|W^\varepsilon_{\ell,n+1}(-\mathcal{H}_{x_{n+1}})^{-1/2}\|_{\mathcal{L}^\infty_{n+1}}\text{Tr}(\tilde{\rho}^{n+1}_{\infty,n+1}) \\
&\lesssim \|O^{I,\mathcal{J}} \otimes [(-\mathcal{H}_{\delta_k,x_{n+1}})^{-1/2} - (-\mathcal{H}_{x_{n+1}})^{-1/2}]\|_{\mathcal{L}^\infty_{n+1}} \to 0
\end{aligned}$$

as $k \to \infty$, since by Lemma 2.13, $(-\mathcal{H}_{\delta_k,x_{n+1}})^{-1/2} - (-\mathcal{H}_{x_{n+1}})^{-1/2} \to 0$ in operator norm.

Lastly, we have

$$|(\text{II})_k| = |\text{Tr}(O^{I,\mathcal{J}} \otimes (-\mathcal{H}_{\delta_k,x_{n+1}})^{-1/2} W^\varepsilon_{\ell,n+1}(-\mathcal{H}_{x_{n+1}})^{-1/2}(\tilde{\rho}^{n+1}_{k,n+1} - \tilde{\rho}^{n+1}_{\infty,n+1}))|,$$

where $(O^{I,\mathcal{J}} \otimes (-\mathcal{H}_{\delta_k,x_{n+1}})^{-1/2} W^\varepsilon_{\ell,n+1}(-\mathcal{H}_{x_{n+1}})^{-1/2})_{k \geq 0}$ is a convergent sequence in $(\mathcal{K}, \|\cdot\|_{\mathcal{L}(L^2_{\mathbb{T}^3})})$ and $\tilde{\rho}^{n+1}_{k,n+1} \to \tilde{\rho}^{n+1}_{\infty,n+1}$ weakly-* as $k \to \infty$ in $L^\infty(\mathbb{R}_+; \mathcal{L}^1_{n+1})$, so that the last limit in Lemma C.1 applies, and $\lim_{k \to +\infty}|(\text{II})_k| = 0$.

Putting all the previous limit together we obtain

$$\limsup_{k \to +\infty}|\text{Tr}(O^{I,\mathcal{J}}\text{Tr}_{n+1}[V_{\ell,n+1}, \rho^{n+1}_k(s) - \rho^{n+1}_\infty(s)])| \leq f(\varepsilon), \tag{6.47}$$

where the function $f$ is the decreasing function introduced in inequality (6.44).

Since $\varepsilon > 0$ can be chosen arbitrarily, and $\lim_{\varepsilon \to 0^+} f(\varepsilon) = 0$ we have obtained (6.39). This means that every term of (6.35) converges to the corresponding term in (6.30) and thus, by Lemma 6.13, Theorem 6.3 is proved. □

## 6.3 Proof of uniqueness: Theorem 6.4

We can now state rigorously our results for the derivation of the Anderson NLS with Hartree nonlinearity from a Bosonic many-body systems with a bounded interaction potential $V \in L^\infty_{\mathbb{T}^3}$. We recall that under the setting of Theorem 6.4, the mollification parameter $\delta = \delta(N) \to 0^+$ in the definition of the mollified Anderson Hamiltonian.

**Lemma 6.14. (Spohn, [Spo80])** *For any $n \geq 1$ we have the following bound*

$$\left\| \sum_{k \leq n} \mathrm{Tr}_{n+1}([V_{k,n+1}, \cdot]) \right\|_{\mathcal{L}(\mathcal{L}^1_{n+1}; \mathcal{L}^1_n)} \leq 2n \|V\|_{L^\infty_{\mathbb{T}^3}}. \tag{6.48}$$

**Proof.** The assertion directly follows from

$$\left\| \sum_{k \leq n} \mathrm{Tr}_{n+1}([V_{k,n+1}, \rho^{n+1}]) \right\|_{\mathcal{L}^1_n} \leq n \mathrm{Tr}|\mathrm{Tr}_{n+1}([V_{k,n+1}, \rho^{n+1}])| \leq 2n \|V\|_{L^\infty_{\mathbb{T}^3}} \|\rho^{n+1}\|_{\mathcal{L}^1_{n+1}}$$

for any $\rho^{n+1} \in \mathcal{L}^1_{n+1}$. □

**Lemma 6.15.** *The evolution $\mathbb{R} \ni t \to \mathcal{U}^n_t \in \mathcal{L}(\mathcal{L}^1_n)$ is a unitary group and furthermore for each $t \in \mathbb{R}$, $\mathcal{U}^n_t$ is an isometry on $\mathcal{L}^1_n$.*

The proof of Lemma 6.15 is postponed to Appendix C. One ingredient for the uniqueness of the solution of BBGKY hierarchy (6.9) is that from the Heisenberg picture, the unitary group $(\mathcal{U}^n_t)_{t \in \mathbb{R}}$ is uniformly bounded on $\mathcal{L}^1_n$. As Lemma 6.15 shows, this property is satisfied by $\mathcal{U}^n_t$ generated by the Anderson Hamiltonian $\mathcal{H}$ in fact also in much higher generality.

For this reason, it is not surprising that the proof due to Spohn [Spo80] where the interaction potential $V$ is assumed to be *bounded* works almost verbatim in our case.

**Proof. (Theorem 6.4)** Suppose that $(\rho^n_{\infty,1}(t)), (\rho^n_{\infty,2}(t))$ are two sequences of density matrices that are mild solutions to the hierarchy equation (6.15) with the same initial consition $(\rho^n_{\infty,0})_{n \in \mathbb{N}}$. This means that, by the linearity of the hierarchy equation and applying inequality (6.48) that

$$\|\rho^n_{\infty,1}(t) - \rho^n_{\infty,2}(t)\|_{\mathcal{L}^1_n} \leq \int_0^t \sum_{k=1}^n \|\mathcal{U}^n_{t-s}(\mathrm{Tr}_{n+1}[V_{k,n+1}, \rho^{n+1}_{\infty,1}(t) - \rho^{n+1}_{\infty,2}(s)])\|_{\mathcal{L}^1_{n+1}} ds$$

$$\leq 2n \|V\|_{L^\infty_{\mathbb{T}^3}} \int_0^t \|\rho^{n+1}_{\infty,1}(s) - \rho^{n+1}_{\infty,2}(s)\|_{\mathcal{L}^1_{n+1}}. \tag{6.49}$$

By iterating the previous inequality we obtain that

$$\|\rho^n_{\infty,1}(t) - \rho^n_{\infty,2}(t)\|_{\mathcal{L}^1_n} \leq \frac{(n+k-1)!}{(n-1)!} (2\|V\|_{L^\infty_{\mathbb{T}^3}})^k \int_{S_k(t)} \|\rho^n_{\infty,1}(s_k) - \rho^n_{\infty,2}(s_k)\|_{\mathcal{L}^1_{n+k}} ds_{1:k}, \tag{6.50}$$

where the simplex $S_k(t)$ is given by $S_k(t) := \{0 \leq t_1 \leq \cdots \leq t_k \leq t\}$. Recalling that the simplex $S_k(t)$ has size $|S_k(t)| = \frac{t^k}{k!}$, and that

$$\|\rho_{\infty,1}^{n+k}(s_k) - \rho_{\infty,2}^{n+k}(s_k)\|_{\mathcal{L}_{n+k}^1} \leq \|\rho_{\infty,1}^{n+k}(s_k)\|_{\mathcal{L}_{n+k}^1} + \|\rho_{\infty,1}^{n+k}(s_k)\|_{\mathcal{L}_{n+k}^1} = 2,$$

inequality (6.49) gives

$$\|\rho_{\infty,1}^{n}(t) - \rho_{\infty,2}^{n}(t)\|_{\mathcal{L}_n^1} \leq \frac{(n+k-1)!}{(n-1)!k!}(4t\|V\|_{L^\infty_{\mathbb{T}^3}})^k = \binom{n+k-1}{k}(4t\|V\|_{L^\infty_{\mathbb{T}^3}})^k \leq 2^{n-1}(8t\|V\|_{L^\infty_{\mathbb{T}^3}})^k.$$

So if we choose $t < \frac{1}{8\|V\|_{L^\infty_{\mathbb{T}^3}}}$ and then send $k \to \infty$ we get

$$\|\rho_{\infty,1}^{n}(t) - \rho_{\infty,2}^{n}(t)\|_{\mathcal{L}_n^1} = 0, \quad \text{for all } 0 < t < \frac{1}{8\|V\|_{L^\infty_{\mathbb{T}^3}}} \text{ and } n \in \mathbb{N}.$$

We conclude that $\rho_{\infty,1}^n(t) = \rho_{\infty,2}^n(t)$ for all $t \geq 0$ thanks to the flow property of the equation. □

We conclude the paper with the proof of Corollary 6.5.

**Proof of Corollary 6.5.** Let $u_0 \in \mathcal{D}(\sqrt{-\mathcal{H}})$ such that $\|u_0\|_{L^2_{\mathbb{T}^3}} = 1$ then by Theorem 5.15 there is a unique solution $u$ to the NLS equation with Hartree nonlinearity (5.25) with $V \in L^\infty_{\mathbb{T}^3}$ such that $u \in C(\mathbb{R}_+, \mathcal{D}(\sqrt{-\mathcal{H}}))$.

$(\Pi_{(u(t,\cdot))^{\otimes n}})_{n \in \mathbb{N}}$ is easily verified to be a mild solution to the hierachy equation (6.15) essentially due to (6.12) of $\mathcal{W}^{1,1}$-type as $u \in C(\mathbb{R}_+, \mathcal{D}(\sqrt{-\mathcal{H}}))$ and with initial condition

$$\Pi_{u_0^{\otimes n}}(x_{1:n}, y_{1:n}) = u_0^{\otimes n}(x_{1:n})\overline{u_0^{\otimes n}(y_{1:n})},$$

and thus $(\Pi_{(u(t,\cdot))^{\otimes n}})_{n \in \mathbb{N}}$ is also unique by Theorem 6.4.

If $\Psi_{N,\delta}$ is a solution to the Schrödinger equation (6.4), satisfying the hypotheses of Corollary 6.5, then by Theorem 6.3 the corresponding density matrices $\rho_{N,\delta}^n$ –up to a subsequence– must converge to a solution to the hierarchy equation (6.15) if its initial condition converges $\rho_{N,\delta}^n(0) \to \Pi_{u_0^{\otimes n}}$ as $N \to \infty$ and $\delta \to 0^+$.

Since any convergent subsequence converges to the same limit, we get

$$\lim_{N \to +\infty, \delta \to 0^+} \rho_{N,\delta}^n(t, x_{1:n}, y_{1:n}) \to (u(t,\cdot))^{\otimes n}(x_{1:n})\overline{(u(t,\cdot))^{\otimes n}(y_{1:n})}$$

weakly* in $L^\infty(\mathbb{R}_+, \mathcal{L}_n^1)$, and thus weakly* in $L^\infty(\mathbb{R}_+, L_s^2((\mathbb{T}^3)^n) \otimes L_s^2((\mathbb{T}^3)^n))$ for all $n \in \mathbb{N}$. □

# Appendix A   Analytic Setup

We collect some elementary results about paraproducts and Besov spaces and refer the readers for example to [BCD11],[GIP15] and [vZ22] for more details on related topics and especially the application to singular SPDEs.

We first recall the definition of Littlewood-Paley blocks and denoted by $\chi$ and $\rho$ two non-negative smooth and compactly supported radial functions $\mathbb{R}^d \to \mathbb{C}$ such that

1. The support of $\chi$ is contained in a ball and the support of $\rho$ is contained in an annulus $\{x \in \mathbb{R}^d : a \le |x| \le b\}$

2. For all $\xi \in \mathbb{R}^d$, $\chi(\xi) + \sum_{j \ge 0} \rho(2^{-j}\xi) = 1$;

3. For $j \ge 1$, $\chi(\cdot)\rho(2^{-j}\cdot) = 0$ and $\rho(2^{-j}\cdot)\rho(2^{-i}\cdot) = 0$ for $|i-j| > 1$.

The Littlewood-Paley blocks $(\Delta_j)_{j \ge -1}$ associated to $f \in \mathcal{S}'(\mathbb{T}^d)$ are defined by $\Delta_{-1} f := \mathcal{F}^{-1}(\chi \mathcal{F} f)$ and $\Delta_j f := \mathcal{F}^{-1}(\rho(2^{-j}\cdot)\mathcal{F} f)$ for $j \ge 0$. We also set, for $f \in \mathcal{S}'(\mathbb{T}^d)$ and $j \ge -1$ that $S_j f := \sum_{i=-1}^{j-1} \Delta_i f$. Then the Besov space with parameters $p, q \in [1, \infty)$, $\alpha \in \mathbb{R}$ can now be defined as $B_{p,q}^\alpha(\mathbb{T}^d) := \{u \in \mathcal{S}'(\mathbb{T}^d) : \|u\|_{B_{p,q}^\alpha} < \infty\}$, where the norm is defined as

$$\|u\|_{B_{p,q}^\alpha} := \left( \sum_{k \ge -1} \left( (2^{\alpha k} \|\Delta_k u\|_{L^p})^q \right) \right)^{\frac{1}{q}},$$

with the obvious modification for $q = \infty$. We also define the *Besov-Hölder* spaces $C^\alpha := B_{\infty,\infty}^\alpha$, which for $\alpha \in (0,1)$ agree with the usual Hölder spaces $C^\alpha$. Using this notation, we can formally decompose the product $f \cdot g$ of two distributions $f$ and $g$ as

$$f \cdot g = f \olessthan g + f \odot g + f \ogreaterthan g,$$

where

$$f \olessthan g := \sum_{j \ge -1} S_{j-1} f \Delta_j g \quad \text{and} \quad f \ogreaterthan g := \sum_{j \ge -1} \Delta_j f S_{j-1} g$$

are referred to as the *paraproducts*, whereas

$$f \odot g := \sum_{j \ge -1} \sum_{|i-j| \le 1} \Delta_i f \Delta_j g$$

is called the *resonant product*. An important point is that the paraproduct terms are always well defined whatever the regularity of $f$ and $g$. The resonant product, on the other hand, is a priori only well defined if the sum of their regularities is positive. We collect some results.

**Lemma A.1.** ( [MW17a], Theorem 3.17) *Let $\alpha, \alpha_1, \alpha_2 \in \mathbb{R}$ and $p, p_1, p_2, q \in [1, \infty]$ be such that*

$$\alpha_1 \ne 0 \quad \alpha = (\alpha_1 \wedge 0) + \alpha_2 \quad \text{and} \quad \frac{1}{p} = \frac{1}{p_1} + \frac{1}{p_2}.$$

*Then we have the bound*

$$\|f \olessthan g\|_{B_{p,q}^\alpha} \lesssim \|f\|_{B_{p_1,\infty}^{\alpha_1}} \|g\|_{B_{p_2,q}^{\alpha_2}}$$

*and in the case where $\alpha_1 + \alpha_2 > 0$ we have the bound*

$$\|f \odot g\|_{B_{p,q}^{\alpha_1+\alpha_2}} \lesssim \|f\|_{B_{p_1,\infty}^{\alpha_1}} \|g\|_{B_{p_2,q}^{\alpha_2}}.$$

**Remark A.2.** In this paper, we only consider the cases $W^{s,p} = B_{p,p}^s$ *fractional Sobolev*, $H^s = B_{2,2}^s$ and *Hölder Besov* spaces $C^s = B_{\infty,\infty}^s$ using the convention $W^{s,2} = H^s$ and $W^{s,\infty} = C^s$.

**Lemma A.3. (Bernstein Inequality)** *Let $\mathcal{A}$ be an annulus and $\mathcal{B}$ be a ball. For any $k \in \mathbb{N}, \lambda > 0,$ and $1 \leq p \leq q \leq \infty$ we have*

1. *if $u \in L^p(\mathbb{T}^d)$ is such that $\mathrm{supp}(\mathcal{F}u) \subset \lambda \mathcal{B}$ then*

$$\max_{\mu \in \mathbb{N}^d : |\mu| = k} \|\partial^\mu u\|_{L^q} \lesssim_k \lambda^{k + d\left(\frac{1}{p} - \frac{1}{q}\right)} \|u\|_{L^p}$$

2. *if $u \in L^p(\mathbb{T}^d)$ is such that $\mathrm{supp}(\mathcal{F}u) \subset \lambda \mathcal{A}$ then*

$$\lambda^k \|u\|_{L^p} \lesssim_k \max_{\mu \in \mathbb{N}^d : |\mu| = k} \|\partial^\mu u\|_{L^p}.$$

**Lemma A.4. (Besov Embedding)** *Let $\alpha < \beta \in \mathbb{R}$, $q_1 \leq q_2$, and $p > r \in [1, \infty]$ be such that $\beta = \alpha + d(1/r - 1/p)$, then we have the following bound*

$$\|f\|_{B^\alpha_{p,q_2}(\mathbb{T}^d)} \lesssim \|f\|_{B^\beta_{r,q_1}(\mathbb{T}^d)}.$$

*Moreover, we have the following embedding between fractional Sobolev spaces and Bessel potential spaces*

$$\|f\|_{W^{\alpha,p}_{\mathbb{T}^d}} \lesssim \left\|(1-\Delta)^{\frac{\beta}{2}} f\right\|_{L^p_{\mathbb{T}^d}} \lesssim \|f\|_{W^{\beta,p}_{\mathbb{T}^d}} \quad \text{for } \alpha < \beta \quad \text{and } p \geq 2. \tag{A.1}$$

**Proposition A.5. (Commutator in [GIP15], see Proposition A.9 [MW17b])**

*Given $\alpha \in (0,1)$, $\beta, \gamma \in \mathbb{R}$ such that $\beta + \gamma < 0$ and $\alpha + \beta + \gamma > 0$, the following trilinear operator $C$ defined for any smooth functions $f, g, h$ by*

$$\mathrm{com}(f, g, h) := (f \varobslash g) \varodot h - f \cdot (g \varodot h)$$

*can be extended continuously to the product space $W^{\alpha,p} \times C^\beta \times C^\gamma$. Moreover, we have the following bound*

$$\|\mathrm{com}(f, g, h)\|_{W^{\alpha+\beta+\gamma-\delta,p}} \lesssim \|f\|_{W^{\alpha,p}} \|g\|_{C^\beta} \|h\|_{C^\gamma}$$

*for all $f \in W^{\alpha,p}$, $g \in C^\beta$ and $h \in C^\gamma$, and every $\delta > 0, p \geq 1$.*

**Lemma A.6. (Fractional Leibniz [GK96])** *Let $1 < p < \infty$ and $p_1, p_2, p'_1, p'_2$ such that*

$$\frac{1}{p_1} + \frac{1}{p_2} = \frac{1}{p'_1} + \frac{1}{p'_2} = \frac{1}{p}.$$

*Then for any $s, \alpha \geq 0$ there exists a constant s.t.*

$$\|\langle \nabla \rangle^s (fg)\|_{L^p} \leq C \|\langle \nabla \rangle^{s+\alpha} f\|_{L^{p_2}} \|\langle \nabla \rangle^{-\alpha} g\|_{L^{p_1}} + C \|\langle \nabla \rangle^{-\alpha} f\|_{L^{p'_2}} \|\langle \nabla \rangle^{s+\alpha} g\|_{L^{p'_1}}.$$

# Appendix B  Anderson Hamiltonian

## B.1  Anderson Hamiltonian: Technical bounds

The remainder terms in Proposition 2.3 are given by

$$-(1-\Delta)B(u^\flat) = 2\nabla W^M \varobslash \nabla u^\flat + [\nabla u^\flat \varobslash, 1-\Delta](1-\Delta)^{-1} \nabla W^M + 2u^\flat \varobslash R^2 + 4\nabla u^\flat \varobslash R^1 \tag{B.1}$$

and

$$\begin{aligned}
G(u^\flat) &= u^\flat \succcurlyeq [(1-\Delta)Z] - \Delta P_{\leq N}(u^\flat \varolessthan Z + 2\nabla u^\flat \varolessthan (1-\Delta)^{-1}\nabla W^M + B(u^\flat)) \\
&+ u^\flat \varolessthan Z + 2\nabla u^\flat \varolessthan (1-\Delta)^{-1}\nabla W^M + B(u^\flat) \\
&+ u^\flat \succcurlyeq R^2 + \nabla u^\flat \succcurlyeq R^1 + 2[\nabla u^\flat \varolessthan Z + 2\nabla^2 u^\flat \varolessthan (1-\Delta)^{-1}\nabla W^M] \varolessthan \nabla W^M \\
&+ \text{com}(\nabla u^\flat, (1-\Delta)^{-1}\nabla^2 W^M, \nabla W^M) + \text{com}(u^\flat, \nabla Z, \nabla W^M) \\
&- \nabla P_{\leq N}(u^\flat \varolessthan Z + 2\nabla u^\flat \varolessthan (1-\Delta)^{-1}\nabla W^M + B(u^\flat)) \varolessthan \nabla W^M \\
&+ \nabla B(u^\flat) \varolessthan \nabla W^M + [u^\flat \varolessthan, 1-\Delta](1-\Delta)^{-1}Z.
\end{aligned} \tag{B.2}$$

The continuity estimates of the functions above are indicated by those of paraproducts and resonant products and the regularity of the stochastic objects provided by Lemma 2.1 immediately.

**Lemma B.1.** *For any $\varepsilon > 0$ small enough, $2 \leq p \leq \infty$, $\gamma < 2$ and $\gamma' \in \left[\frac{3}{2} + \varepsilon, 2\right)$, we have*

$$\|B\|_{\mathcal{L}\left(W^{\gamma-1/2+\varepsilon,p}_{\mathbb{T}^3}, W^{\gamma,p}_{\mathbb{T}^3}\right)} \lesssim \|\Xi\|_{\mathcal{X}} \tag{B.3}$$

*and*

$$\|G\|_{\mathcal{L}\left(W^{\gamma',p}_{\mathbb{T}^3}, W^{\gamma'-3/2-\varepsilon,p}_{\mathbb{T}^3}\right)} \lesssim \|\Xi\|_{\mathcal{X}}. \tag{B.4}$$

*In particular we have $\|G(u^\flat)\|_{L^2_{\mathbb{T}^3}} \lesssim_\Xi \|u^\flat\|_{H^{\frac{3}{2}+\varepsilon}_{\mathbb{T}^3}}$.*

**Proof.** For the first bound, it suffices to see in (B.1) that

$$[\nabla u^\flat \varolessthan, 1-\Delta](1-\Delta)^{-1}\nabla W^M = (\Delta \nabla u^\flat) \varolessthan (1-\Delta)^{-1}\nabla W^M + 2\nabla^2 u^\flat \varolessthan (1-\Delta)^{-1}\nabla^2 W^M,$$

while for the second one, we have

$$\begin{aligned}
\|\text{com}(\nabla u^\flat, (1-\Delta)^{-1}\nabla^2 W^M, \nabla W^M)\|_{W^{\alpha-1-\frac{\varepsilon}{2},p}_{\mathbb{T}^3}} &\lesssim \|u^\flat\|_{W^{\alpha,p}_{\mathbb{T}^3}}, \quad \text{for } \alpha \in (1+\varepsilon, 2), \\
\|\text{com}(u^\flat, \nabla Z, \nabla W^M)\|_{W^{\alpha-\frac{\varepsilon}{2},p}_{\mathbb{T}^3}} &\lesssim \|u^\flat\|_{W^{\alpha,p}_{\mathbb{T}^3}}, \quad \text{for } \alpha \in (\varepsilon, 1), \\
\|\nabla u^\flat \succcurlyeq R^1\|_{W^{\alpha-\varepsilon,p}_{\mathbb{T}^3}} &\lesssim \|u^\flat\|_{W^{\alpha+1,p}_{\mathbb{T}^3}}, \quad \text{for } \alpha \in (\varepsilon, \infty), \\
\|\nabla B(u^\flat) \varolessthan \nabla W^M\|_{W^{\frac{1}{2}-\varepsilon,p}_{\mathbb{T}^3}} &\lesssim \|u^\flat\|_{W^{\alpha,p}_{\mathbb{T}^3}}, \quad \text{for } \alpha \in \left[\frac{3}{2}+\varepsilon, 2\right),
\end{aligned}$$

while the other terms are more regular. □

The properties of the ansatz are summarized as follows. The exponential ansatz (i.e. multiplication with $e^W$) and the paracontrolled ansatz (i.e. the $\Gamma$ map) transform the domain $\mathcal{D}(\mathcal{H})$ to $H^2_{\mathbb{T}^3}$. Here we recall the construction of the map $\Gamma$ and its continuity properties.

**Proof.** (Lemma 2.5) The first claim follows from Proposition 2.57 in [GUZ20], while the continuity of the $\Gamma$-map can be proven with similar argument as in Proposition 2.56 in [GUZ20]. The continuity of $\Gamma^{-1}$ can be directly deduced from the definition, since

$$\Gamma^{-1} u^\flat = u^\flat - P_{>N}[u^\flat \varolessthan Z + \nabla u^\flat \varolessthan (1-\Delta)^{-1}\nabla W^M + B(u^\flat)],$$

for any $u^\flat \in e^{-W}\mathcal{D}(\mathcal{H})$ or $W^{\beta,p}_{\mathbb{T}^3}$ as in the statement, where

$$Z, (1-\Delta)^{-1}\nabla W^M \in \mathcal{C}^{\frac{3}{2}-}_{\mathbb{T}^3}$$

and $B$ is continuous from $W_{\mathbb{T}^3}^{\beta+,p}$ to $W_{\mathbb{T}^3}^{\beta+\frac{1}{2},p}$ for any $2 \leq p \leq \infty$ and $\beta < \frac{3}{2}$, see (B.3), which is a slight modification of Lemma 2.41 in [GUZ20]. Lastly, since

$$(\Gamma - \mathrm{id})u^\sharp = P_{>N}[\Gamma u^\sharp \otimes Z + \nabla \Gamma u^\sharp \otimes (1-\Delta)^{-1}\nabla W^M + B(\Gamma u^\sharp)],$$

the bound is a direct consequence of (B.3), the continuity of the $\Gamma$-map, the regularity estimates of the stochastic objects given in Lemma 2.1 and the condition that $0 < \gamma \leq 1$. □

For the Strichartz estimate Lemma 3.4 for the full regularity regime as claimed, it is essential to have a characterization of Bessel potential spaces using $\mathcal{H}^\sharp$, extending Proposition 2.57 in [GUZ20] which says that

$$\|\mathcal{H}^\sharp u^\sharp\|_{L^2_{\mathbb{T}^3}} \eqsim \|u^\sharp\|_{H^2_{\mathbb{T}^3}}. \tag{B.5}$$

Since we will only use the following bounds for fixed $p$, the dependence of the coefficients on $p$ is omitted.

**Lemma B.2.** *The boundary terms defined in (5.22) is bounded in $L^{\frac{10}{3}}_{[0,T]\times\mathbb{T}^3}$ by*

$$\|\mathrm{Bd}[u^\sharp]\|_{L^{\frac{10}{3}}_{[0,T]\times\mathbb{T}^3}} \lesssim \|u^\sharp\|^3_{L^\infty_{[0,T]}L^6_{\mathbb{T}^3}}\|V_\beta\|_{L^{\frac{20}{19}}_{\mathbb{T}^3}}.$$

**Proof.** The boundary terms can be bounded by

$$\|\mathrm{Bd}[u^\sharp]\|_{L^{\frac{10}{3}}_{[0,T]\times\mathbb{T}^3}} \lesssim \|\Lambda^{-1}[\Lambda u^\sharp(t)(|\Lambda u^\sharp(t)|^2 * V_\beta)]\|_{L^{\frac{10}{3}}_{[0,T]}W^{-\frac{9}{20}-,\frac{10}{3}}_{\mathbb{T}^3}}$$
$$+ \left\|(-\mathcal{H}^\sharp)^{-\frac{9}{40}-}\Lambda^{-1}[\Lambda u^\sharp(0)(|\Lambda u^\sharp(0)|^2 * V_\beta)]\right\|_{H^{\frac{9}{20}+}_{\mathbb{T}^3}},$$

where the first can be bounded by

$$\|\Lambda^{-1}[\Lambda u^\sharp(t)(|\Lambda u^\sharp(t)|^2 * V_\beta)]\|_{L^{\frac{10}{3}}_{[0,T]}W^{-\frac{9}{20}-,\frac{10}{3}}_{\mathbb{T}^3}} \lesssim \|\Lambda u^\sharp(t)(|\Lambda u^\sharp(t)|^2 * V_\beta)\|_{L^{\frac{10}{3}}_{[0,T]}L^{\frac{20}{9}}_{\mathbb{T}^3}}$$
$$\lesssim \|u^\sharp(t)\|^3_{L^\infty_{[0,T]}L^6_{\mathbb{T}^3}}\|V_\beta\|_{L^{\frac{20}{19}}_{\mathbb{T}^3}},$$

using Lemma 2.5, while

$$\left\|(-\mathcal{H}^\sharp)^{-\frac{9}{40}-}\Lambda^{-1}[\Lambda u^\sharp(0)(|\Lambda u^\sharp(0)|^2 * V_\beta)]\right\|_{H^{\frac{9}{20}+}_{\mathbb{T}^3}} \lesssim \|\Lambda u^\sharp(0)(|\Lambda u^\sharp(0)|^2 * V_\beta)\|_{L^2_{\mathbb{T}^3}}$$
$$\lesssim \|u^\sharp(0)\|^3_{L^6_{\mathbb{T}^3}}\|V_\beta\|_{L^1_{\mathbb{T}^3}}$$

having used Lemma 3.5. □

We now prove the heat kernel estimate of $(e^{t\mathcal{H}^\sharp})_{t\geq 0}$ in Lemma 2.12.

**Proof. (Lemma 2.12)** We start by showing (2.30) By functional calculus and using the boundedness of $\Lambda$, one can see that for any $t > 0$

$$\|e^{t\mathcal{H}^\sharp}u\|_{H^2_{\mathbb{T}^3}} \eqsim \|e^{t\mathcal{H}}\mathcal{H}u\|_{L^2_{\mathbb{T}^3}} \lesssim \sup_{\lambda > 0}(\lambda e^{-t\lambda})\|u\|_{L^2_{\mathbb{T}^3}} \lesssim t^{-1}\|u\|_{L^2_{\mathbb{T}^3}}$$

and $\|e^{t\mathcal{H}^\sharp}u\|_{L^2} \lesssim \|u\|_{L^2}$ so by interpolation

$$\|e^{t\mathcal{H}^\sharp}u\|_{H^\gamma_{\mathbb{T}^3}} \lesssim t^{-\frac{\gamma}{2}}\|u\|_{L^2_{\mathbb{T}^3}} \quad 0 < \gamma < 2.$$

Next we want to prove (2.31) and write

$$e^{t\mathcal{H}^{\natural}}u = e^{t\Delta}u + \int_0^t e^{(t-s)\Delta}(\mathcal{H}^{\natural} - \Delta)e^{s\mathcal{H}^{\natural}}u\,ds.$$

Then for $p \geq 2$ and $\sigma \in \left(\frac{3}{2}, 2\right]$ and for $\varepsilon, \varepsilon' > 0$ small enough and $0 \leq t \leq 1$

$$\|e^{t\mathcal{H}^{\natural}}u\|_{W^{2,p}} \lesssim \|e^{t\Delta}u\|_{W^{2,p}} + \int_0^t |t-s|^{-\frac{3}{4}+\varepsilon}\|(\mathcal{H}^{\natural}-\Delta)e^{s\mathcal{H}^{\natural}}u\|_{W^{\frac{1}{2}-\varepsilon,p}}$$

$$\lesssim \|e^{t\Delta}u\|_{W^{2,p}} + \int_0^t |t-s|^{-\frac{3}{4}+\varepsilon}\|e^{s\mathcal{H}^{\natural}}u\|_{W^{2,p}}ds$$

$$\lesssim \|e^{t\Delta}u\|_{W^{2,p}} + C(\varepsilon)\left(\int_0^t \|e^{s\mathcal{H}^{\natural}}u\|_{W^{\sigma,p}}^{4+\varepsilon'}ds\right)^{\frac{1}{4+\varepsilon'}}$$

$$\Rightarrow$$

$$\|e^{t\mathcal{H}^{\natural}}u\|_{W^{2,p}}^{4+\varepsilon'} \lesssim \|e^{t\Delta}u\|_{W^{2,p}}^{4+\varepsilon'} + C(\varepsilon,\varepsilon')\int_0^t \|e^{s\mathcal{H}^{\natural}}u\|_{W^{\sigma,p}}^{4+\varepsilon'}ds$$

Gronwall implies

$$\|e^{t\mathcal{H}^{\natural}}u\|_{W^{2,p}}^{4+\varepsilon'} \lesssim \|e^{t\Delta}u\|_{W^{2,p}}^{4+\varepsilon'} + C(\varepsilon,\varepsilon')\int_0^t \|e^{s\Delta}u\|_{W^{\sigma,p}}^{4+\varepsilon'}ds$$

$$\lesssim t^{-\frac{\gamma}{2}(4+\varepsilon')}\|u\|_{W^{2-\gamma,p}}^{4+\varepsilon'} + C(\varepsilon,\varepsilon')\|u\|_{W^{2-\gamma,p}}^{4+\varepsilon'}\int_0^t s^{-\frac{\gamma}{2}(4+\varepsilon')}ds$$

$$\lesssim C(\varepsilon,\varepsilon')t^{-\frac{\gamma}{2}(4+\varepsilon')}\|u\|_{W^{2-\gamma,p}}^{4+\varepsilon'}$$

$$\Rightarrow$$

$$\|e^{t\mathcal{H}^{\natural}}u\|_{W^{2,p}} \lesssim t^{-\frac{\gamma}{2}}\|u\|_{W^{2-\gamma,p}}$$

for $\frac{\gamma}{2}(4+\varepsilon') < 1$ and $0 < t < 1$. So the flow regularises at high regularities by almost $\frac{1}{2}$.

By an analogous proof (or setting $\gamma = 0$) we can show for $\sigma \in \left(\frac{3}{2}, 2\right)$ and $p \geq 2$

$$\|e^{t\mathcal{H}^{\natural}}u\|_{W^{\sigma,p}} \lesssim \|e^{t\Delta}u\|_{W^{\sigma,p}} + \int_0^t |t-s|^{-\frac{3}{4}+\varepsilon}\|(\mathcal{H}^{\natural}-\Delta)e^{s\mathcal{H}^{\natural}}u\|_{W^{\sigma-\frac{3}{2}-\varepsilon,p}}$$

$$\lesssim \|u\|_{W^{\sigma,p}} + \int_0^t |t-s|^{-\frac{3}{4}+\varepsilon}\|e^{s\mathcal{H}^{\natural}}u\|_{W^{\sigma,p}}$$

$$\Rightarrow$$

$$\|e^{t\mathcal{H}^{\natural}}u\|_{W^{\sigma,p}} \lesssim \|u\|_{W^{\sigma,p}} \text{ for } 0 \leq t \leq 1$$

which is (2.32).

Combining these results, we can prove that for fixed time, the semigroup regularises up to $W^{2,\infty}$. By iteratively applying Besov embedding, we can bound for $K = 4\lceil t \rceil + 1$

$$\|e^{t\mathcal{H}^{\natural}}u\|_{W^{2,\infty}_{\mathbb{T}^3}} \lesssim \left(\frac{t}{K}\right)^{-\frac{1}{4}+}\|e^{\frac{K-1}{K}t\mathcal{H}^{\natural}}u\|_{W^{\frac{3}{2}+,\infty}_{\mathbb{T}^3}} \lesssim \left(\frac{t}{K}\right)^{-\frac{1}{4}+}\|e^{\frac{K-1}{K}t\mathcal{H}^{\natural}}u\|_{W^{2,6+}_{\mathbb{T}^3}}$$

$$\lesssim \left(\frac{t}{K}\right)^{-\frac{1}{2}+}\|e^{\frac{K-2}{K}t\mathcal{H}^{\natural}}u\|_{W^{\frac{3}{2}+,6+}_{\mathbb{T}^3}} \lesssim \left(\frac{t}{K}\right)^{-\frac{1}{2}+}\|e^{\frac{K-2}{K}t\mathcal{H}^{\natural}}u\|_{W^{2,3+}_{\mathbb{T}^3}}$$

$$\lesssim \left(\frac{t}{K}\right)^{-\frac{3}{4}+}\|e^{\frac{K-3}{K}t\mathcal{H}^{\natural}}u\|_{W^{\frac{3}{2}+,3+}_{\mathbb{T}^3}} \lesssim \left(\frac{t}{K}\right)^{-\frac{3}{4}+}\|e^{\frac{K-3}{K}t\mathcal{H}^{\natural}}u\|_{W^{2,2+}_{\mathbb{T}^3}}$$

$$\lesssim \left(\frac{t}{K}\right)^{-1+}\|e^{\frac{K-4}{K}t\mathcal{H}^{\natural}}u\|_{H^{\frac{3}{2}+}_{\mathbb{T}^3}} \lesssim K^{1+}\left(\frac{K}{K-4}\right)^{\frac{3}{2}+}t^{-\frac{7}{4}-}\|u\|_{L^2_{\mathbb{T}^3}}$$

$$\lesssim \lceil t \rceil^{1+}t^{-\frac{7}{4}-}\|u\|_{L^2_{\mathbb{T}^3}}$$

which for short times is only slightly worse than what one gets for the usual heat semigroup. □

# Appendix C  Topology of Density Matrices

**Lemma C.1.** *We fix some Hilbert space H.*

a) *Let $(B_n)_{n \geq 0}$ be a sequence of operators converging strongly to $B$ and $O$ be a compact operator on some Hilbert space, then $S_n O \to SO$ and $OS_n \to OS$ both under operator norm.*

b) *If in addition there are bounded operators $(A_n)_{n \geq 0}$ on $H$ which verify that $A_n \to A$ in $(\mathcal{L}^1(H), w^*)$, we have $A_n B_n \to AB$ in $(\mathcal{L}^1(H), w^*)$.*

c) *Lastly, if we have further some sequence $(C_n)_{n \geq 0} \subset \mathcal{K}(H)$ that converges under operator norm $C_n \to C \in \mathcal{K}(H)$, then $\mathrm{Tr}(C_n A_n) \to \mathrm{Tr}(CA)$.*

**Proof.**

a) Since $B_n$ converge strongly then $\sup \|B_n\|_{\mathcal{L}} < +\infty$. Let $O_\varepsilon$ be a finite-rank operator such that $\|O_\varepsilon - O\|_{\mathcal{L}} < \varepsilon$ then we have

$$\limsup_{n \to \infty} \|S_n O - SO\|_{\mathcal{L}} \leq \lim_{n \to \infty} \|S_n O_\varepsilon - SO_\varepsilon\|_{\mathcal{L}} + (\sup \|S_n\|_{\mathcal{L}} + \|S\|_{\mathcal{L}}) \varepsilon \lesssim \varepsilon$$

since $\lim_{n \to \infty} \|S_n O_\varepsilon - SO_\varepsilon\|_{\mathcal{L}} = 0$ due to the fact that $O_\varepsilon$ has finite rank and $(S_n)$ converges strongly. Since $\varepsilon > 0$ is arbitrary, we get the convergence of $S_n O$. For the case $OS_n$, we only need to apply the first result to the adjoint $(OS_n)^* = S_n^* O^*$.

b) For any compact operator $K \in \mathcal{K}(H)$, we have

$$|\mathrm{Tr}((A_n B_n - AB) K)| \leq |\mathrm{Tr}((A_n - A) BK)| + |\mathrm{Tr}(A_n (B - B_n) K)|,$$

where the first term vanishes as $n \to \infty$ since $\mathcal{K}(H)$ is a two-sided ideal of bounded operators on $H$; while for the second term, we have $\sup_{n \geq 1} \|A_n\|_{\mathcal{L}^1} < \infty$ due to Banach-Steinhaus theorem, and $\sup_{n \geq 1} \|B - B_n\| < \infty$ by strong convergence, so that we can assume without loss of generality that $K$ is finite-rank, since we can always approximate compact $K$ by finite-rank $K_\varepsilon$ under operator norm, so that

$$\sup_{n \geq 1} |\mathrm{Tr}(A_n (B - B_n)(K - K_\varepsilon))| \leq \sup_{n \geq 1} \left[ \mathrm{Tr}(|A_n|) \|B - B_n\|_{\mathcal{L}} \right] \|K - K_\varepsilon\|_{\mathcal{L}} \to 0,$$

as $\varepsilon \to 0^+$. In the case that $K$ is finite-rank, $\|(B-B_n)K\| \to 0$ as $n \to \infty$ by strong convergence of $(B_n)$, and the claim follows.

c) Lastly, we have

$$|\text{Tr}(C_n A_n) - \text{Tr}(CA)| \leq |\text{Tr}(C(A_n - A))| + |\text{Tr}((C_n - C)A_n)|,$$

where $|\text{Tr}(C(A_n - A))| \to 0$ since $C$ is compact and $|\text{Tr}((C_n - C)A_n)| \leq 2\sup_n(\|A_n\|_{\mathcal{L}^1})\|C_n - C\|_{\mathcal{L}} \to 0$, as $n \to \infty$. $\square$

**Proof. (Lemma 6.15)** The group structure is clear. Fix $\rho \in \mathcal{L}_n^1$, then we have by the polar decomposition that

$$\mathcal{U}_t^n(\rho) = \exp\left(it\sum_{l=1}^n \mathcal{H}_{x_l}\right)\rho\exp\left(-it\sum_{l=1}^n \mathcal{H}_{x_l}\right) = \exp\left(it\sum_{l=1}^n \mathcal{H}_{x_l}\right)U|\rho|\exp\left(-it\sum_{l=1}^n \mathcal{H}_{x_l}\right),$$

for some unitary $U$, which implies

$$|\mathcal{U}_t^n(\rho)|^2 = \exp\left(it\sum_{l=1}^n \mathcal{H}_{x_l}\right)U|\rho|^2 U^*\exp\left(-it\sum_{l=1}^n \mathcal{H}_{x_l}\right)$$

$$\Rightarrow$$

$$|\mathcal{U}_t^n(\rho)| = \exp\left(it\sum_{l=1}^n \mathcal{H}_{x_l}\right)U|\rho|U^*\exp\left(-it\sum_{l=1}^n \mathcal{H}_{x_l}\right).$$

Therefore,

$$\|\mathcal{U}_t^n(\rho)\|_{\mathcal{L}_n^1} = \text{Tr}\left(\exp\left(it\sum_{l=1}^n \mathcal{H}_{x_l}\right)U|\rho|U^*\exp\left(-it\sum_{l=1}^n \mathcal{H}_{x_l}\right)\right)$$

$$= \sum_\alpha \left(|\rho|U^*\exp\left(-it\sum_{l=1}^n \mathcal{H}_{x_l}\right)f_\alpha, U^*\exp\left(-it\sum_{l=1}^n \mathcal{H}_{x_l}\right)f_\alpha\right)$$

$$= \sum_{\alpha \in \mathbb{N}^n} \exp\left(-i\sum_{l=1}^n \lambda_{\alpha_l}\right)\exp\left(-i\sum_{l=1}^n \lambda_{\alpha_l}\right)(|\rho|U^* f_\alpha, U^* f_\alpha)$$

$$= \sum_{\alpha \in \mathbb{N}^n} (|\rho|U^* f_\alpha, U^* f_\alpha) = \text{Tr}(|\rho|) = \|\rho\|_{\mathcal{L}_n^1},$$

where $f_\alpha := \bigotimes_{j=1}^n f_{\alpha_j}$ and $(f_i)_{i \in \mathbb{N}}$ is the orthonormal basis of $L^2(\mathbb{T}^3)$ of eigenvectors of $\mathcal{H}$ and $(\lambda_i)_{i \in \mathbb{N}}$ are the corresponding eigenvalues, since $(U^* f_\alpha)_\alpha$ is also a family of orthonormal basis. $\square$

**Acknowledgments.**

F.C.D. is partially funded by INdAM (Istituto Nazionale di Alta Matematica, Gruppo Nazionale per l'Analisi Matematica, la Probabilità e le loro Applicazioni) through the project "Analisi spettrale, armonica e stocastica in presenza di potenziali magnetici".

X.J. was funded by the DRS-Stipendium through Nicolas Perkowski from FU Berlin, and part of this work was carried out when X.J. was affiliated with the Max-Planck Institute for Mathematics in the Sciences in Leipzig, where the hospitality from the Research Group "Pattern Formation, Energy Landscapes and Scaling Laws" is sincerely acknowledged.

I.Z. acknowledges support from DFG CRC/TRR 388 "Rough Analysis, Stochastic Dynamics and Related Fields" Project A01. He also thanks Nicolas Rougerie for some helpful comments about mean field Quantum dynamics.